%%%%%%%%%%%%%%%%%%%%%%%%%%%%%%%%%%%%%%%%%%%%%%%%%%%%%%%%%%%%%%%%%%%%%%%%%%%
%% @texfile{
%%     filename="review1a.tex",
%%     version="2.1b",
%%     date="24-DEC-1991",
%%     filetype="AMS-TeX: user documentation",
%%     author="Sergiu I. Vacaru",
%%     address="Institute of Applied Physics,
%%              Academy of Sciences,
%%              5 Academy str., Chisinau 2028,
%%              Republic of Moldova ",
%%     telephone="011-3732-729294",
%%           fax="011-3732-738149",
%%     email="Internet: lises@cc.acad.md  ",
%%     keywords="spinors, generalized Finsler geometry, nonlinear connections,
%%     nearly autoparallel maps",
%%     abstract="This is the file of an article for Memoirs of
%%     the American Mathematical Society:
%%       "Spinors, Nonlinear Connections, and Nearly Autoparallel Maps
%%       of Generalized Finsler Spaces",  by  Sergiu I. Vacaru
%%
% AMSTEX.TEX (version 2.1)
% AMSPPT.STY (version 2.1)
% AMSSYM.TEX (loaded by AMSPPT.STY)
% MSAM10
% MSBM10
% EUFM10
% CMEX8
% CMEX7
% CMBSY7
% CMCSC8
%
%%%%%%%%%%%%%%%%%%%%%%%%%%%%%%%%%%%%%%%%%%%%%%%%%%%%%%%%%%%%%%%%%%%%%%%%

\input amstex
\documentstyle{amsppt}
\Monograph
\pagewidth{30pc}
\pageheight{48.0pc}
\refstyle{B}
\loadbold
\topmatter
\title  Spinors, Nonlinear Connections,\\ and Nearly Autoparallel Maps
of Generalized Finsler Spaces  \endtitle
\rightheadtext{Spinors in Generalized Finsler Spaces}
\leftheadtext{Sergiu I. Vacaru}
\author Sergiu I. Vacaru \endauthor
\affil Department of Statistical and Nuclear Physics, \\
       Institute of Applied Physics,\\
       Academy of Sciences of Moldova, \\
  5 Academy str., Chi\c sin\v au 2028, \\
	  Republic of Moldova \\
 E--mail address:\  lises\@cc.acad.md  \endaffil

\abstract We study the geometric setting of the field theory with locally
anisot\-rop\-ic interactions. The concept of locally anisotropic space
is introduced as a general one for various type of extensions of Lagrange
and Finsler geometry and  higher dimension (Kaluza--Klein type) spaces.
The problem of definition of  spinors  on generalized
 Finsler spaces is solved in the framework of the geometry of Clifford
 bundles provided with compatible nonlinear and distinguished connections
 and metric structures. We construct the spinor differential geometry of
locally anisotropic spaces and discuss some related issues connected
 with the geometric aspects of locally anisotropc interactions for
 gravitational, gauge, spinor, Dirac spinor and Proca fields. Motion
 equations in generalized Finsler spaces, of the mentioned type of fields,
 are defined in a geometric manner by using bundles of linear and affine frames
 locally adapted to the nonlinear connection structure. The nearly
 autoparallel maps are introduced as maps  with deformation of connections
  extending the class of geodesic and conformal transforms. Using this
 approach we propose two variants of solution of the  problem of definition of
  conservation laws on locally anisotropic spaces.
\vskip15pt
{\sl 1991 Mathematics Subject Classification:}\ 53A50, 53B15, 53B40, 53B50,
 53C05, 70D10, 81D27, 81E13, 81E20, 83C40, 83C47, 83D05, 83E15
\vskip10pt
{\sl Key words and phrases:}\ Spinors, nonlinear connections, generalized
 Lagrange and Finsler spaces, locally anisotropic fields, nearly autoparallel
  maps, conservation laws.
 \endabstract
\endtopmatter
\vskip80pt
\newpage
\topmatter  \title Contents \endtitle \endtopmatter
%\document
\toc
\title PREFACE \page{4} \endtitle
\vskip5pt
\title \chapter{1} LOCALLY ANISOTROPIC \newline
  {\qquad}{\qquad} SPINOR SPACES \page{6} \endtitle
\vskip5pt
\head {} {\bf I.1 Connections and Metric in Locally Anisotropic Spaces}
\page{6} \endhead
\subhead {} I.1.1\ Nonlinear connections and
 distinghuishing of geometric objects
\page{6} \endsubhead
\subhead{} I.1.2\ Torsions and curvatures of nonlinear and distinguished
 connections \page{12} \endsubhead
\subhead{} I.1.3\ Field equations for locally anisotropic gravity
 \page{15} \endsubhead
\vskip5pt
\head{} {\bf I.2\ Distinguished Clifford Algebras} \page{16} \endhead
\vskip5pt
\head{} {\bf I.3\ Clifford Bundles and Distinguished Spinor Structures}
     \page{19}   \endhead
\subhead{} I.3.1\ Clifford distinguished module structures in vector bundles
        \page{20} \endsubhead
\subhead{} I.3.2\ Clifford fibration \page{21} \endsubhead
\vskip5pt
\head{} {\bf I.4\ Almost Complex Spinor Structures} \page{22} \endhead
\vskip5pt
\head{} {\bf I.5\ Spinor Techniques for Distinguished Vector Spaces} \page{25}
       \endhead
\subhead{} I.5.1\ Clifford d--algebra, d--spinors and d--twistors \page{25}
     \endsubhead
\subhead{} I.5.2\ Mutual transforms of d--tensors and d--spinors \page{30}
      \endsubhead
\subsubhead{} I.5.2.1\ Transformation of d--tensors into d--spinors \page{30}
   \endsubsubhead
\subsubhead{} I.5.2.2\ Transformation of d--spinors into d--tensors; fundamental
 d--spinors \page{31} \endsubsubhead
\vskip5pt
\head{} {\bf I.6\ The Differential Geometry of Locally Anisotropic
  \newline Spinors}
 \page{32} \endhead
\subhead{} I.6.1\ D--covariant derivation on la--spaces \page{33} \endsubhead
\subhead{} I.6.2\ Infeld---van der Waerden coefficients and d--connections
 \page{34} \endsubhead
\subhead{} I.6.3\ D--spinors of la--space curvature and torsion \page{36}
      \endsubhead
\vskip5pt
\head{} {\bf I.7\ Field Equations on Locally Anisotropic Spaces}
    \page{38}\endhead
\subhead{} I.7.1\ Locally anisotropic scalar field equations \page{38}
 \endsubhead
\subhead{} I.7.2\ Proca equations on la--spaces \page{39}\endsubhead
\subhead{} I.7.3\ La--gravitons on la--backgrounds \page{40}\endsubhead
\subhead{} I.7.4\ Locally anisotropic Dirac equations \page{40}\endsubhead
\subhead{} I.7.5\ D--spinor Yang--Mills equations \page{41}\endsubhead
\vskip5pt
\title\chapter{2} GAUGE FIELDS AND
  \newline {\qquad}{\qquad} LOCALLY ANISOTROPIC GRAVITY \page{43}
       \endtitle
\vskip5pt
\head{} {\bf II.1\ Gauge Fields on Locally Anisotropic Spaces} \page{43}
     \endhead
\vskip5pt
\head{} {\bf II.2\ Yang--Mills Equations on Locally Anisotropic Spaces}
   \page{46}\endhead
\vskip5pt
\head{} {\bf II.3\ Locally Anisotropic Gauge Gravity}\page{49} \endhead
\subhead{} II.3.1\ The bundle of linear locally anisotropic frames \page{49}
   \endsubhead
\subhead{} II.3.2\ The bundle of affine locally adapted frames \page{50}
 \endsubhead
\vskip5pt
\head{}  {\bf II.4\ Nonlinear De Sitter Gauge \newline
Locally Anisotropic  Gravity} \page{52} \endhead
\subhead{} II.4.1\ Nonlinear gauge gravity of
the de Sitter group \page{52}   \endsubhead
\subhead{} II.4.2\ Dynamics of the nonlinear locally anisotropic de Sitter
 gravity \page{53}\endsubhead
\vskip5pt
\title\chapter{3}NEARLY AUTOPARALLEL MAPS AND CONSERVATION LAWS
\page{56}\endtitle
\vskip5pt
\head{} {\bf III.1\ Nearly Autoparallel Maps of Locally Anisotropic Spaces}
   \page{56} \endhead
\vskip5pt
\head{} {\bf III.2\ Classification of Nearly Autoparallel Maps of LA--Spaces}
  \page{59} \endhead
\vskip5pt
\head{} {\bf III.3\ Nearly Autoparallel Tensor Integral on LA--Spaces}
   \page{67}\endhead
\vskip5pt
\head{} {\bf III.4\ On Conservation Laws on LA--Spaces} \page{70} \endhead
\subhead{} {III.4.1\ Nonzero divergence of energy--momentum d--tensor}
   \page{70} \endsubhead
\subhead{} III.4.2\ Deformations d--tensors and tensor integral conservation
 laws \page{71} \endsubhead
\vskip5pt
\head{} {\bf III.5\ NA--Conservation Laws in LA--Gravity} \page{72} \endhead
\vskip5pt
\head{} {\bf REFERENCES} \page{75} \endhead
\endtoc
\vskip40pt
\newpage
\topmatter
 \title  Preface \endtitle
\endtopmatter

The generalizations of Finsler spaces \cite{Finsler 1918} (see fundamental
 contributions and references in \cite{Cartan 1934} and \cite{Rund 1959})
 with various applications in physics,  chemistry, biology, ecology etc
 have been a subject for extensive study in the last two decades
 \cite{Miron and Anastasiei 1987, 1994},
 \cite{Asanov 1985}, \cite{Matsumoto 1986}, \cite{Bejancu 1990},
\cite{Miron and Kavaguchi 1996} and \cite{Antonelli and Miron (eds) 1996}.
 There were developed a number of theories of Finsler gravity and gauge
 fields or of stochastic processes in generalized Finsler spaces.
 Recently, a deal of attention is attracted by the problem of formulating of
the geometric background of classical and quantum field interactions with
 local anisotropy.  We note some of our contributions based on modeling of
 locally anisotropic physical theories on vector bundles provided with
 nonlinear connection structure \cite{Vacaru 1996},
\cite{Vacaru and Goncharenko 1995} and \cite{Vacaru and Ostaf 1996a}
 which make up a starting point of this work.
 Our approach stands out because it allows a rigorous  definition of
 spinors in generalized Finsler spaces and a corresponding geometric
 treatment of fundamental fields interactions on such spaces.

 Spinor variables on Finsler  spaces have been considered in a heuristic
 manner, for instance, by \cite{Asanov and Ponomarenko 1988} and
\cite{Ono and Takano 1993}. A substantial difficulty that arises in every
 attempt to construct physical models on Finsler spaces is the absence of
 groups of automorphisms even locally on spaces with generic local anisotropy
 and the existence of a new fundamental geometric object the nonlinear
 connection. In consequence it is a chalenging task to introduce spinors as
 Clifford structures or to define conservation laws by using usual variational
 calculus. Our key idea was to concentrate  efforts not on definition of
  various types  of spinor spaces for every particular class of Finsler, and
 theirs generalizations, spaces but to apply the fact that spaces with local
anisotropy can be in general modeled on tangent \cite{Yano and Ishihara 1973}
 or vector \cite{Miron and Anastasiei 1987, 1994} bundles provided with
 nonlinear and distinguished connections and metric structures. If the
 mentioned connections and metric structure are compatible, the Clifford
 and spinor locally anisotropic bundles can be introduced similarly as for
 curved spaces or vector bundles \cite{Geroch 1958}, 
 \cite{Karoubi 1978} and \cite{Turtoi 1989} with that distinction that the
 geometric constructions must be adapted to the nonlinear connection (see a
 study of  Clifford and spinor structures on generalized Lagrange and Finsler
  spaces in \cite{Vacaru 1996} and an geometric approach to locally
 anisotropic gauge fields and gauge like gravity in
 \cite{Vacaru and Goncharenko 1995}).

 Regarding the definition of conservation laws on locally anisotropic spaces,
 we propose to introduce into consideration a new class of maps with
deformation of connections of both type of locally isotropic or anisotropic
 curved spaces. We are inspired by the idea \cite{Petrov 1970} (see also
 H. Poincare works \cite{Poincare 1905, 1954} on conventionality of concepts of
 geometrical space--time and physical theories) that geometric
 constuctions and physical models can be locally equivalently modeled on
 arbitrary  curved, or (if some well defined conditions are satisfied) flat
 space
 by using corresponding generalizations of conformal transforms. Partially the
 geometric aspect of the Petrov's program on modeling of field interactions
 was realized in the monograph  \cite{Sinyukov 1979} where the
 theory of nearly geodesic maps  of Riemannian and affine connection spaces
 is formulated. Nearly geodesic
 maps generalizes the class of conformal, geodesic and concircular transforms
 \cite{Schouten and Struik 1938}, \cite{Yano 1940}, \cite{Vr\^anceanu 1977}
 and \cite{Mocanu 1955}. We extended
(see \cite{Vacaru 1987}, \cite{Vacaru 1992}, \cite{Vacaru 1994} and
 \cite{Vacaru and Ostaf 1996b}) the Sinyukov's theory for spaces
 with torsion and nonmetricity, the so--called Einstein--Cartan--Weyl spaces,
 by considering nearly autoparallel maps, developed some approaches to
 formulation of conservation laws for physical interactions and definition of
 twistors \cite{Penrose and Rindler 1986} on such spaces and proposed a
 nearly autoparallel map classification of Lagrange and Finsler space in
 paper \cite{Vacaru and Ostaf 1996a}.

After presenting this informal discussion of some basic ideas and results to be
 used in our further considerations, we now turn to a more detailed description
of the content of this work.

The presentation in the Chapter I is organized as follows: The geometry of
nonlinear connections in vector bundles is revewed in section I.1, were the
 explicit formulas for torsions and curvatures on locally anisotropic spaces
 are given and motion equations for fundamental field interactions with local
 anisotropy  are introduced in a geometric manner. Section I.2 is devoted to
 the distinguished Clifford algebras. Clifford bundles and distinguished by
 nonlinear connections spinor structures are defined in section I.3. Almost
complex spinor structures for generalized Lagrange spaces are analyzed in
 section I.4. In section I.5 the spinor techniques is developed for
 distinguished vector spaces. The differential geometry of locally anisotropic
 spinors is considered in section I.6 where distinguished spinor formulas for
 connections torsions and curvatures are presented. The spinor form of field
 equations on locally anisotropic spaces is analyzed in section I.7.

The Capter II is devoted to the geometry of gauge fields and gauge gravity in
 locally anisotropic spaces. In section II.1 we give a geometrical
 interpretation of gauge (Yang-Mills) fields on general locally anisotropic
 spaces.  Section II.2 contains a geometrical definition of anisotropic
 Yang-Mills equations; the variational proof of gauge field
equations is considered  in connection with the "pure" geometrical method of
 introducing field equations. In section II.3 the locally anisotropic gravity
 is reformulated as a gauge theory for nonsemisimple  groups. A model of
  nonlinear de Sitter gauge gravity with local anisotropy is formulated in
 section II.4.

The problem of formulation of conservation laws on locally anisotropic spaces
 is investigated in the framework of the geometry of local 1--1 maps of vector
 bundles provided with nonlinear connection structures and by developing the
 formalism of tensor--integral for locally anisotropic multispaces in the
 Chapter III. Section III.1 is  devoted to the formulation of the theory of
 nearly autoparallel maps of  locally anisotropic spaces. The classification
 of na--maps and formulation of their invariant conditions are given in
section III.2. In section III.3 we define the nearly  autoparallel
 tensor--integral on locally anisotropic multispaces. The question of
 formulation of conservation laws on
spaces with local anisotropy is studied in section III.4. We present a
definition of conservation laws for locally anisotropic gravitational fields
 on nearly autoparallel images of locally anisotropic spaces in section III.5.
\vskip80pt
\newpage
\topmatter
\title\chapter{1} Locally Anisotropic Spinor Spaces\endtitle
\endtopmatter

The purpose of this Chapter is to present an introduction into the geometry
of spinors in generalized Finsler spaces and to propose a geometric framework
 for the theory  field interactions on locally anisotropic (la) spaces starting
 from papers \cite{Vacaru 1996} and \cite{Vacaru and Goncharenko 1995}
(in brief we shall write la--spaces, la--geometry, la--spinors and so on).
 The geometric constructions will be adapted to the nonlinear connection
  (N--connection) structure. We consider the reader to be familiar with
 basic results on differential geometry of bundle spaces
 \cite{Bishop and Crittenden 1964} and \cite{Kobayashi and Nomizu 1963, 1969}
 and note that as a rule all geometric constructions in this work will be
 local in nature.

For our considerations on the la--spinor theory it will be convenient to
extend the \cite{Penrose and Rindler 1984, 1986} abstract index approach
 (see also the \cite{Luehr and Rosenbaum 1974} index free methods) proposed
 for spinors on locally isotropic spaces. We note that for applications in
 mathematical physics we usually we have  dimensions $d>4$ for spaces into
 consideration. In this case the 2--spinor calculus does not play a
 preferential role.
\vskip25pt

\head  I.1 Connections and Metrics in Locally Anisotropic Spaces \endhead

As a preliminary to a discussion of la--spinor formalism we review some
 results and methods of the differential geometry of tangent and vector
bundles provided with nonlinear and distinguished connections and metric
 structures (see regorous results and details in
\cite{Miron and Anastasiei 1987, 1994} and \cite{Yano and Ishihara 1973}).\
 This section serves the twofold purpose of establishing of abstract index
 denotations and starting the geometric backgrounds which are used in the
 next sections of the Chapter. Combersome proofs and  calculations
 will not be presented.
\vskip15pt

\subhead I.1.1\ Nonlinear connections and distinguishing of geometric objects
\endsubhead

Let ${\Cal E=}$ $\left( E,p,M,Gr,F\right) $ be a locally trivial vector
bundle, v--bundle, where $F={\Bbb R }^m$ (a real vector space of dimension $%
m,\dim F=m,$ ${\Bbb R }$ denotes the real number field) is the typical
fibre, the structural group is chosen to be the group of automorphisms of $%
{\Bbb R }^m$ , i.e. $Gr=GL( m,{\Bbb R }),$ and $p:E\rightarrow M$
is a differentiable surjection of a differentiable manifold $E$ (total
space, $\dim E=n+m)$ to a differentiable manifold $M$ $\left( \text{base
space, }\dim M=n\right) .$ Local coordinates on ${\Cal E}$ are denoted as $%
u^{{\boldsymbol \alpha }}=\left( x^{{\bold i}},y^{{\bold a\ }}\right) ,$
 or in brief $%
{\bold u} = \left( x,y\right) ,$
where boldfaced indices will be considered as
coordinate ones for which the Einstein summation rule holds (Latin indices
$\bold i$,$\bold j$,$\bold k$,$...=1,2,...,n$ will parametrize coordinates of
 geometrical objects with respect to a base space $M,$ Latin indices
$\bold a$,$\bold b$,$\bold c$,$...=1,2,...,m$ will parametrize fibre
coordinates of geometrical objects and
Greek indices $\boldsymbol \alpha$,$\boldsymbol \beta$,
$\boldsymbol \gamma$,$...$ are considered as
cumulative ones for coordinates of objects defined on the total space of a
v-bundle). We shall correspondingly use abstract indices $\alpha =(i,a),$ $%
\beta =(j,b),\gamma =(k,c),...$ in the Penrose manner
\cite{Penrose and Rindler 1984, 1986} in order to
mark geometical objects and theirs (base, fibre)-components or, if it will
be convenient, we shall consider boldfaced letters (in the main for pointing
to the operator character of tensors and spinors into consideration) of type
${\bold A\equiv }A=\left( A^{(h)},A^{(v)}\right) , {\bold b} = \left(
b^{(h)},b^{(v)}\right) ,...,{\bold R},{\boldsymbol \omega},
 {\boldsymbol \Gamma},...$ for geometrical
objects on ${\Cal E}$ and theirs splitting into horizontal (h), or base, and
vertical (v), or fibre, components. For simplicity, we shall prefer writing
out of abstract indices instead of boldface ones if this will not
give rise to ambiguities.

Coordinate trans\-forms
$u^{{\boldsymbol \alpha}^{\prime }\ } = u^{{\boldsymbol \alpha} ^{\prime
}\ }\left( u^{{\boldsymbol \alpha }}\right) $ on ${\Cal E}$ are writ\-ten as
$$\left(x^{{\bold i}},
y^{{\bold a}}\right) \rightarrow \left( x^{{\bold i}^{\prime }\ }, y^{%
{\bold a}^{\prime }}\right) ,$$ where $$x^{{\bold i}^{\prime }\ } =
x^{{\bold i}^{\prime }\ }(x^{{\bold i}}), y^{{\bold a}^{\prime }\ } =
 M_{{\bold a\ }}^{{\bold a}^{\prime }}(x^{{\bold i\ }})
 y^{{\bold a}}$$ and matrices $M_{{\bold a\ }}^{%
{\bold a}^{\prime }}(x^{{\bold i\ }})\in GL\left( m,{\Bbb R}\right) $ are
functions of necessary smoothness class.

A local coordinate parametrization of ${\Cal E}$ naturally defines a
coordinate basis%
$$
\frac \partial {\partial u^\alpha }=\left( \frac \partial {\partial
x^i},\frac \partial {\partial y^a}\right) , \tag 1.1$$
in brief we shall write $\partial _\alpha =(\partial _i,\partial _a),$ and
the reciprocal to (1.1) coordinate basis
$$
du^\alpha =(dx^i,dy^a), \tag 1.2$$
or, in brief, $d^\alpha =(d^i,d^a),$ which is uniquely defined from the
equations%
$$
d^\alpha \circ \partial _\beta =\delta _\beta ^\alpha ,
$$
where $\delta _\beta ^\alpha $ is the Kronecher symbol and by ''$\circ $$"$
we denote the inner (scalar) product in the tangent bundle ${\Cal TE}.$

The concept of {\bf nonlinear connection,}  is
fundamental in the geometry of la--spaces. It came from Finsler
 geometry, \cite{Cartan 1934} and \cite{Kawaguchi 1952, 1956}, and was globally
 defined by \cite{Barthel 1963} (see a detailed study and basic
references in \cite{Miron and Anastasiei 1987, 1994}).

In a v--bundle ${\Cal E}$ we can consider the distribution
$\{N:E_u\rightarrow H_uE,T_uE=H_uE\oplus V_uE\}$ on $E$ being a global
decomposition, as a Whitney sum, into horizontal,${\Cal HE},$ and
vertical, ${\Cal VE,}$ subbundles of the tangent bundle ${\Cal TE}:$
$$
{\Cal TE}={\Cal HE}\oplus{\Cal VE}. \tag 1.3 $$

\proclaim{\bf Definition 1.1 }
A nonlinear connection in the vector bundle $(E,p,M)$ is a splitting on the
 left of the exact sequence
$$0 \longrightarrow \quad VE \quad
{\overset i \to \longrightarrow} \quad TE \quad
{ \longrightarrow} \quad TE/VE \quad \longrightarrow 0$$
that is a morphism of vector bundles $C: TE \to VE$ such that $C \circ i$ is
 the identity on $VE.$
\endproclaim

Locally a N--connection in ${\Cal E}$ is given by it components $N_{{\bold i}%
}^{{\bold a}}({\bold u})=N_{{\bold i}}^{{\bold a}}({\bold x},{\bold y})$
(in brief we shall write $N_i^a(u)=N_i^a(x,y)$ )
with respect to bases (1.1) and (1.2)):
$$
{\bold N}=N_i^a(u)d^i\otimes \partial _a.
$$

We note that a linear connection in a v-bundle ${\Cal E}$ can be considered
as a particular case of a N--connection when $N_i^a(x,y)=K_{bi}^a\left(
x\right) y^b,$ where functions $K_{ai}^b\left( x\right) $ on the base $M$
are called the Christoffel coefficients.

To coordinate locally geometric constructions with the global splitting of $%
{\Cal E}$ defined by a N-connection structure, we have to introduce a
locally adapted basis (la--basis, la--frame):
$$
\frac \delta {\delta u^\alpha }=\left( \frac \delta {\delta x^i}=\partial
_i-N_i^a\left( u\right) \partial _a,\frac \partial {\partial y^a}\right) ,
\tag 1.4 $$
or, in brief,$\delta _\alpha =\left( \delta _i,\partial _a\right) ,$ and it
dual la--basis%
$$
\delta u^\alpha =\left( \delta x^i=dx^i,\delta y^a+N_i^a\left( u\right)
dx^i\right) , \tag 1.5$$
or, in brief, $\ \delta \ ^\alpha =$ $\left( d^i,\delta ^a\right) .$

The{\bf \ nonholonomic coefficients }
${\bold w}=\{w_{\beta \gamma }^\alpha
\left( u\right) \}$ of la--frames are defined as%
$$
\left[ \delta _\alpha ,\delta _\beta \right] =\delta _\alpha \delta _\beta
-\delta _\beta \delta _\alpha =w_{\beta \gamma }^\alpha \left( u\right)
\delta _\alpha . \tag 1.6 $$

The {\bf algebra of tensorial distinguished fields} $DT\left( {\Cal E}%
\right) $ (d--fields, d--tensors, d--objects) on ${\Cal E}$ is introduced as
the tensor algebra ${\Cal T} =\{ {\Cal T}_{qs}^{pr}\}$ of the v--bundle
${\Cal E}_{\left( d\right) },$
 $p_d:{\Cal HE\oplus VE\rightarrow E.\,\ }$ An
element ${\bold t}\in {\Cal T}_{qs}^{pr},$ d--tensor field of type %$\left(
%\begin{array}{cc}
%p & r \\
%q & s
%\end{array}
%\right) ,$
$\left(\smallmatrix p&r\\ q&s\endsmallmatrix\right)$
can be written in local form as%
$$
{\bold t}=t_{j_1...j_qb_1...b_r}^{i_1...i_pa_1...a_r}\left( u\right) \delta
_{i_1}\otimes ...\otimes \delta _{i_p}\otimes \partial _{a_1}\otimes
...\otimes \partial _{a_r}\otimes d^{j_1}\otimes ...\otimes d^{j_q}\otimes
\delta ^{b_1}...\otimes \delta ^{b_r}.
$$

We shall respectively use denotations ${\Cal X\left( E\right) }$ (or ${\Cal %
X\ } {\left( M\right) ),\ } \Lambda ^p\left( {\Cal E}\right) $ $\left(
\text{ or }\Lambda ^p\left( M\right) \right) $ and ${\Cal F\left( E\right) }$
(or ${\Cal F}$ $\left( M\right) $) for the module of d--vector fields on $%
{\Cal E}$ (or $M$ ), the exterior algebra of p--forms on ${\Cal E\ }$(or $M)$
and the set of real functions on ${\Cal E\ }$(or $M).$

In general, d-objects on ${\Cal E\ }$ are introduced as geometric objects
with various group and coordinate transforms coordinated with the
N-connection structure on ${\Cal E}.$\ For example, a d-connection $D$ on
 ${\Cal E}$ is defined as a linear connection $D$ on $E$ conserving
under a parallelism the global decomposition (1.3) into horizontal and
vertical subbundles of ${\Cal TE}$ .

A N-connection in ${\Cal E}$ induces a corresponding decomposition of
d--tensors into sums of horizontal and vertical parts, for example, for every
d--vector $X\in {\Cal X\left( E\right) }$ and 1-form $\widetilde{X}\in
\Lambda ^1\left( {\Cal E}\right) $ we have respectively
$$
X=hX+vX \text{\quad and \quad }\widetilde{X}=h\widetilde{X}+v%
\widetilde{X}~. \tag 1.7 $$
In consequence, we can associate to every d-covariant derivation along the
 d-vector (1.7), $D_X  =X\circ D ,$
two new operators of h- and v-covariant derivations defined
respectively as
$$
D_X^{(h)}Y=D_{hX}Y\quad \text{ and \quad }D_X^{\left( v\right) }Y=D_{vX}Y,
 \forall Y{ \in }{\Cal X\left( E\right) ,}
$$
for which the following conditions hold:%
$$
D_XY = D_X^{(h)}Y{ + }D_X^{(v)}Y, \tag 1.8$$ $$
D_X^{(h)}f=(hX) f\text{ \quad and\quad }D_X^{(v)}f=(vX) f,\quad X,Y%
 \in {\Cal X\left( E\right) ,}f\in {\Cal F}\left(M\right) .
$$

We define a {\bf metric structure }${\bold G}$ in the total space $E$ of
v-bundle ${\Cal E}$ = $\left( E,p,M\right) $ over a connected and paracompact
base $M$ as a symmetric covariant tensor field of type $\left( 0,2\right) $,
$G_{\alpha \beta ,}$ being nondegenerate and of constant signature on
$E.$
\proclaim {\bf Definition 1.2}
Nonlinear connection ${\bold N}$ and metric ${\bold G}$ structures on
${\Cal E}$
are mutually compatible it there are satisfied the conditions:
$$
{\bold G}\left( \delta _i,\partial _a\right) =0, \text{ or equivalently, }%
G_{ia}\left( u\right) - N_i^b\left( u\right) h_{ab}\left( u\right) =0,
\tag 1.9$$
where $h_{ab}={\bold G}\left( \partial _a,\partial _b\right) $ and $G_{ia}=%
{\bold G}\left( \partial _i,\partial _a\right) .$ \endproclaim

 From (1.9) one follows
$$
N_i^b\left( u\right) =h^{ab}\left( u\right) G_{ia}\left( u\right)
\tag 1.10$$
(the matrix $h^{ab}$ is inverse to $h_{ab}).$ In consequence one
obtains the following decomposition of metric :%
$$
{\bold G}(X,Y) = h{\bold G}(X,Y)+ v{\bold G}(X,Y),
\tag 1.11$$
where the d--tensor $h{\bold G}(X,Y) = G (hX,hY)$ is of type
%$\left(
%\begin{array}{cc}
%0 & 0 \\
%2 & 0
%\end{array}
%\right) $
$\left(\smallmatrix 0&0\\2&0\endsmallmatrix\right)$
and the d--tensor $v{\bold G}(X,Y) = {\bold  G}(vX,vY)$ is of type
% $\left(
%\begin{array}{cc}
%0 & 0 \\
%0 & 2
%\end{array}
%\right) .$
 With respect to la--basis (1.5) the d--metric (1.11) is written as%
$$
{\bold G}=g_{\alpha \beta }\left( u\right) \delta ^\alpha \otimes \delta
^\beta =g_{ij}\left( u\right) d^i\otimes d^j+h_{ab}\left( u\right) \delta
^a\otimes \delta ^b,
\tag 1.12 $$
where $g_{ij}={\bold G}\left( \delta _i,\delta _j\right) .$

A metric structure of type (1.11) (equivalently, of type (1.12)) or a metric
on $E$ with components satisfying constraints (1.9), equivalently (1.10))
defines an adapted to the given N-connection inner (d-scalar) product on the
tangent bundle ${\Cal TE}.$

\proclaim{\bf Definition 1.3}
We shall say that a d-connection $\widehat{D}_X$ is compatible with the
d-scalar product on ${\Cal TE}$\ (i.e. is a standard d-connection) if
$$
\widehat{D}_X\left( {\bold X\cdot Y}\right) =\left( \widehat{D}_X{\bold Y}%
\right) \cdot {\bold Z}+Y{\cdot }\left( \widehat{D}_X{\bold Z}\right) ,\forall
{\bold X},{\bold Y},{\bold Z}{ \in }{\Cal X}\left( E\right) .
\tag 1.13$$
\endproclaim

An arbitrary d-connection $D_X$ differs from the standard one $\widehat{D}_X$
by an operator $\widehat{P}_X\left( u\right) =\{X^\alpha \widehat{P}_{\alpha
\beta }^\gamma \left( u\right) \},$ called the deformation d-tensor with
respect to $\widehat{D}_X,$ which is just a d-linear transform of
 ${\Cal E}_u,$ $\forall u\in {\Cal E}.$ The explicit form of $\widehat{P}_X$
can be found by using the corresponding axiom defining linear connections
\cite{Luehr and Rosenbaum 1974}
$$
\left( D_X-\widehat{D}_X\right) fZ=f\left( D_X-\widehat{D}_X\right) Z,
$$
written with respect to la--bases (1.4) and (1.5). From the last expression
we obtain
$$
\widehat{P}_X\left( u\right) =\left[ (D_X-\widehat{D}_X)\delta _\alpha
\left( u\right) \right] \delta ^\alpha \left( u\right) ,
\tag 1.14 $$
therefore
$$
D_XZ =\widehat{D}_XZ + \widehat{P}_XZ  .
\tag 1.15$$

\proclaim{\bf Definition 1.4}
A d--connection $D_X$ is {\bf metric (}or {\bf compatible } with metric
${\bold G}$) on ${\Cal E}$ if%
$$
D_X{\bold G} =0,\forall X \in {\Cal X\left( E\right) } .
\tag 1.16 $$ \endproclaim

If there is a complex structure $J_{\alpha}^{\ \beta},\ J\dot J=-I,$ being
compatible with a metric of type (1.12) and a d--connection $D$ on
tangent bundle $TM,$  when conditions
$$ J_{\alpha}^{\ \beta} J_{\gamma}^{\ \delta} G_{\beta \delta} =
G_{\alpha \gamma} \ \text{ and } \ D_{\alpha} J^{\gamma}_{\ \beta}=0 $$
are satisfied, one considers that on $TM$ it is defined an almost Hermitian
model, $H^{2n}$-model, of a generalized Lagrange space,
GL-space \cite{Miron and Anastasiei 1987, 1994}.

Locally adapted components $\Gamma _{\beta \gamma }^\alpha $ of a
d-connection $D_\alpha =(\delta _\alpha \circ D)$ are defined by the
equations%
$$
D_\alpha \delta _\beta =\Gamma _{\alpha \beta }^\gamma \delta _\gamma ,
$$
from which one immediately follows%
$$
\Gamma _{\alpha \beta }^\gamma \left( u\right) =\left( D_\alpha \delta
_\beta \right) \circ \delta ^\gamma .
\tag 1.17$$

The operations of h- and v-covariant derivations, $D_k^{(h)}=%
\{L_{jk}^i,L_{bk\;}^a\}$ and $D_c^{(v)}=\{C_{jk}^i,C_{bc}^a\}$ (see (1.8)),
are introduced as corresponding h- and v-paramet\-ri\-za\-ti\-ons of (1.17):%
$$
L_{jk}^i=\left( D_k\delta _j\right) \circ d^i,\quad L_{bk}^a=\left(
D_k\partial _b\right) \circ \delta ^a
\tag 1.18$$
and%
$$
C_{jc}^i=\left( D_c\delta _j\right) \circ d^i,\quad C_{bc}^a=\left(
D_c\partial _b\right) \circ \delta ^a.
\tag 1.19$$
A set of components (1.18) and (1.19), $D\Gamma =\left(
L_{jk}^i,L_{bk}^a,C_{jc}^i,C_{bc}^a\right) ,$ completely defines the local
action of a d--connection $D$ in ${\Cal E}.$ For instance, taken a d--tensor
field of type\
%$\left(
%\begin{array}{cc}
%1 & 1 \\
%1 & 1
%\end{array}
%\right) ,$
$\left(\smallmatrix a&b\\ c&d\endsmallmatrix\right)$\
${\bold t}=t_{jb}^{ia}\delta _i\otimes \partial _a\otimes \partial
^j\otimes \delta ^b,$ and a d-vector ${\bold X}=X^i\delta _i+X^a\partial _a$
we have%
$$
D_X{\bold t=}D_X^{(h)}{\bold t+}D_X^{(v)}{\bold t=}\left(
X^kt_{jb|k}^{ia}+X^ct_{jb\perp c}^{ia}\right) \delta _i\otimes \partial
_a\otimes d^j\otimes \delta ^b,
$$
where the h-covariant derivative is written as%
$$
t_{jb|k}^{ia}=\frac{\delta t_{jb}^{ia}}{\delta x^k}%
+L_{hk}^it_{jb}^{ha}+L_{ck}^at_{jb}^{ic}-L_{jk}^ht_{hb}^{ia}-L_{bk}^ct_{jc}^{ia}
$$
and the v-covariant derivative is written as%
$$
t_{jb\perp c}^{ia}=\frac{\partial t_{jb}^{ia}}{\partial y^c}%
+C_{hc}^it_{jb}^{ha}+C_{dc}^at_{jb}^{id}-C_{jc}^ht_{hb}^{ia}-
C_{bc}^dt_{jd}^{ia}.
$$
For a scalar function $f\in {\Cal F\left( E\right) }$ we have
$$
D_k^{(h)}=\frac{\delta f}{\delta x^k}=\frac{\partial f}{\partial x^k}-N_k^a
\frac{\partial f}{\partial y^a}\text{ and }D_c^{(v)}f=\frac{\partial f}{%
\partial y^c}.
$$

We emphasize that the geometry of connections in a v-bundle ${\Cal E}$ is
very reach. If a triple of fundamental geometric objects $\left( N_i^a\left(
u\right) ,\Gamma _{\beta \gamma }^\alpha \left( u\right) ,G_{\alpha \beta
}\left( u\right) \right) $ is fixed on ${\Cal E},$ really, a
multiconnection structure (with corresponding different rules of covariant
derivation, which are, or not, mutually compatible and with the same, or
not, induced d-scalar products in ${\Cal TE)}$ is defined on this v-bundle.
For instance, we enumerate some of connections and covariant derivations
which can present interest in investigation of locally anisotropic
gravitational and matter field interactions:

\item 1.
Every N-connection in ${\Cal E,}$ with coefficients $N_i^a\left(
x,y\right) $ being differentiable on y-variables, induces a structure of
linear connection $\widetilde{N}_{\beta \gamma }^\alpha ,$ where $\widetilde{%
N}_{bi}^a=\frac{\partial N_i^a}{\partial y^b}$ and $\widetilde{N}%
_{bc}^a\left( x,y\right) =0.$ For some $Y\left( u\right) =Y^i\left( u\right)
\partial _i+Y^a\left( u\right) \partial _a$ and $B\left( u\right) =B^a\left(
u\right) \partial _a$ one writes%
$$
D_Y^{(\widetilde{N})}B=\left[ Y^i\left( \frac{\partial B^a}{\partial x^i}+%
\widetilde{N}_{bi}^aB^b\right) +Y^b\frac{\partial B^a}{\partial y^b}\right]
\frac \partial {\partial y^a}.
$$

\item 2. The d-connection of Berwald type \cite{Berwald 1926}
$$
\Gamma _{\beta \gamma }^{(B)\alpha }=\left( L_{jk}^i,\frac{\partial N_k^a}{%
\partial y^b},0,C_{bc}^a\right) , \tag 1.20 $$
where
$$ \split
L_{.jk}^i\left( x,y\right) =\frac 12g^{ir}\left( \frac{\delta g_{jk}}{\delta
x^k}+\frac{\delta g_{kr}}{\delta x^j}-\frac{\delta g_{jk}}{\delta x^r}%
\right) ,\\
C_{.bc}^a\left( x,y\right) =\frac 12h^{ad} \left( \frac{\partial h_{bd}}
{\partial y^c}+ \frac{\partial h_{cd}}{\partial y^b}-\frac{\partial h_{bc}}
{\partial y^d}\right) ,\endsplit \tag 1.21 $$

which is hv-metric, i.e. $D_k^{(B)}g_{ij}=0$ and $D_c^{(B)}h_{ab}=0.$

\item 3.  The canonical d-connection ${\bold \Gamma}^{(c)}$ associated to a
metric ${\bold G}$ of type (1.12) $\Gamma _{\beta \gamma }^{(c)\alpha } =
  \left(
L_{jk}^{(c)i},L_{bk}^{(c)a},C_{jc}^{(c)i},C_{bc}^{(c)a}\right) ,$ with
coefficients%
$$ \multline
L_{jk}^{(c)i}=L_{.jk}^i,C_{bc}^{(c)a}=C_{.bc}^a \text{ (see (1.21)}\\
L_{bi}^{(c)a}=\widetilde{N}_{bi}^a+\frac 12h^{ac}\left( \frac{\delta h_{bc}}{%
\delta x^i}-\widetilde{N}_{bi}^dh_{dc}-\widetilde{N}_{ci}^dh_{db}\right) ,\
C_{jc}^{(c)i}=\frac 12g^{ik}\frac{\partial g_{jk}}{\partial y^c}.\endmultline
\tag 1.22$$
This is a metric d--connection which satisfies conditions
$$
D_k^{(c)}g_{ij}=0,D_c^{(c)}g_{ij}=0,D_k^{(c)}h_{ab}=0,D_c^{(c)}h_{ab}=0.
$$

\item 4.  We can consider N--adapted Christoffel d--symbols%
$$
\widetilde{\Gamma }_{\beta \gamma }^\alpha =\frac 12G^{\alpha \tau }\left(
\delta _\gamma G_{\tau \beta }+\delta _\beta G_{\tau \gamma }-\delta
G_{\beta \gamma }\right) ,
\tag 1.23$$
which have the components of d-connection $\widetilde{\Gamma }_{\beta \gamma
}^\alpha =\left( L_{jk}^i,0,0,C_{bc}^a\right) ,$ with $L_{jk}^i$ and $%
C_{bc}^a$ as in (1.21) if $G_{\alpha \beta }$ is taken in the form
(1.12).

Arbitrary linear connections on a v--bundle ${\Cal E}$ can be also
characterized by theirs deformation tensors (see (1.15)) with respect, for
instance, to d-connect\-i\-on (1.23):%
$$
\Gamma _{\beta \gamma }^{(B)\alpha }=\widetilde{\Gamma }_{\beta \gamma
}^\alpha +P_{\beta \gamma }^{(B)\alpha },\Gamma _{\beta \gamma }^{(c)\alpha
}=\widetilde{\Gamma }_{\beta \gamma }^\alpha +P_{\beta \gamma }^{(c)\alpha }
$$
or, in general,%
$$
\Gamma _{\beta \gamma }^\alpha =\widetilde{\Gamma }_{\beta \gamma }^\alpha
+P_{\beta \gamma }^\alpha ,
$$
where $P_{\beta \gamma }^{(B)\alpha },P_{\beta \gamma }^{(c)\alpha }$ and $%
P_{\beta \gamma }^\alpha $ are respectively the deformation d-tensors of
d-connect\-i\-ons (1.20),\ (1.22), or of a general one.
\vskip15pt

\subhead{I.1.2 Torsions and curvatures of nonlinear and
 distinguished connections}\endsubhead

The curvature ${\bold \Omega },$ of a nonlinear connection ${\bold N}$ in a
v--bundle ${\Cal E}$ can be defined as the Nijenhuis tensor field $N_v\left(
X,Y\right) $ associated to ${\bold N}:$
$$
{\bold \Omega }=N_v={ \left[ vX,vY\right] +v\left[ X,Y\right] -v\left[
vX,Y\right] -v\left[ X,vY\right] ,X,Y}\in {\Cal X\left( E\right) .}
$$
In local form one has%
$$
{\bold \Omega }=\frac 12\Omega _{ij}^ad^i\bigwedge d^j\otimes \partial _a,
$$
where%
$$
\Omega _{ij}^a=\frac{\partial N_i^a}{\partial x^j}-\frac{\partial N_j^a}{%
\partial x^i}+N_i^b\widetilde{N}_{bj}^a-N_j^b\widetilde{N}_{bi}^a.
\tag 1.24$$

The torsion ${\bold T}$ of a d--connection ${\bold D}$ in ${\Cal E}$ is defined
by the equation%
$$
{\bold T\left( X,Y\right) =XY_{\circ }^{\circ }T\doteq }D_X{\bold Y-}D_Y
{\bold X\ - \left[ X,Y\right] .} \tag 1.25$$
One holds the following h- and v--decompositions
$$
{\bold T}\left( X,Y\right) =T\left( hX,hY\right) +T\left( hX,vY\right) +T\left(
vX,hY\right) +T\left( vX,vY\right) . \tag 1.26$$
We consider the projections:$ h{\bold T}\left( X,Y\right) ,vT\left(
hX,hY\right) ,hT\left( hX,hY\right) ,...$ and say that, for instance,
$h{\bold T}\left( hX,hY\right) $ is the h(hh)-torsion of
${\bold D}$ , $v{\bold T}\left(hX,hY\right) $ is the
v(hh)-torsion of ${\bold D}$ and so on.

The torsion (1.25) is locally determined by five d--tensor fields, torsions,
defined as
$$\multline
T_{jk}^i=h{\bold T}\left( \delta _k,\delta _j\right) \cdot d^i,\quad T_{jk}^a=%
 v{\bold T}\left( \delta _k,\delta _j\right) \cdot \delta ^a,\\
P_{jb}^i=h{\bold T}\left( \partial _b,\delta _j\right) \cdot d^i,\quad
P_{jb}^a=v{\bold T}\left( \partial _b,\delta _j\right) \cdot \delta ^a,\quad
S_{bc}^a=v{\bold T}\left( \partial _c,\partial _b\right) \cdot \delta ^a.
\endmultline \tag 1.27$$
Using formulas (1.4),(1.5),(1.24) and (1.25) we can computer in explicit
form the components of torsions (1.26) for a d--connection of type (1.18) and
(1.19):
$$\split
T_{.jk}^i=T_{jk}^i=L_{jk}^i-L_{kj}^i,\ T_{ja}^i=C_{.ja}^i,\
T_{aj}^i=-C_{ja}^i,\ T_{.bc}^a=S_{.bc}^a=C_{bc}^a-C_{cb}^a,\\
 T_{.ja}^i=0,\
T_{.ij}^a=\frac{\delta N_i^a}{\delta x^j}-\frac{\delta N_j^a}{\delta x^i}%
,\ T_{.bi}^a=P_{.bi}^a=\frac{\partial N_i^a}{\partial y^b} - L_{.bj}^a,\
 T_{.ib}^a=-P_{.bi}^a.  \endsplit \tag 1.28 $$

The curvature ${\bold R}$ of a d--connection in ${\Cal E}$ is defined by the
equation
$$
{\bold R}\left( X,Y\right) Z=XY_{\bullet }^{\bullet }R\bullet Z =
 D_XD_Y{\bold Z} -D_YD_X{\bold Z-}D_{[X,Y]}{\bold Z}. \tag 1.29$$
One holds the next properties for the h- and v-decompositions of curvature:%
$$\split
v{\bold R}\left( X,Y\right) hZ=0,\ h{\bold R}\left( X,Y\right) vZ=0,\\
{\bold R}\left( X,Y\right) Z=hR\left( X,Y\right) hZ+vR\left( X,Y\right) vZ.%
 \endsplit \tag 1.30 $$
From (1.29) and the equation ${\bold R}\left( X,Y\right) =-R\left( Y,X\right)$
we get that the curvature of a d-con\-nec\-ti\-on ${\bold D}$ in ${\Cal E}$ is
completely determined by the following six d--tensor fields:%
$$\split
R_{h.jk}^{.i}=d^i\cdot {\bold R}\left( \delta _k,\delta _j\right) \delta
_h,\ R_{b.jk}^{.a}=\delta ^a\cdot {\bold R}\left( \delta _k,\delta _j\right)
\partial _b,\
P_{j.kc}^{.i}=d^i\cdot {\bold R}\left( \partial _c,\partial _k\right) \delta
_j,\\
P_{b.kc}^{.a}=\delta ^a\cdot {\bold R}\left( \partial _c,\partial
_k\right) \partial _b,\
S_{j.bc}^{.i}=d^i\cdot {\bold R}\left( \partial _c,\partial _b\right) \delta
_j,\ S_{b.cd}^{.a}=\delta ^a\cdot {\bold R}\left( \partial _d,\partial
_c\right) \partial _b. \endsplit \tag 1.31$$
By a direct computation, using (1.4),(1.5),(1.18),(1.19) and (1.31) we get
$$ \split
R_{h.jk}^{.i}=\frac{\delta L_{.hj}^i}{\delta x^h}-\frac{\delta L_{.hk}^i}{%
\delta x^j}+L_{.hj}^mL_{mk}^i-L_{.hk}^mL_{mj}^i+C_{.ha}^iR_{.jk}^a,\\
R_{b.jk}^{.a}=\frac{\delta L_{.bj}^a}{\delta x^k}-\frac{\delta L_{.bk}^a}{%
\delta x^j}+L_{.bj}^cL_{.ck}^a-L_{.bk}^cL_{.cj}^a+C_{.bc}^aR_{.jk}^c,\\
P_{j.ka}^{.i}=\frac{\partial L_{.jk}^i}{\partial y^k}-\left( \frac{\partial
C_{.ja}^i}{\partial x^k}%
+L_{.lk}^iC_{.ja}^l-L_{.jk}^lC_{.la}^i-L_{.ak}^cC_{.jc}^i\right)
+C_{.jb}^iP_{.ka}^b,\\
P_{b.ka}^{.c}=\frac{\partial L_{.bk}^c}{\partial y^a}-\left( \frac{\partial
C_{.ba}^c}{\partial x^k}+L_{.dk}^{c%
\,}C_{.ba}^d-L_{.bk}^dC_{.da}^c-L_{.ak}^dC_{.bd}^c\right)
+C_{.bd}^cP_{.ka}^d,\\
S_{j.bc}^{.i}=\frac{\partial C_{.jb}^i}{\partial y^c}-\frac{\partial
C_{.jc}^i}{\partial y^b}+C_{.jb}^hC_{.hc}^i-C_{.jc}^hC_{hb}^i,\\
S_{b.cd}^{.a}=\frac{\partial C_{.bc}^a}{\partial y^d}-\frac{\partial
 C_{.bd}^a}{\partial y^c}+C_{.bc}^eC_{.ed}^a-C_{.bd}^eC_{.ec}^a.
\endsplit \tag 1.32
$$

 We
note that torsions (1.28) and curvatures (1.32) can be computed by particular
 cases of d-connections when d-connections (1.20), (1.22) or (1.24) are used
 instead of (1.18) and (1.19).

The components of the Ricci d-tensor
$$
R_{\alpha \beta }=R_{\alpha .\beta \tau }^{.\tau }
$$
with respect to locally adapted frame (1.5) are as follows:%
$$\split
R_{ij}=R_{i.jk}^{.k},\quad R_{ia}=-^2P_{ia}=-P_{i.ka}^{.k},\\
R_{ai}=^1P_{ai}=P_{a.ib}^{.b},\quad R_{ab}=S_{a.bc}^{.c}.
\endsplit \tag 1.33$$
We point out that because, in general, $^1P_{ai}\neq ~^2P_{ia}$ the Ricci
d--tensor is non symmetric.

Having defined a d-metric of type (1.12) in ${\Cal E}$ we can introduce the
scalar curvature of d--connection ${\bold D}:$
$$
{\overleftarrow{R}}=G^{\alpha \beta }R_{\alpha \beta }=R+S,
\tag 1.34$$
where $R=g^{ij}R_{ij}$ and $S=h^{ab}S_{ab}.$

For our further considerations it will be also useful to consider an alternative
way of definition torsion (1.25) and curvature (1.29) by using the
commutator
$$
\Delta _{\alpha \beta }\doteq \nabla _\alpha \nabla _\beta -\nabla _\beta
\nabla _\alpha =2\nabla _{[\alpha }\nabla _{\beta ]}.
\tag 1.35$$
For components (1.28) of d-torsion we have
$$
\Delta _{\alpha \beta }f=T_{.\alpha \beta }^\gamma \nabla _\gamma f
\tag 1.36$$
for every scalar function $f\,$ on ${\Cal E}.$ Curvature can be introduced
as an operator acting on arbitrary d-vector $V^\delta :$
$$
(\Delta _{\alpha \beta }-T_{.\alpha \beta }^\gamma \nabla _\gamma )V^\delta
=R_{~\gamma .\alpha \beta }^{.\delta }V^\gamma
\tag 1.37$$
(we note that in this Chapter we shall follow conventions of
 \cite{Miron and Anastasiei 1987, 1994} on d--tensors; we can obtain
corresponding  abstract index formulas from
\cite{Penrose and Rindler 1984, 1986} just for a trivial N--connection
 structure and by changing denotations for components of torsion and curvature
 in this manner:\
$T_{.\alpha \beta }^\gamma \rightarrow T_{\alpha \beta }^{\quad
\gamma }$ and $R_{~\gamma .\alpha \beta }^{.\delta }\rightarrow R_{\alpha
\beta \gamma }^{\qquad \delta }).$

Here we also note that torsion and curvature of a d-connection on ${\Cal E}$
satisfy generalized for la--spaces Ricci and Bianchi identities
\cite{Miron and Anastasiei 1987, 1994} which in
terms of components (1.36) and (1.37) are written respectively as%
$$
R_{~[\gamma .\alpha \beta ]}^{.\delta }+\nabla _{[\alpha }T_{.\beta \gamma
]}^\delta +T_{.[\alpha \beta }^\nu T_{.\gamma ]\nu }^\delta =0
\tag 1.38$$
and%
$$
\nabla _{[\alpha }R_{|\nu |\beta \gamma ]}^{\cdot \sigma }+T_{\cdot [\alpha
\beta }^\delta R_{|\nu |.\gamma ]\delta }^{\cdot \sigma }=0.
\tag 1.39$$
Identities (1.38) and (1.39) can be proved by straightforward calculations.

We can also consider a la--generalization of the so--called conformal Weyl
tensor (see, for instance, \cite{Penrose and Rindler 1984}):
$$
C_{\quad \alpha \beta }^{\gamma \delta }=R_{\quad \alpha \beta }^{\gamma
\delta }-\frac 4{n+m-2}R_{\quad [\alpha }^{[\gamma }~\delta _{\quad \beta
]}^{\delta ]}+
\frac 2{(n+m-1)(n+m-2)}{\overleftarrow{R}~\delta _{\quad
[\alpha }^{[\gamma }~\delta _{\quad \beta ]}^{\delta ]}.}
\tag 1.40$$
This d--tensor is conformally invariant on la-spaces provided with
d-connection generated by d--metric structures.
\vskip15pt

\subhead{I.1.3 Field equations for locally anisotropic gravity} \endsubhead

The Einstein equations and the problem of conservation laws on v--bundles
 provided with N--connection structures are considered in
\cite{Miron and Anastasiei 1994}. In work \cite{Vacaru and Goncharenko 1995}
 it was  proved
that the la--gravity can be formulated in a gauge like manner and the
conditions when the Einstein la--gravitational field equations are equivalent
 to a corresponding form of Yang-Mills equations where analized.
In this subsection we shall express the la--gravitational field equations in
 a form more convenient for theirs equivalent reformulation in la--spinor
 variables.

We define d-tensor $\Phi _{\alpha \beta }$ as to satisfy conditions
$$
-2\Phi _{\alpha \beta }\doteq R_{\alpha \beta }-\frac 1{n+m}\overleftarrow{R}%
g_{\alpha \beta }
\tag 1.41$$
which is the torsionless part of the Ricci tensor for locally isotropic
spaces \cite{Penrose and Rindler 1984},
 i.e. $\Phi _\alpha ^{~~\alpha }\doteq 0$.\ The Einstein
equations on la--spaces
$$
\overleftarrow{G}_{\alpha \beta }+\lambda g_{\alpha \beta }=\kappa E_{\alpha
\beta }, \tag 1.42$$
where%
$$
\overleftarrow{G}_{\alpha \beta }=R_{\alpha \beta }-\frac 12\overleftarrow{R}%
g_{\alpha \beta } \tag 1.43$$
is the Einstein d-tensor, $\lambda $ and $\kappa $ are correspondingly the
cosmological and gravitational constants and by $E_{\alpha \beta }$
is denoted the locally anisotropic energy--momentum d--tensor
\cite{Miron and Anastasiei 1987, 1994}, can
be rewritten in equivalent form:%
$$
\Phi _{\alpha \beta }=-\frac \kappa 2(E_{\alpha \beta }-\frac 1{n+m}E_\tau
^{~\tau }~g_{\alpha \beta }). \tag 1.44$$

Because la--spaces generally have nonzero torsions we shall add to (1.44)
(equivalently to (1.42)) a system of algebraic d--field equations with the
source $S_{~\beta \gamma }^\alpha $ being the locally anisotropic spin
density of matter (if we consider a variant of locally anisotropic
Einstein--Cartan theory):%
$$
T_{~\alpha \beta }^\gamma +2\delta _{~[\alpha }^\gamma T_{~\beta ]\delta
}^\delta =\kappa S_{~\alpha \beta .}^\gamma \tag 1.45$$
From (1.38 ) and (1.45) one follows the conservation law of locally
anisotropic spin matter:%
$$
\nabla _\gamma S_{~\alpha \beta }^\gamma -T_{~\delta \gamma }^\delta
S_{~\alpha \beta }^\gamma =E_{\beta \alpha }-E_{\alpha \beta }.
$$

Finally, in this section, we remark that all presented geometric
constructions contain those elaborated for generalized Lagrange spaces
\cite{Miron and Anastasiei 1987, 1994}
(for which a tangent bundle $TM$ is considered instead of a v-bundle ${\Cal E%
}$ ). Here we note that the Lagrange (Finsler) geometry is characterized by
a metric of type (1.12) with components parametized as $g_{ij}=\frac 12%
\frac{\partial ^2{\Cal L}}{\partial y^i\partial y^j}$ $\left( g_{ij}=\frac 12%
\frac{\partial ^2\Lambda ^2}{\partial y^i\partial y^j}\right) $ and $%
h_{ij}=g_{ij},$ where ${\Cal L=L}$ $(x,y)$ $\left( \Lambda =\Lambda \left(
x,y\right) \right) $ is a Lagrangian ( Finsler metric) on $TM$
(see details in \cite{Miron and Anastasiei 1994}, \cite{Matsumoto 1986} and
 \cite{Bejancu 1990}).
\vskip25pt

\head{I.2 Distinguished Clifford Algebras} \endhead

The typical fiber of a v--bundle $\xi _d\ ,\ \pi _d:\ HE\oplus VE\rightarrow E$
is a d-vector space, ${\Cal F}=h{\Cal F}\oplus v{\Cal F},$ split
into horizontal $h{\Cal F}$ and vertical $v{\Cal F}$ subspaces, with metric
 $G(g,h)$ induced by
v-bundle metric (1.12). Clifford algebras (see, for example,
\cite{Karoubi 1978} and \cite{Penrose and Rindler 1986})
 formulated for d--vector spaces will be called
Clifford d--algebras \cite{Vacaru 1996}.
In this section we shall consider the main properties
of Clifford d--algebras. The proof of theorems will be based on the technique
developed in  \cite{Karoubi 1978} correspondingly adapted
to the distinguished character of spaces in consideration.

Let $k$ be a number field (for our purposes $k={\Bbb R}$ or $k={\Bbb C},$\
$\Bbb R$ and ${\Bbb C}$ are respectively real and complex number fields) and
define ${\Cal F}$ as a d--vector space on $k$ provided with a nondegenerate
symmetric quadratic form (metric)\ $G.$ Let $C$ be an algebra on $k$ (not
necessarily commutative) and $j\ :\ {\Cal F}$ $\rightarrow C$ a homomorphism
of underlying vector spaces such that $j(u)^2=\;G(u)\cdot 1\ (1$ is the
unity in algebra $C$ and d-vector $u\in {\Cal F}).$ We are interested in
definition of the pair $\left( C,j\right) $ satisfying the next
universitality conditions. For every $k$-algebra $A$ and arbitrary
homomorphism $\varphi :{\Cal F}\rightarrow A$ of the underlying d-vector
spaces, such that $\left( \varphi (u)\right) ^2\rightarrow G\left( u\right)
\cdot 1,$ there is a unique homomorphism of algebras $\psi \ :\ C\rightarrow
A$ transforming the diagram
%\begin{figure}[htbp]
%\begin{center}
%\begin{picture}(100,50) \setlength{\unitlength}{1pt}
%\thinlines
%\put(0,45){${\cal F}$}
%\put(96,45){C}
%\put(48,2){A}
%\put(50,38){j}
%\put(50,48){ \vector(-1,0){45}}
%\put(50,48){\vector(1,0){45}}
%\put(5,42){\vector(3,-2){45}}
%\put(98,42){\vector(-3,-2){45}}
%\put(15,20){$\varphi$}
%\put(80,20){$\tau$}
%\end{picture}
%\end{center}
%\caption{Diagram 1}
%\end{figure}
%$$\CD
%{\Cal F} \quad {\overset j \to {\longleftrightarrow}}\quad C \\
%@V{\varphi}VV                       @VV{\tau}V\\
%A\quad \Longleftrightarrow \quad A
%\endCD$$
%$$\CD
%{\Cal F}\ \overset j \to  @>>> }\ C \\
%@V{\varphi}VV                         @VV{\tau}V\\
%A
%\endCD $$
$$\CD
{\Cal F} @. \overset j \to {\longleftrightarrow} C\\
@V{\varphi}VV      \qquad  @VV{\tau}V\\
A   @. \Longleftrightarrow A
\endCD$$
into a commutative one.\ The algebra solving this problem will be denoted as
$C\left( {\Cal F},A\right) $ [equivalently as $C\left( G\right) $ or $%
C\left( {\Cal F}\right) $ and called as Clifford d--algebra associated with
pair $\left( {\Cal F},G\right) .$
\proclaim{\bf Theorem 1.1}
The above-presented diagram has a unique solution $\left(
C,j\right) $ up to isomorphism.
\endproclaim

{\bf Proof:} (We adapt for d-algebras that of
\cite{Karoubi 1978}, p. 127.)
For a universal problem the uniqueness is obvious if we
prove the existence of solution $C\left( G\right) $ . To do this we use
tensor algebra ${\Cal L}^{(F)}=\oplus {\Cal L}_{qs}^{pr}\left( {\Cal %
F}\right) $ =$\oplus _{i=0}^\infty T^i\left( {\Cal F}\right) ,$ where $%
T^0\left( {\Cal F}\right) =k$ and $T^i\left( {\Cal F}\right) =k$ and $%
T^i\left( {\Cal F}\right) ={\Cal F}\otimes ...\otimes {\Cal F}$ for $i>0.$
Let $%
I\left( G\right) $ be the bilateral ideal generated by elements of form $%
\epsilon \left( u\right) =u\otimes u-G\left( u\right) \cdot 1$ where $u\in
{\Cal F}$ and $1$ is the unity element of algebra ${\Cal L}\left( {\Cal F}%
\right) .$ Every element from $I\left( G\right) $ can be written as $%
\sum\nolimits_i\lambda _i\epsilon \left( u_i\right) \mu _i,$ where $\lambda
_{i},\mu _i\in {\Cal L}({\Cal F})$ and
$u_i\in {\Cal F}.$ Let $C\left( G\right) $
=${\Cal L}({\Cal F})/I\left( G\right) $ and define
$j:{\Cal F}\rightarrow C\left(
G\right) $ as the composition of monomorphism $i:{{\Cal F}\rightarrow L}^1%
 ({\Cal F})\subset {\Cal L}({\Cal F})$ and projection
$p:{\Cal L}\left( {\Cal F}\right)
\rightarrow C\left( G\right) .$ In this case pair $\left( C\left( G\right)
,j\right) $ is the solution of our problem. From the general properties of
tensor algebras the homomorphism $\varphi :{\Cal F}\rightarrow A$ can be
extended to ${\Cal L}({\Cal F})$ , i.e., the diagram
%\begin{figure}[htbp]
%\begin{center}
%\begin{picture}(100,50) \setlength{\unitlength}{1pt}
%\thinlines
%\put(0,45){${\cal F}$}
%\put(96,45){${\cal L}({\cal F})$}
%\put(48,2){A}
%\put(50,38){i}
%\put(50,48){ \line(-1,0){45}}
%\put(50,48){\vector(1,0){45}}
%\put(5,42){\vector(3,-2){45}}
%\put(98,42){\vector(-3,-2){45}}
%\put(15,20){$\varphi$}
%\put(80,20){$\rho$}
%\end{picture}
%\end{center}
%\caption{Diagram 2}
%\end{figure}
$$\CD
{\Cal F} @. \overset j \to {\longrightarrow} {\Cal L}({\Cal F})\\
@V{\varphi}VV      \qquad  @VV{\rho}V\\
A   @. \Longleftrightarrow A
\endCD$$
is commutative, where $\rho $ is a monomorphism of algebras. Because $\left(
\varphi \left( u\right) \right) ^2=G\left( u\right) \cdot 1,$ then $\rho $
vanishes on ideal $I\left( G\right) $ and in this case the necessary
homomorphism $\tau $ is defined. As a consequence of uniqueness of $\rho ,$
the homomorphism $\tau $ is unique. $\blacksquare$
\vskip10pt
Tensor d-algebra ${\Cal L}(\Cal F )$
 can be considered as a ${\Bbb Z}/2$ graded algebra.
Really, let us in\-tro\-duce ${\Cal L}^{(0)}({\Cal F}) =  \sum_{i=1}^\infty
T^{2i}\left( {\Cal F}\right) $ and ${\Cal L}^{(1)}({\Cal F}) =%
\sum_{i=1}^\infty T^{2i+1}\left( {\Cal F}\right) .$ Setting $I^{(\alpha
)}\left( G\right) =I\left( G\right) \cap {\Cal L}^{(\alpha )}({\Cal F}).$
Define $C^{(\alpha )}\left( G\right) $ as $p\left( {\Cal L}^{(\alpha )}
({\Cal F})\right) ,$
where $p:{\Cal L}\left( {\Cal F}\right) \rightarrow C\left(
G\right) $ is the canonical projection. Then $C\left( G\right)
=C^{(0)}\left( G\right) \oplus C^{(1)}\left( G\right) $ and in consequence
we obtain that the Clifford d--algebra is ${\Bbb Z}/2$ graded.

It is obvious that Clifford d--algebra functorially depends on pair $\left(
{\Cal F},G\right) .$ If $f:{\Cal F}\rightarrow{\Cal F}^{\prime }$ is a
homomorphism of k-vector spaces, such that $G^{\prime }\left( f(u)\right)
=G\left( u\right) ,$ where $G$ and $G^{\prime }$ are, respectively, metrics
on ${\Cal F}$ and ${\Cal F}^{\prime },$ then $f$ induces an homomorphism of
d--algebras%
$$
C\left( f\right) :C\left( G\right) \rightarrow C\left( G^{\prime }\right)
$$
with identities $C\left( \varphi \cdot f\right) =C\left( \varphi \right)
C\left( f\right) $ and $C\left( Id_{{\Cal F}}\right) =Id_{C({\Cal F)}}.$

If ${\Cal A}^{\alpha}$ and ${\Cal B}^{\beta}$ are ${\Bbb Z}/2$--graded
 d--algebras,
then their graded tensorial product $%
{\Cal A}^\alpha \otimes {\Cal B}^\beta $ is defined as a d-algebra for
k-vector d-space ${\Cal A}^\alpha \otimes {\Cal B}^\beta $ with the graded
product induced as $\left( a\otimes b\right) \left( c\otimes d\right)
=\left( -1\right) ^{\alpha \beta }ac\otimes bd,$ where $b\in {\Cal B}^\alpha
$ and $c\in {\Cal A}^\alpha \quad \left( \alpha ,\beta =0,1\right) .$

Now we reformulate for d--algebras the theorem \cite{Chevalley 1955}:

\proclaim{\bf Theorem 1.2}
The Clifford d-algebra $C\left( h{\Cal F}\oplus v{\Cal F},%
g+h\right) $ is naturally isomorphic to $C(g)\otimes C\left( h\right) .$
\endproclaim

{\bf Proof. }Let $n:h{\Cal F}\rightarrow C\left( g\right) $ and
$n^{\prime }:%
v{\Cal F}\rightarrow C\left( h\right) $ be canonical maps and map
$m:h{\Cal F}\oplus v{\Cal F}\rightarrow C(g)\otimes C\left( h\right) $
is defined as $%
m(x,y)=n(x)\otimes 1+1\otimes n^{\prime }(y),$ $x\in
h{\Cal F}, y\in v{\Cal F}.$
We have $\left( m(x,y)\right) ^2=\left[ \left( n\left( x\right) \right)
^2+\left( n^{\prime }\left( y\right) \right) ^2\right] \cdot 1=[g\left(
x\right) +h\left( y\right) ].$ \ Taking into account the universality
property of Clifford d--algebras we conclude that $m$ induces the homomorphism%
$$
\zeta :C\left(
h{\Cal F}\oplus v{\Cal F}, g+h\right) \rightarrow C\left( h{\Cal F},g\right)
\widehat{\otimes }C\left( v{\Cal F},h\right) .
$$
We also can define a homomorphism%
$$
\upsilon :C\left( h{\Cal F},g\right) \widehat{\otimes }C\left( v{\Cal F}%
,h\right) \rightarrow C\left(h{\Cal F}\oplus v{\Cal F}, g+h\right)
$$
by using formula $\upsilon \left( x\otimes y\right) =\delta \left( x\right)
\delta ^{\prime }\left( y\right) ,$ where homomorphysms $\delta $ and $%
\delta ^{\prime }$ are, respectively, induced by imbeddings of $h{\Cal F}$
and $v{\Cal F}$ into $h{\Cal F}\oplus v{\Cal F} :$%
$$
\delta :C\left( h{\Cal F},g\right) \rightarrow C\left( h{\Cal F}\oplus
 v{\Cal F} , g+h\right) ,
$$
$$
\delta ^{\prime }: C\left(
v{\Cal F},h\right) \rightarrow C\left( h{\Cal F}\oplus
v{\Cal F} , g+h\right) .
$$
Because $x\in C^{(\alpha )}\left( g\right) $ and $y\in C^{(\alpha )}\left(
g\right) ,~\delta \left( x\right) \delta ^{\prime }\left( y\right) =\left(
-1\right) ^{\left( \alpha \right) }\delta ^{\prime }\left( y\right) \delta
\left( x\right) .$

Superpositions of homomorphisms $\zeta $ and $\upsilon $ lead to identities%
$$
\upsilon \zeta =Id_{C\left( h{\Cal F}, g\right) \widehat{\otimes }C\left(
v{\Cal F},h\right)} , \tag 1.46
$$
$$
\zeta \upsilon =Id_{C\left( h{\Cal F},g\right) \widehat{\otimes }C\left(
v{\Cal F},h\right)} .
$$
Really, d--algebra $C\left( h{\Cal F}\oplus v{\Cal F},g+h\right) $
is generated by
elements of type $m(x,y).$ Calculating
$$
\upsilon \zeta \left( m\left( x,y\right) \right) =\upsilon \left( n\left(
x\right) \otimes 1+1\otimes n^{\prime }\left( y\right) \right) =\delta
\left( n\left( x\right) \right) \delta \left( n^{\prime }\left( y\right)
\right) =
$$
$$
m\left( x,0\right) +m(0,y)=m\left( x,y\right) ,
$$
we prove the first identity in (1.46).

On the other hand, d-algebra $C\left( h{\Cal F},g\right) \widehat{\otimes }%
C\left( v{\Cal F},h\right) $ is generated by elements of type $n\left(
x\right) \otimes 1$ and $1\otimes n^{\prime }\left( y\right) ,$ we prove the
second identity in (1.46).

Following from the above--mentioned properties of homomorphisms $\zeta $ and
$\upsilon $ we can assert that the natural isomorphism is explicitly
constructed.$\blacksquare$
\vskip10pt

In consequence of theorem 1.2 we conclude that all operations with Clifford
d--algebras can be reduced to calculations for
$C\left( h{\Cal F},g\right) $
and $C\left( v{\Cal F},h\right) $ which are usual Clifford algebras of
dimension $2^n$ and, respectively, $2^m$
\cite{Atiyah, Bott and Shapiro 1964} and \cite{Karoubi 1978}.

Of special interest is the case when $k={\Cal R}$ and ${\Cal F}$ is
isomorphic to vector space ${\Bbb R}^{p+q,a+b}$ provided with quadratic form
$-x_1^2-...-x_p^2+x_{p+q}^2-y_1^2-...-y_a^2+...+y_{a+b}^2.$ In this case,
the Clifford algebra, denoted as $\left( C^{p,q},C^{a,b}\right) ,\,$ is
generated by symbols $%
e_1^{(x)},e_2^{(x)},...,e_{p+q}^{(x)},e_1^{(y)},e_2^{(y)},...,e_{a+b}^{(y)}$
satisfying properties $\left( e_i\right) ^2=-1~\left( 1\leq i\leq p\right)
,\left( e_j\right) ^2=-1~\left( 1\leq j\leq a\right) ,\left( e_k\right)
^2=1~(p+1\leq k\leq p+q),$

$\left( e_j\right) ^2=1~(n+1\leq s\leq a+b),~e_ie_j=-e_je_i,~i\neq j.\,$
Explicit calculations of $C^{p,q}$ and $C^{a,b}$ are possible by using
isomorphisms \cite{Karoubi 1978} and \cite{Penrose and Rindler 1986}
$$
C^{p+n,q+n}\simeq C^{p,q}\otimes M_2\left( {\Bbb R}\right) \otimes
...\otimes M_2\left( {\Bbb R}\right) \cong C^{p,q}\otimes M_{2^n}\left(
{\Bbb R}\right) \cong M_{2^n}\left( C^{p,q}\right) ,
$$
where $M_s\left( A\right) $ denotes the ring of quadratic matrices of order $%
s$ with coefficients in ring $A.$ Here we write the simplest isomorphisms $%
C^{1,0}\simeq {\Bbb C}, C^{0,1}\simeq {\Bbb R}\oplus \Bbb R ,$ and
 $C^{2,0}={\Bbb H},$ where by ${\Bbb H}$ is denoted the body of quaternions.
 We summarize this calculus as
$$
C^{0,0}={\Bbb R}, C^{1,0}={\Bbb C}, C^{0,1}={\Bbb R}\oplus{ \Bbb R},
C^{2,0}={\Bbb H}, C^{0,2}= M_2\left( {\Bbb R}\right) ,
$$
$$
C^{3,0}={\Bbb H}\oplus {\Bbb H} , C^{0,3} = M_2\left( {\Bbb R}\right),
C^{4,0}=M_2\left( {\Bbb H}\right) , C^{0,4}=M_2\left( {\Bbb H}\right) ,
$$
$$
C^{5,0}=M_4\left( {\Bbb C}\right) ,~C^{0,5}=M_2\left( {\Bbb H}\right) \oplus
M_2\left( {\Bbb H}\right) ,~C^{6,0}=M_8\left( {\Bbb R}\right)
,~C^{0,6}=M_4\left( {\Bbb H}\right) ,
$$
$$
C^{7,0}=M_8\left( {\Bbb R}\right) \oplus M_8\left( {\Bbb R}\right)
,~C^{0,7}=M_8\left( {\Bbb C}\right) ,~C^{8,0}=M_{16}\left( {\Bbb R}\right)
,~C^{0,8}=M_{16}\left( {\Bbb R}\right) .
$$
One of the most important properties of real algebras $C^{0,p}~\left(
C^{0,a}\right) $ and $C^{p,0}~\left( C^{a,0}\right) $ is eightfold
periodicity of $p(a).$

Now, we emphasize that $H^{2n}$-spaces admit locally a structure of Clifford
algebra on complex vector spaces. Really, by using almost \ Hermitian
structure $J_\alpha ^{\quad \beta }$ and considering complex space ${\Bbb C}%
^n$ with nondegenarate quadratic form $\sum_{a=1}^n\left| z_a\right|
^2,~z_a\in {\Bbb C}^2$ induced locally by metric (1.12)
(rewritten in complex coordinates $z_a=x_a+iy_a)$ we define Clifford algebra
$\overleftarrow{C}^n=
\underbrace{\overleftarrow{C}^1\otimes ...\otimes \overleftarrow{C}%
^1}_{n},$
where
$\overleftarrow{C}^1={\Bbb C}{\otimes }_R{\Bbb C}={\Bbb C}\oplus {\Bbb C}$
or in
consequence, $\overleftarrow{C}^n\simeq C^{n,0}\otimes _{{\Bbb R}}{\Bbb C}%
\approx C^{0,n}\otimes _{{\Bbb R}}{\Bbb C}.$ Explicit calculations lead to
isomorphisms $\overleftarrow{C}^2=C^{0,2}\otimes _{{\Bbb R}}{\Bbb C}\approx
M_2\left( {\Bbb R}\right) \otimes _{{\Bbb R}}{\Bbb C}\approx M_2\left(
\overleftarrow{C}^n\right) ,~C^{2p}\approx M_{2^p}\left( {\Bbb C}\right) $
and $\overleftarrow{C}^{2p+1}\approx M_{2^p}\left( {\Bbb C}\right) \oplus
M_{2^p}\left( {\Bbb C}\right) ,$ which show that complex Clifford algebras,
defined locally for $H^{2n}$-spaces, have periodicity 2 on $p.$

Considerations presented in the proof of theorem 1.2 show that map
$j:{\Cal F}\rightarrow C\left( {\Cal F}\right) $ is monomorphic,
so we can identify
space ${\Cal F}$ with its image in $C\left( {\Cal F},G\right) ,$
denoted as $u\rightarrow \overline{u},$ if
$u\in C^{(0)}\left( {\Cal F},
G\right) ~\left( u\in C^{(1)}\left( {\Cal F},G\right) \right) ;$ then $u=%
\overline{u}$ ( respectively, $\overline{u} = -u ) .$

\proclaim{\bf Definition 1.5}
The set of elements $u\in C\left( G\right) ^{*},$ where
$C\left( G\right) ^{*}$ denotes the multiplicative group of invertible
elements of $C\left( {\Cal F},G\right) $ satisfying $\overline{u}{\Cal F}%
u^{-1}\in {\Cal F},$ is called the twisted Clifford d--group, denoted as $%
\widetilde{\Gamma }\left( {\Cal F}\right) .$
\endproclaim

Let $\widetilde{\rho }:\widetilde{\Gamma }\left( {\Cal F}\right) \rightarrow
GL\left( {\Cal F}\right) $ be the homorphism given by $u\rightarrow \rho
\widetilde{u},$ where $\widetilde{\rho }_u\left( w\right) =\overline{u}%
wu^{-1}.$ We can verify that $\ker \widetilde{\rho }={\Bbb R}^{*}$is a
subgroup in $\widetilde{\Gamma }\left( {\Cal F}\right) .$

Canonical map $j:{\Cal F}\rightarrow C\left( {\Cal F}\right) $ can be
interpreted as the linear map ${\Cal F}\rightarrow C\left( {\Cal F}\right)
^0 $ satisfying the universal property of Clifford d-algebras. This leads to
a homomorphism of algebras, $C\left( {\Cal F}\right) \rightarrow C\left(
{\Cal F}\right) ^t,$ considered by an anti--involution of $C\left( {\Cal F}%
\right) $ and denoted as $u\rightarrow ~^tu.$ More exactly, if $u_1...u_n\in
{\Cal F,}$ then $t_u=u_n...u_1$ and $^t\overline{u}=\overline{^tu}=\left(
-1\right) ^nu_n...u_1.$

\proclaim{\bf Definition 1.6}
 The spinor norm of arbitrary $u\in C\left( {\Cal F}%
\right) $ is defined as\newline
  $S\left( u\right) =~^t\overline{u}\cdot u\in C\left(
{\Cal F}\right) .$
\endproclaim

It is obvious that if $u,u^{\prime },u^{\prime \prime }\in \widetilde{\Gamma
}\left( {\Cal F}\right) ,$ then $S(u,u^{\prime })=S\left( u\right) S\left(
u^{\prime }\right) $ and  $S\left( uu^{\prime }u^{\prime \prime }\right)
=S\left( u\right) S\left( u^{\prime }\right) S\left( u^{\prime \prime
}\right) .$ For $u,u^{\prime }\in {\Cal F} S\left( u\right) =-G\left(
u\right) $ and $S\left( u,u^{\prime }\right) =S\left( u\right) S\left(
u^{\prime }\right) =S\left( uu^{\prime }\right) .$

Let us introduce the orthogonal group $O\left( G\right) \subset GL\left(
G\right) $ defined by metric $G$ on ${\Cal F}$ and denote sets $SO\left(
G\right) =\{u\in O\left( G\right) ,\det \left| u\right| =1\},~Pin\left(
G\right) =\{u\in \widetilde{\Gamma }\left( {\Cal F}\right) ,S\left( u\right)
=1\}$ and $Spin\left( G\right) =Pin\left( G\right) \cap C^0\left( {\Cal F}%
\right) .$ For ${{\Cal F}\cong {\Bbb R}}^{n+m}$ we write
$Spin\left( n+m\right) .$ By
straightforward calculations (see similar considerations in
\cite{Karoubi 1978}) we can verify the exactness of these sequences:%
$$
1\rightarrow {\Bbb Z}/2\rightarrow Pin\left( G\right) \rightarrow O\left(
G\right) \rightarrow 1,
$$
$$
1\rightarrow {\Bbb Z}/2\rightarrow Spin\left( G\right) \rightarrow SO\left(
G\right) \rightarrow 0,
$$
$$
1\rightarrow
{\Bbb Z}/2\rightarrow Spin\left( n+m\right) \rightarrow SO\left( n+m\right)
\rightarrow 1.
$$
We conclude this section by emphasizing that the spinor norm was defined
with respect to a quadratic form induced by a metric in v-bundle $\xi _d$
(or by an $H^{2n}$-metric in the case of GL-spaces). This approach differs
from those presented in  \cite{Asanov and Ponomarenko 1988} and
\cite{Ono and Takano 1993}.
\vskip25pt

\head{I.3 Clifford Bundles and Distinguished Spinor Structures}\endhead

There are two possibilities for generalizing our spinor constructions
defined for d-vector spaces to the case of vector bundle spaces enabled with
the structure of N--connection. The firs is to use the extension to the
category of vector bundles. The second is to define the Clifford fibration
associated with compatible linear d-connection and metric $G$ on a vector
bundle (or with an $H^{2n}$-metric on GL-space). Let us consider both
variants.
\vskip15pt

\subhead{I.3.1 Clifford distinguished module structures in vector bundles}
\endsubhead

Because functor ${\Cal F}\to C({\Cal F})$
is smooth we can extend it to the category of
vector bundles of type $\xi _d=\{\pi _d:HE\oplus VE\rightarrow E\}.$ Recall
that by ${\Cal F}$ we denote the typical fiber of such bundles. For $\xi _d$
we obtain a bundle of algebras, denoted as $C\left( \xi _d\right) ,\,$ such
that $C\left( \xi _d\right) _u=C\left( {\Cal F}_u\right) .$ Multiplication
in every fibre defines a continuous map $C\left( \xi _d\right) \times
C\left( \xi _d\right) \rightarrow C\left( \xi _d\right) .$ If $\xi _d$ is a
vector bundle on number field $k$,\thinspace \thinspace the structure of the
$C\left( \xi _d\right) $-module, the d-module, the d-module, on $\xi _d$ is
given by the continuous map $C\left( \xi _d\right) \times _E\xi
_d\rightarrow \xi _d$ with every fiber ${\Cal F}_u$ provided with the
structure of the $C\left( {\Cal F}_u\right)$--module, correlated with its $k$%
--module structure, Because ${\Cal F}\subset C\left( {\Cal F}\right) ,$ we
have a fiber to fiber map ${\Cal F}\times _E\xi _d\rightarrow \xi _d,$
inducing on every fiber the map ${\Cal F}_u\times _E\xi _{d\left( u\right)
}\rightarrow \xi _{d\left( u\right) }$ (${\Bbb R}$--linear on the first
factor and $k$-linear on the second one ). Inversely, every such bilinear
map defines on $\xi _d$ the structure of the $C\left( \xi _d\right) $-module
by virtue of universal properties of Clifford d--algebras. Equivalently, the
above-mentioned bilinear map defines a morphism of v--bundles $m:\xi
_d\rightarrow HOM\left( \xi _d,\xi _d\right) \quad [HOM\left( \xi _d,\xi
_d\right) $ denotes the bundles of homomorphisms] when $\left( m\left(
u\right) \right) ^2=G\left( u\right) $ on every point.

Vector bundles $\xi _d$ provided with $C\left( \xi _d\right) $--structures
are objects of the category with morphisms being morphisms of v--bundles,
which induce on every point $u\in \xi $ morphisms of $C\left( {\Cal F}%
_u\right)$--modules. This is a Banach category contained in the category of
finite-dimensional d-vector spaces on filed $k.$ We shall not use category
formalism in this work, but point to its advantages in further formulation
of new directions of K--theory (see , for example, an introduction in
\cite{Karoubi 1978}) concerned with la--spaces.

Let us denote by $H^s\left( \xi ,GL_{n+m}\left( {\Bbb R}\right) \right) \,$
the s-dimensional cohomology group of the algebraic sheaf of germs of
continuous maps of v--bundle $\xi $ with group $GL_{n+m}\left( {\Bbb R}\right)
$ the group of automorphisms of ${\Bbb R}^{n+m}\,$ (for the language of
algebraic topology see, for example,
\cite{Karoubi 1978} and \cite{Godbillon 1971}). We shall also use the
group $SL_{n+m}\left( {\Bbb R}\right) =\{A\subset GL_{n+m}\left( {\Bbb R}%
\right) ,\det A=1\}.$ Here we point out that cohomologies $H^s(M,Gr)$
characterize the class of a principal bundle $\pi :P\rightarrow M$ on $M$
with structural group $Gr.$ Taking into account that we deal with bundles
distinguished by an N--connection we introduce into consideration
cohomologies $H^s\left( \xi ,GL_{n+m}\left( {\Bbb R}\right) \right) $ as
distinguished classes (d--classes) of bundles $\xi $ provided with a global
N--connection structure.

For a real vector bundle $\xi _d$ on compact base $\xi $ we can define the
orientation on $\xi _d$ as an element $\alpha _d\in H^1\left( \xi
,GL_{n+m}\left( {\Bbb R}\right) \right) $ whose image on map%
$$
H^1\left( \xi ,SL_{n+m}\left( {\Bbb R}\right) \right) \rightarrow H^1\left(
\xi ,GL_{n+m}\left( {\Bbb R}\right) \right)
$$
is the d--class of bundle $\xi .$

\proclaim{\bf Definition 1.7}
The spinor structure on $\xi _d$ is defined as an
element $\beta _d\in H^1\left( \xi ,Spin\left( n+m\right) \right) $ whose
image in the composition%
$$
H^1\left( \xi ,Spin\left( n+m\right) \right) \rightarrow H^1\left( \xi
,SO\left( n+m\right) \right) \rightarrow H^1\left( \xi ,GL_{n+m}\left( {\Bbb %
R}\right) \right)
$$
is the d--class of $\xi .$
\endproclaim

The above definition of spinor structures can be reformulated in terms of
principal bundles. Let $\xi _d$ be a real vector bundle of rank n+m on a
compact base $\xi .$ If there is a principal bundle $P_d$ with structural
group $SO\left( n+m\right) ~[$ or $Spin\left( n+m\right) ],$ this bundle $%
\xi _d$ can be provided with orientation (or spinor) structure. The bundle $%
P_d$ is associated with element $\alpha _d\in H^1\left( \xi ,SO(n+m)\right) $
[or $\beta _d\in H^1\left( \xi ,Spin\left( n+m\right) \right) .$

We remark that a real bundle is oriented if and only if its first
Stiefel--Whitney d--class vanishes,
$$
w_1\left( \xi _d\right) \in H^1\left( \xi ,{\Bbb Z}/2\right) =0,
$$
where $H^1\left( \xi ,{\Bbb Z}/2\right) $ is the first group of Chech
cohomology with coefficients in ${\Bbb Z}/2 .$ Considering the second
Stiefel-Whitney class $w_2\left( \xi _d\right) \in H^{21}\left( \xi ,{\Bbb Z}%
/2\right) $ it is well known that vector bundle $\xi _d$ admits the
spinor structure if and only if $w_2\left( \xi _d\right) =0.$ Finally, in
this subsection, we emphasize that taking into account that base space $\xi $
is also a v-bundle, $p:E\rightarrow M,$ we have to make explicit
calculations in order to express cohomologies $H^s\left( \xi
,GL_{n+m}\right) \,$ and $^{}H^s\left( \xi ,SO\left( n+m\right) \right) $
through cohomologies $H^s\left( M,GL_n\right) ,H^s\left( M,SO\left( m\right)
\right) $ , which depends on global topological structures of spaces $M$ and
$\xi .$ For general bundle and base spaces this requires a cumbersome
cohomological calculus.
\vskip15pt

\subhead{I.3.2 Clifford fibration} \endsubhead

Another way of defining the spinor structure is to use Clifford fibrations.
Consider the principal bundle with the structural group $Gr$ being a
subgroup of orthogonal group $O\left( G\right) ,$ where $G$ is a quadratic
nondegenerate form (see(1.12)) defined on the base (also being a bundle space)
space $\xi .$  The fibration associated
to principal fibration $P\left( \xi ,Gr\right) \left[ \text{or }P\left(
H^{2n},Gr\right) \right] $ with a typical fiber having Clifford algebra $%
C\left( G\right) $ is, by definition, the Clifford fibration $PC\left( \xi
,Gr\right) .$ We can
always define a metric on the Clifford fibration if every fiber is isometric
to $PC\left( \xi ,G\right) $ (this result is proved for arbitrary quadratic
forms $G$ on pseudo-Riemannian bases \cite{Turtoi 1989}). If,
additionally, $Gr\subset SO\left( G\right) $ a global section can be defined
on $PC\left( G\right) .$

Let ${\Cal P}\left( \xi ,Gr\right) $ be the set of principal bundles with
differentiable base $\xi $ and structural group $Gr.$ If
$g:Gr\rightarrow Gr^{\prime }$ is an homomorphism of Lie groups and $P\left(
\xi ,Gr\right) \subset {\Cal P}\left( \xi ,Gr\right)  $
(for simplicity in this section we shall
denote mentioned bundles and sets of bundles as $P,P^{\prime }$ and
respectively, ${\Cal P},{\Cal P}^{\prime }),$
we can always construct a principal
bundle with the property that there is as homomorphism $f:P^{\prime
}\rightarrow P$ of principal bundles which can be projected to the identity
map of $\xi $  and corresponds to isomorphism $g:Gr\rightarrow
Gr^{\prime }.$ If the inverse statement also holds, the bundle $P^{\prime }$
is called as the extension of $P$ associated to $g$ and $f$ is called the
extension homomorphism denoted as $\widetilde{g.}$

Now we can define distinguished spinor structures on bundle spaces (compare
with definition 1.7 ).

\proclaim{\bf Definition 1.8}
Let $P\in {\Cal P}\left( \xi ,O\left( G\right) \right)
$ be a principal bundle. A distinguished spinor structure of $P,$
equivalently a ds-structure of $\xi $ is an extension $%
\widetilde{P}$ of $P$ associated to homomorphism $h:PinG\rightarrow O\left(
G\right) $ where $O\left( G\right) $ is the group of orthogonal rotations,
generated by metric $G,$ in bundle $\xi .$
\endproclaim

So, if $\widetilde{P}$ is a spinor structure of the space $\xi ,$
 then $\widetilde{P}\in {\Cal P}\left( \xi ,PinG\right) .$

On the definition of spinor structures on varieties we cite
\cite{Geroch 1958}. It is proved that a necessary and sufficient
condition for a space time to be orientable is to admit a global field of
orthonormalized frames. We mention that spinor structures can be also
defined on varieties modeled on Banach spaces \cite{Anastasiei 1977}.
 As we have shown in this subsection, similar constructions are
possible for the cases when space time has the structure of a v--bundle with
an N--connection.

\proclaim{\bf Definition 1.9}
 A special distinguished spinor structure, ds-structure,
of principal bundle $P=P\left( \xi ,SO\left( G\right) \right) $
 is a principal bundle $%
\widetilde{P}=\widetilde{P}\left( \xi ,SpinG\right) $
 for which a homomorphism of
principal bundles $\widetilde{p}:\widetilde{P}\rightarrow P,$ projected on
the identity map of $\xi $ (or of $H^{2n})$ and corresponding to
representation%
$$
R:SpinG\rightarrow SO\left( G\right) ,
$$
is defined.
\endproclaim

In the case when the base space variety is oriented, there is a natural
bijection between tangent spinor structures with a common base. For special
ds--structures we can define, as for any spinor structure, the concepts of
spin tensors, spinor connections, and spinor covariant derivations.
\vskip25pt

\head{I.4 Almost Complex Spinor Structures} \endhead

Almost complex structures are an important characteristic of $H^{2n}$%
-spaces.  We can rewrite the almost Hermitian metric
 $H^{2n}$-metric (see
considerations from subsection 1.1.1 with respect to conditions of type (1.10)
 and (1.16) for a metric (1.12)),
  in complex form:
$$
G=H_{ab}\left( z,\xi \right) dz^a\otimes dz^b,\tag 1.47
$$
where
$$
z^a=x^a+iy^a,~\overline{z^a}=x^a-iy^a,~H_{ab}\left( z,\overline{z}\right)
=g_{ab}\left( x,y\right) \mid _{y=y\left( z,\overline{z}\right) }^{x=x\left(
z,\overline{z}\right) },
$$
and define almost complex spinor structures. For given metric (1.47) on $%
H^{2n}$-space there is always a principal bundle $P^U$ with unitary structural
group U(n) which allows us to transform $H^{2n}$-space
into v-bundle $\xi ^U\approx P^U\times _{U\left( n\right) }%
{\Bbb R}^{2n}.$ This statement will be proved after we introduce complex
%\begin{figure}[htbp]
%\begin{center}
%\begin{picture}(255,50) \setlength{\unitlength}{1pt}
%\thinlines
%\put(0,45){$U(n)$}
%\put(116,45){$SO(2n)$}
%\put(53,2){${Spin}^c (2n)$}
%\put(70,38){i}
%\put(70,48){ \line(-1,0){45}}
%\put(70,48){\vector(1,0){45}}
%\put(25,42){\vector(3,-2){45}}
%\put(78,13) {\vector(3,2){45}}
%\put(245,10){(1.48)}
%\put(118,42){\vector(-3,-2){45}}
%\put(35,20){$\sigma$}
%\put(100,20){${\rho}^c$}
%\end{picture}
%\end{center}
%\caption{Diagram 3}
%\end{figure}
spinor structures on oriented real vector bundles \cite{Karoubi 1978}.

Let us consider momentarily $k={\Cal C}$ and introduce into consideration
[instead of the group $Spin(n)]$ the group $Spin^c\times _{{\Bbb Z}/2}U\left(
1\right) $ being the factor group of the product $Spin(n)\times U\left(
1\right) $ with the respect to equivalence%
$$
\left( y,z\right) \sim \left( -y,-a\right) ,\quad y\in Spin(m).
$$
This way we define the short exact sequence%
$$
1\rightarrow U\left( 1\right) \rightarrow Spin^c\left( n\right)
{\overset{\rho ^c} \to \longrightarrow} SO\left( n\right) \rightarrow 1,
$$
where $\rho ^c\left( y,a\right) =\rho ^c\left( y\right) .$ If $\lambda $ is
oriented , real, and rank $n,$ $\gamma $-bundle $\pi :E_\lambda \rightarrow
M^n,$ with base $M^n,$ the complex spinor structure, spin structure, on $%
\lambda $ is given by the principal bundle $P$ with structural group $%
Spin^c\left( m\right) $ and isomorphism $\lambda \approx P\times
_{Spin^c\left( n\right) }{\Bbb R}^n.$ For such bundles the categorial
equivalence can be defined as
$$
\epsilon ^c:{\Cal E}_{{\Bbb C}}^T\left( M^n\right) \rightarrow {\Cal E}_{%
{\Bbb C}}^\lambda \left( M^n\right) , \tag 1.48
$$
where $\epsilon ^c\left( E^c\right) =P\bigtriangleup _{Spin^c\left( n\right)
}E^c$ is the category of trivial complex bundles on $M^n,{\Cal E}_{{\Bbb C}%
}^\lambda \left( M^n\right) $ is the category of complex v-bundles on $M^n$
with action of Clifford bundle $C\left( \lambda \right) ,P\bigtriangleup
_{Spin^c(n)}$ and $E^c$ is the factor space of the bundle product $P\times
_ME^c$ with respect to the equivalence $\left( p,e\right) \sim \left( p%
\widehat{g}^{-1},\widehat{g}e\right) ,p\in P,e\in E^c,$ where $\widehat{g}%
\in Spin^c\left( n\right) $ acts on $E$ by via the imbedding $Spin\left(
n\right) \subset C^{0,n}$ and the natural action $U\left( 1\right) \subset
{\Cal C}$ on complex v-bundle $\xi ^c,E^c=tot\xi ^c,$ for bundle $\pi
^c:E^c\rightarrow M^n.$

Now we return to the bundle $\xi .$ A real v-bundle (not being a spinor
bundle) admits a complex spinor structure if and only if there exist a
homomorphism $\sigma :U\left( n\right) \rightarrow Spin^c\left( 2n\right) $
making the diagram
$$\CD
{U(n)} @. \overset \sigma \to {\longrightarrow} Spin^c (2n)\\
@V{\varphi}VV      \qquad  @VV{\rho^c}V\\
SO(2n)  @. \Longleftrightarrow SO(2n)
\endCD \tag 1.49$$
commutative. The explicit construction of $\sigma $ for arbitrary $\gamma $%
--bundle is given in
\cite{Karoubi 1978}  and \cite{Atiyah, Bott and Shapiro 1964}.
For $H^{2n}$-spaces it is obvious that a
diagram similar to (1.49) can be defined for the tangent bundle $TM^n ,$ which
directly points to the possibility of defining the $^cSpin$-structure on $%
H^{2n}$-spaces.

Let $\lambda $ be a complex, rank\thinspace $n,$ spinor bundle with
$$
\tau :Spin^c\left( n\right) \times _{{\Bbb Z}/2}U\left( 1\right) \rightarrow
U\left( 1\right) \tag 1.50
$$
the homomorphism defined by formula $\tau \left( \lambda ,\delta \right)
=\delta ^2.$ For $P_s$ being the principal bundle with fiber $Spin^c\left(
n\right) $ we introduce the complex linear bundle $L\left( \lambda ^c\right)
=P_S\times _{Spin^c(n)}{\Bbb C}$ defined as the factor space of $P_S\times
{\Bbb C}$ on equivalence relation$$\left( pt,z\right) \sim
\left( p,l\left(
t\right) ^{-1}z\right) ,$$ where $t\in Spin^c\left( n\right) .$ This linear
bundle is associated to complex spinor structure on $\lambda ^c.$

If $\lambda ^c$ and $\lambda ^{c^{\prime }}$ are complex spinor bundles, the
Whitney sum $\lambda ^c\oplus \lambda ^{c^{\prime }}$ is naturally provided
with the structure of the complex spinor bundle. This follows from the
holomorphism%
$$
\omega ^{\prime }:Spin^c\left( n\right) \times Spin^c\left( n^{\prime
}\right) \rightarrow Spin^c\left( n+n^{\prime }\right) ,\tag 1.51
$$
given by formula $\left[ \left( \beta ,z\right) ,\left( \beta ^{\prime
},z^{\prime }\right) \right] \rightarrow \left[ \omega \left( \beta ,\beta
^{\prime }\right) ,zz^{\prime }\right] ,$ where $\omega $ is the
homomorphism making the following diagram  commutative:%
%\begin{figure}[htbp]
%\begin{center}
%\begin{picture}(190,50) \setlength{\unitlength}{1pt}
%\thinlines
%\put(0,45){$Spin(n)\times Spin(n')$}
%\put(160,45){$Spin(n+n')$}
%\put(15,2){$O(n)\times O(n')$}
%\put(168,02){$O(n+n')$}
%\put(100,45){ \vector(1,0){55}}
%\put(48,40){\vector(0,-1){25}}
%\put(90,5){\vector(1,0){70}}
%\put(185,40) {\vector(0,-1){25}}
%\end{picture}
%\end{center}
%\caption{Diagram 4}
%\end{figure}
$$\CD
Spin(n) \times Spin(n') @>{}>> Spin(n+n') \\
@V{}VV                      @VV{}V\\
O(n) \times O(n') @>{}>> O(n+n')
\endCD $$
Here, $z,z^{\prime }\in U\left( 1\right) .$ It is obvious that $L\left(
\lambda ^c\oplus \lambda ^{c^{\prime }}\right) $ is isomorphic to $L\left(
\lambda ^c\right) \otimes L\left( \lambda ^{c^{\prime }}\right) .$

We conclude this section by formulating our main result on complex spinor
structures for $H^{2n}$-spaces:
\proclaim{\bf Theorem 1.3}
Let $\lambda ^c$ be a complex spinor bundle of rank $n$ and $%
H^{2n}$-space considered as a real vector bundle $\lambda ^c\oplus \lambda
^{c^{\prime }}$ provided with almost complex structure $J_{\quad \beta
}^\alpha ;$ multiplication on $i$ is given by $\left(
\smallmatrix
0 & -\delta _j^i \\
\delta _j^i & 0
\endsmallmatrix \right) $.
Then, the diagram
$$\CD
{\Cal E}^{0,2n}_{\Bbb C} (M^{2n}) \qquad @.
	\overset {\varepsilon}^c \to {\longrightarrow}\qquad
	{\Cal E}^{{\lambda}^c \oplus {\lambda}^c} (M^n) \\
@V{\tilde{\varepsilon}}^cVV      \qquad  @VV{\Cal H}V\\
{\Cal E}^W_{\Bbb C}(M^n)  @. \Longleftrightarrow \qquad
{\Cal E}^W_{\Bbb C}(M^n)
\endCD $$
 is commutative up to isomorphisms
 $\epsilon ^c$ and $\widetilde{\epsilon }^c$  defined as in (1.48), $%
{\Cal H}$ is functor $E^c\rightarrow E^c\otimes L\left( \lambda ^c\right) $
and ${\Cal E}_{{\Bbb C}}^{0,2n}\left( M^n\right) $ is defined by functor $%
{\Cal E}_{{\Bbb C}}\left( M^n\right) \rightarrow {\Cal E}_{{\Bbb C}%
}^{0,2n}\left( M^n\right) $ given as correspondence $E^c\rightarrow \Lambda
\left( {\Bbb C}^n\right) \otimes E^c$ (which is a categorial equivalence), $%
\Lambda \left( {\Bbb C}^n\right) $ is the exterior algebra on ${\Bbb C}^n.$ $%
W$ is the real bundle $\lambda ^c\oplus \lambda ^{c^{\prime }}$ provided
with complex structure.
\endproclaim

{\bf Proof:} We use composition of homomorphisms%
$$
\mu :Spin^c\left( 2n\right) \overset{\pi}\to \longrightarrow SO\left( n\right)
\overset{r}\to \longrightarrow
U\left( n\right) \overset{\sigma}\to \longrightarrow
 Spin^c\left( 2n\right)
 \times _{{\Bbb Z}/2}U\left( 1\right) ,
$$
commutative diagram
$$\CD
Spin(2n) @. \subset \qquad Spin^c(2n) \\
@A{\beta}AA                    @AA{}A\\
SO(n) \qquad   @>{}>> SO(2n)
\endCD $$
and introduce composition of homomorphisms%
$$
\mu :Spin^c\left( n\right) \overset{\Delta}\to \longrightarrow
	 Spin^c\left( n\right) \times
Spin^c\left( n\right) \overset{{\omega}^c}\to \longrightarrow
	 Spin^c\left( n\right) ,
$$
where $\Delta $ is the diagonal homomorphism and $\omega ^c$ is defined as
in (1.51). Using homomorphisms (1.50) and (1.51) we obtain formula $\mu \left(
t\right) =\mu \left( t\right) r\left( t\right) .$

Now consider bundle $P\times _{Spin^c\left( n\right) }Spin^c\left( 2n\right)
$ as the principal $Spin^c\left( 2n\right) $-bundle, associated to $M\oplus
M $ being the factor space of the product $P\times Spin^c\left( 2n\right) $
on the equivalence relation $\left( p,t,h\right) \sim \left( p,\mu \left(
t\right) ^{-1}h\right) .$ In this case the categorial equivalence (1.49) can
be rewritten as
$$
\epsilon ^c\left( E^c\right) =P\times _{Spin^c\left( n\right) }Spin^c\left(
2n\right) \Delta _{Spin^c\left( 2n\right) }E^c
$$
and seen as factor space of $P\times Spin^c\left( 2n\right) \times _ME^c$ on
equivalence relation $$\left( pt,h,e\right) \sim \left( p,\mu \left( t\right)
^{-1}h,e\right) \text{and} \left( p,h_1,h_2,e\right) \sim \left(
p,h_1,h_2^{-1}e\right)$$ (projections of elements $p$ and $e$ coincides on
base $M).$ Every element of $\epsilon ^c\left( E^c\right) $ can be
represented as $P\Delta _{Spin^c\left( n\right) }E^c,$ i.e., as a factor
space $P\Delta E^c$ on equivalence relation $\left( pt,e\right) \sim \left(
p,\mu ^c\left( t\right) ,e\right) ,$ when $t\in Spin^c\left( n\right) .$
The complex line bundle $L\left( \lambda ^c\right) $ can be interpreted as
the factor space of  $P\times _{Spin^c\left( n\right) }{\Bbb C}$ on
equivalence relation $\left( pt,\delta \right) \sim \left( p,r\left(
t\right) ^{-1}\delta \right) .$

Putting $\left( p,e\right) \otimes \left( p,\delta \right) \left( p,\delta
e\right) $ we introduce morphism%
$$
\epsilon ^c\left( E\right) \times L\left( \lambda ^c\right) \rightarrow
\epsilon ^c\left( \lambda ^c\right)
$$
with properties $\left( pt,e\right) \otimes \left( pt,\delta \right)
\rightarrow \left( pt,\delta e\right) =\left( p,\mu ^c\left( t\right)
^{-1}\delta e\right) ,$

$\left( p,\mu ^c\left( t\right) ^{-1}e\right) \otimes \left( p,l\left(
t\right) ^{-1}e\right) \rightarrow \left( p,\mu ^c\left( t\right) r\left(
t\right) ^{-1}\delta e\right) $ pointing to the fact that we have defined
the isomorphism correctly and that it is an isomorphism on every fiber.
$\blacksquare$
\vskip25pt

\head{I.5 Spinor Techniques for Distinguished Vector Spaces}\endhead

The purpose of this section is to show how a corresponding abstract spinor
technique entailing notational and calculations advantages can be developed
for arbitrary splits of dimensions of a d-vector space
${\Cal F}= h{\Cal F} \oplus v{\Cal F}$,
 where $\dim h{\Cal F} = n$ and $\dim v{\Cal F} = m.$ For convenience we
shall also present some necessary coordinate expressions.

The problem of a rigorous definition of spinors on la-spaces (la-spinors,
d-spinors) was posed and solved \cite{Vacaru 1996}
(see previous sections 1.2--1.4)
in the framework of the formalism of
Clifford and spinor structures on v-bundles provided with compatible
nonlinear and distinguished connections and metric. We introduced d--spinors
as corresponding objects of the Clifford d--algebra ${\Cal C}\left( {\Cal F},
G\right)$, defined for a d-vector space ${\Cal F}$ in a standard manner
(see, for instance, \cite{Karoubi 1978}) and proved that operations with
${\Cal C}\left({\Cal F},G\right) \ $ can be reduced to
calculations for ${\Cal C}\left(h{\Cal F},g\right) $ and ${\Cal C} \left(
v{\Cal F},h\right) ,$ which are usual Clifford algebras of respective dimensions 2$%
^n$ and 2$^m$ (if it is necessary we can use quadratic forms $g$ and $h$
correspondingly
induced on $h{\Cal F}$ and $v{\Cal F}$ by a metric ${\bold G}$ (1.12)).
 Considering the orthogonal subgroup
 $O{ \left({\bold G}\right) }\subset GL{\left( {\bold G}\right) }$
defined by a metric ${\bold G}$ we can define the d-spinor norm and
parametrize d-spinors by ordered pairs of elements of Clifford algebras $%
{\Cal C}\left(h{\Cal F},g\right) $ and
${\Cal C}\left( v{\Cal F},h\right) .$ We emphasize
that the splitting of a Clifford d-algebra associated to a v-bundle ${\Cal E}
$ is a straightforward consequence of the global decomposition (1.3) defining
a N-connection structure in ${\Cal E}.$

In this section, as a rule, we shall omit proofs which in most cases are
mechanical but rather tedious. We can apply the methods developed in
\cite{Penrose and Rindler 1984, 1986} and \cite{Luehr and Rosenbaum 1974}
  in a straightforward manner on h- and v-subbundles in order to verify
the correctness of affirmations.
\vskip15pt

\subhead{I.5.1 Clifford d-algebra, d-spinors and d-twistors}\endsubhead

In order to relate the succeeding constructions with Clifford d-algebras
we consider a la--frame decomposition of the metric (1.12):%
$$
G_{\alpha \beta }\left( u\right) =l_\alpha ^{\widehat{\alpha }}\left(
u\right) l_\beta ^{\widehat{\beta }}\left( u\right) G_{\widehat{\alpha }
\widehat{\beta }}, \tag 1.52
$$
where the frame d-vectors and constant metric matrices are
respectively distinguished as
$$
l_\alpha ^{\widehat{\alpha }}\left( u\right) =\left(
\matrix
l_j^{\widehat{j}}\left( u\right) & 0 \\
0 & l_a^{\widehat{a}}\left( u\right)
\endmatrix
\right) \ \text { and }\ G_{\widehat{\alpha }\widehat{\beta }}\left(
\matrix
g_{\widehat{i}\widehat{j}} & 0 \\
0 & h_{\widehat{a}\widehat{b}}
\endmatrix
\right) , \tag 1.53
$$
where
$g_{\widehat{i}\widehat{j}}$ and $h_{\widehat{a}\widehat{b}}$ are diagonal
matrices with $g_{\widehat{i}\widehat{i}}=$ $h_{\widehat{a}\widehat{a}}=\pm
1. $

To generate Clifford d-algebras we start with matrix equations%
$$
\sigma _{\widehat{\alpha }}\sigma _{\widehat{\beta }}+\sigma _{\widehat{%
\beta }}\sigma _{\widehat{\alpha }}=-G_{\widehat{\alpha }\widehat{\beta }%
}I, \tag 1.54
$$
where $I$ is the identity matrix, matrices $\sigma _{\widehat{\alpha }%
}\,(\sigma $-objects) act on a d-vector space ${\Cal F} =
h{\Cal F} \oplus v{\Cal F}$ and theirs components are distinguished as
$$
\sigma _{\widehat{\alpha }}\,=\left\{ (\sigma _{\widehat{\alpha }})_{%
\underline{\beta }}^{\cdot \underline{\gamma }}=\left(
\matrix
(\sigma _{\widehat{i}})_{\underline{j}}^{\cdot \underline{k}} & 0 \\
0 & (\sigma _{\widehat{a}})_{\underline{b}}^{\cdot \underline{c}}
\endmatrix
\right) \right\} , \tag 1.55
$$
indices $\underline{\beta }, \underline{\gamma },...$ refer to spin spaces
of type ${\Cal S} = S_{(h)}\oplus S_{(v)}$ and underlined Latin indices
$\underline{j}, \underline{k},...$ and $\underline{b}, \underline{c},...$
refer respectively to a h-spin space ${\Cal S}_{(h)}$ and a v-spin space $%
{\Cal S}_{(v)},\ $which are correspondingly associated to a h- and
v-decomposition of a v-bundle ${\Cal E}_{(d)} .$ The irreducible algebra of
matrices $\sigma _{\widehat{\alpha }}$ of minimal dimension $N\times N,$
where $N=N_{(n)}+N_{(m)},$ $\dim {\Cal S}_{(h)}$=$N_{(n)}$ and
$\dim {\Cal S}_{(v)}$=$N_{(m)},$ has these dimensions%
$$
N_{(n)}=  \cases 2^{(n-1)/2},&\text{ for $n=2k+1$}\\
             {2^{n/2}},&\text{for $n=2k$} \endcases $$

   and
$$N_{(m)}= \cases 2^{(m-1)/2},& \text{ for $m=2k+1$}\\
             {2^{m/2}},&\text{ for $m=2k$}, \endcases$$
where $k=1,2,...$ .

The Clifford d-algebra is generated by sums on $n+1$ elements of form%
$$
A_1I+B^{\widehat{i}}\sigma _{\widehat{i}}+C^{\widehat{i}\widehat{j}}\sigma _{%
\widehat{i}\widehat{j}}+D^{\widehat{i}\widehat{j}\widehat{k}}\sigma _{%
\widehat{i}\widehat{j}\widehat{k}}+...
$$
and sums of $m+1$ elements of form%
$$
A_2I+B^{\widehat{a}}\sigma _{\widehat{a}}+C^{\widehat{a}\widehat{b}}\sigma _{%
\widehat{a}\widehat{b}}+D^{\widehat{a}\widehat{b}\widehat{c}}\sigma _{%
\widehat{a}\widehat{b}\widehat{c}}+...
$$
with antisymmetric coefficients $C^{\widehat{i}\widehat{j}}=C^{[\widehat{i}
\widehat{j}]},C^{\widehat{a}\widehat{b}}=C^{[\widehat{a}\widehat{b]}},D^{%
\widehat{i}\widehat{j}\widehat{k}}=D^{[\widehat{i}\widehat{j}\widehat{k}%
]},D^{\widehat{a}\widehat{b}\widehat{c}}=D^{[\widehat{a}\widehat{b}\widehat{c%
}]},...$ and matrices $\sigma _{\widehat{i}\widehat{j}}=\sigma _{[\widehat{i}%
}\sigma _{\widehat{j}]},\sigma _{\widehat{a}\widehat{b}}=\sigma _{[\widehat{a%
}}\sigma _{\widehat{b}]},\sigma _{\widehat{i}\widehat{j}\widehat{k}}=\sigma
_{[\widehat{i}}\sigma _{\widehat{j}}\sigma _{\widehat{k}]},...$ . Really, we
have 2$^{n+1}$ coefficients $\left( A_1,C^{\widehat{i}\widehat{j}},D^{%
\widehat{i}\widehat{j}\widehat{k}},...\right) $ and 2$^{m+1}$ coefficients
 $ \left( A_2,C^{\widehat{a}\widehat{b}},D^{\widehat{a}\widehat{b}\widehat{c}%
},...\right) $ of the Clifford algebra on ${\Cal F}.$

For simplicity, in this subsection, we shall present the necessary geometric
constructions only for h-spin spaces ${\Cal S}_{(h)}$ of dimension $N_{(n)}.$
Considerations for a v-spin space ${\Cal S}_{(v)}$ are similar but
with proper characteristics for a dimension $N_{(m)}.$

In order to define the scalar (spinor) product on ${\Cal S}_{(h)}$ we
introduce into consideration this finite sum (because of a finite number of
elements $\sigma _{[\widehat{i}\widehat{j}...\widehat{k}]}$ ):%
$$
^{(\pm )}E_{\underline{k}\underline{m}}^{\underline{i}\underline{j}}=\delta
_{\underline{k}}^{\underline{i}}\delta _{\underline{m}}^{\underline{j}%
}+\frac 2{1!}(\sigma _{\widehat{i}})_{\underline{k}}^{.\underline{i}}(\sigma
^{\widehat{i}})_{\underline{m}}^{.\underline{j}}+\frac{2^2}{2!}(\sigma _{%
\widehat{i}\widehat{j}})_{\underline{k}}^{.\underline{i}}(\sigma ^{\widehat{i%
}\widehat{j}})_{\underline{m}}^{.\underline{j}}+\frac{2^3}{3!}(\sigma _{%
\widehat{i}\widehat{j}\widehat{k}})_{\underline{k}}^{.\underline{i}}(\sigma
^{\widehat{i}\widehat{j}\widehat{k}})_{\underline{m}}^{.\underline{j}}+...
\tag 1.56$$
which can be factorized as
$$
^{(\pm )}E_{\underline{k}\underline{m}}^{\underline{i}\underline{j}}=N_{(n)}
{ }^{(\pm )}\epsilon _{\underline{k}\underline{m}}{ }^{(\pm
)}\epsilon ^{\underline{i}\underline{j}}\text{ for }n=2k \tag 1.57
$$
and%
$$
^{(+)}E_{\underline{k}\underline{m}}^{\underline{i}\underline{j}}=2N_{(n)}
\epsilon _{\underline{k}\underline{m}}\epsilon ^{\underline{i}
\underline{j}},{ }^{(-)}E_{\underline{k}\underline{m}}^{\underline{i}
\underline{j}}=0 \text{ for }n=3(mod4),   \tag 1.58
$$
$$
^{(+)}E_{\underline{k}\underline{m}}^{\underline{i}\underline{j}}=0,
{ }^{(-)}E_{\underline{k}\underline{m}}^{\underline{i}\underline{j}}=2N_{(n)}
\epsilon _{\underline{k}\underline{m}}\epsilon ^{\underline{i}
\underline{j}} \text{ for }n=1(mod4).
$$

Antisymmetry of $\sigma _{\widehat{i}\widehat{j}\widehat{k}...}$ and the
construction of the objects (1.56),(1.57) and (1.58) define the properties
of $
\epsilon $-objects $^{(\pm )}\epsilon _{\underline{k}\underline{m}}$ and $%
\epsilon _{\underline{k}\underline{m}}$ which have an eight--fold periodicity
on $n$ (see details in \cite{Penrose and Rindler 1986} and, with respect to
 la--spaces, \cite{Vacaru 1986}).

For even values of $n$ it is possible the decomposition of every h-spin
space ${\Cal S}_{(h)}$into irreducible h-spin spaces ${\bold S}_{(h)}$
and ${\bold S}_{(h)}^{\prime }$ (one considers splitting
of h-indices, for instance, $\underline{l}=L\oplus L^{\prime }, \underline{m%
}=M\oplus M^{\prime },... ;$ for v-indices we shall write $\underline{a}%
=A\oplus A^{\prime },\underline{b}=B\oplus B^{\prime },...)$ and defines
new $\epsilon $-objects
$$
\epsilon ^{\underline{l}\underline{m}}=\frac 12\left( ^{(+)}\epsilon ^{%
\underline{l}\underline{m}}+^{(-)}\epsilon ^{\underline{l}\underline{m}%
}\right) \text{ and }\widetilde{\epsilon }^{\underline{l}\underline{m}%
}=\frac 12\left( ^{(+)}\epsilon ^{\underline{l}\underline{m}}-^{(-)}\epsilon
^{\underline{l}\underline{m}}\right) \tag 1.59
$$
We shall omit similar formulas for $\epsilon $-objects with lower indices.

We can verify, by using expressions (1.58) and straightforward calculations,
these para\-met\-ri\-za\-ti\-ons on symmetry properties of
$\epsilon $-objects (1.59)
$$\multline
\epsilon ^{\underline{l}\underline{m}}=\left(
\matrix
\epsilon ^{LM}=\epsilon ^{ML} & 0 \\
0 & 0
\endmatrix
\right) \text{ and }
\widetilde{\epsilon }^{\underline{l}\underline{m}}=\left(
\matrix
0 & 0 \\
0 & \widetilde{\epsilon }^{LM}=\widetilde{\epsilon }^{ML}
\endmatrix
\right)\\ \text{ for }n=0(mod8);\endmultline \tag 1.60 $$
$$\multline
\epsilon ^{\underline{l}\underline{m}} =
-\frac 12 {}^{(-)}\epsilon ^{\underline{l}\underline{m}} =
\epsilon ^{\underline{m}\underline{l}},\text{ where }%
^{(+)}\epsilon ^{\underline{l}\underline{m}}=0\text{ and } \widetilde{%
\epsilon }^{\underline{l}\underline{m}}=-\frac 12 {}^{(-)}\epsilon ^{\underline{%
l}\underline{m}}=\widetilde{\epsilon }^{\underline{m}\underline{l}},\\ \text{
for }n=1(mod8);\endmultline$$
$$\multline
\epsilon ^{\underline{l}\underline{m}}=\left(
\matrix
0 & 0 \\
\epsilon ^{L^{\prime }M} & 0
\endmatrix
\right) \text{ and } \widetilde{\epsilon }^{\underline{l}\underline{m}%
}=\left(
\matrix
0 & \widetilde{\epsilon }^{LM^{\prime }}=-\epsilon ^{M^{\prime }L} \\ 0 & 0
\endmatrix
\right) \\ \text{ for }n=2(mod8); \endmultline $$
$$\multline
\epsilon ^{\underline{l}\underline{m}}=-\frac 12 {}^{(+)}\epsilon ^{\underline{l%
}\underline{m}}=-\epsilon ^{\underline{m}\underline{l}},\text{ where }%
^{(-)}\epsilon ^{\underline{l}\underline{m}}=0 \text{ and }
\widetilde{%
\epsilon }^{\underline{l}\underline{m}}=\frac 12 {}^{(+)}\epsilon ^{\underline{l%
}\underline{m}}=-\widetilde{\epsilon }^{\underline{m}\underline{l}},
\\ \text{for }n=3(mod8);\endmultline $$
$$\multline
\epsilon ^{\underline{l}\underline{m}}=\left(
\matrix
\epsilon ^{LM}=-\epsilon ^{ML} & 0 \\
0 & 0
\endmatrix
\right) \text{ and }  \widetilde{\epsilon }^{\underline{l}\underline{m}%
}=\left(
\matrix
0 & 0 \\
0 & \widetilde{\epsilon }^{LM}=-\widetilde{\epsilon }^{ML}
\endmatrix
\right) \\ \text{ for }n=4(mod8);\endmultline$$
$$ \multline
\epsilon ^{\underline{l}\underline{m}}=-\frac 12 {}^{(-)}\epsilon ^{\underline{l%
}\underline{m}}=-\epsilon ^{\underline{m}\underline{l}},\text{ where }%
^{(+)}\epsilon ^{\underline{l}\underline{m}}=0 \text{ and }
\widetilde{%
\epsilon }^{\underline{l}\underline{m}}=-\frac 12
{}^{(-)}\epsilon ^{\underline{%
l}\underline{m}} =-\widetilde{\epsilon }^{\underline{m}\underline{l}},
\\ \text{for }n=5(mod8);\endmultline $$
$$\multline
\epsilon ^{\underline{l}\underline{m}}=\left(
\matrix
0 & 0 \\
\epsilon ^{L^{\prime }M} & 0
\endmatrix
\right) \text{ and } \widetilde{\epsilon }^{\underline{l}\underline{m}%
}=\left(
\matrix
0 & \widetilde{\epsilon }^{LM^{\prime }}=\epsilon ^{M^{\prime }L} \\ 0 & 0
\endmatrix
\right)\\  \text{ for }n=6(mod8);\endmultline$$
$$\multline
\epsilon ^{\underline{l}\underline{m}}=
\frac 12 {}^{(-)}\epsilon ^{\underline{l}%
\underline{m}}=\epsilon ^{\underline{m}\underline{l}},\text{ where }%
{}^{(+)}\epsilon ^{\underline{l}\underline{m}}=0 \text{ and }
 \widetilde{%
\epsilon }^{\underline{l}\underline{m}}=
-\frac 12 {}^{(-)}\epsilon ^{\underline{%
l}\underline{m}}=\widetilde{\epsilon }^{\underline{m}\underline{l}},
\\ \text{for }n=7(mod8). \endmultline
$$

Let denote reduced and irreducible h-spinor spaces in a form pointing to the
symmetry of spinor inner products in dependence of values $n=8k+l$ ($%
k=0,1,2,...;l=1,2,...7)$ of the dimension of the horizontal subbundle (we
shall write respectively $\bigtriangleup $ and $\circ $ for antisymmetric
and symmetric inner products of reduced spinors and $\diamondsuit
=(\bigtriangleup ,\circ )$ and $\widetilde{\diamondsuit }=(\circ
,\bigtriangleup )$ for corresponding parametrizations of inner products, in
brief {\it i.p.}, of irreducible spinors; properties of scalar products of
spinors are defined by $\epsilon $-objects (1.60); we shall use $\blacklozenge $
for a general {\it i.p.} when the symmetry is not pointed out):%
$$
{\Cal S}_{(h)}{ }\left( 8k\right) ={\bold S}_{\circ }\oplus {\bold S}_{\circ
}^{\prime };\tag 1.61$$
$${\Cal S}_{(h)}{ }\left( 8k+1\right) ={\Cal S}_{\circ
}^{(-)}\ \text{({\it i.p.} is defined by an }^{(-)}\epsilon \text{-object);}
$$
$$
{\Cal S}_{(h)} { }\left( 8k+2\right) =\cases
{\Cal S}_{\blacklozenge }=
({\bold S}_{\blacklozenge },{\bold S}_{\blacklozenge }),&\text{ or}
\\ {\Cal S}_{\blacklozenge }^{\prime }=({\bold S}_{\widetilde{\blacklozenge }%
}^{\prime }, {\bold S}_{\widetilde{\blacklozenge }}^{\prime })&{};
\endcases$$
$$ {\Cal S}_{(h)} \left( 8k+3\right) ={\Cal S}_{\bigtriangleup
}^{(+)}\ \text{({\it i.p.} is defined by an }^{(+)}\epsilon \text{-object);}
$$
$$
{\Cal S}_{(h)} \left( 8k+4\right) ={\bold S}_{\bigtriangleup }\oplus
{\bold S}_{\bigtriangleup }^{\prime };$$
 $${\Cal S}_{(h)} \left(
8k+5\right) ={\Cal S}_{\bigtriangleup }^{(-)}\ \text{({\it i.p. }is defined
by an }^{(-)}\epsilon \text{-object),}
$$
$$
{\Cal S}_{(h)} \left( 8k+6\right) =\cases
{\Cal S}_{\blacklozenge } =
({\bold S}_{\blacklozenge },{\bold S}_{\blacklozenge }),&\text{ or}
\\ {\Cal S}_{\blacklozenge }^{\prime }=({\bold S}_{\widetilde{\blacklozenge }%
}^{\prime },{\bold S}_{\widetilde{\blacklozenge }}^{\prime })&{};\endcases
 $$
$$
 {\Cal S}_{(h)} \left( 8k+7\right) ={\Cal S}_{\circ }^{(+)}\ %
\text{({\it i.p. } is defined by an }^{(+)}\epsilon \text{-object)}.
 $$
We note that by using corresponding $\epsilon $-objects we can lower and
rise indices of reduced and irreducible spinors (for $n=2,6(mod4)$ we can
exclude primed indices, or inversely, see details in
\cite{Penrose and Rindler 1986}).

The similar v-spinor spaces are denoted by the same symbols as in (1.61)
provided with a left lower mark ''$|"$ and parametrized with respect to the
values $m=8k^{\prime }+l$ (k'=0,1,...; l=1,2,...,7) of the dimension of the
vertical subbundle, for example, as
$$
{\Cal S}_{(v)} ( 8k^{\prime })={\bold S}_{|\circ }\oplus {\bold S}_{|\circ
}^{\prime },{\Cal S}_{(v)} \left( 8k+1\right) ={\Cal S}_{|\circ
}^{(-)},... \tag 1.62
$$
We use '' $\widetilde {}$ ''-overlined symbols,
$$
{\widetilde {\Cal S}}_{(h)}\left( 8k\right) ={\widetilde{\bold S}}_{\circ
}\oplus \widetilde{S}_{\circ }^{\prime },{\widetilde {\Cal S}}_{(h)}\left(
8k+1\right) ={\widetilde{\Cal S}}_{\circ }^{(-)},... \tag 1.63
$$
and
$$
{\widetilde {\Cal S}}_{(v)} ( 8k^{\prime })={\widetilde{\bold S}}_{|\circ
}\oplus {\widetilde S}_{|\circ }^{\prime },{\widetilde {\Cal S}}_{(v)} \left(
8k^{\prime }+1\right) ={\widetilde {\Cal S}}_{|\circ }^{(-)},...
\tag 1.64
$$
respectively for the dual to (1.61) and (1.62) spinor spaces.

The spinor spaces (1.61)-(1.64) are called the prime spinor spaces, in brief
p-spinors. They are considered as building blocks of distinguished
(n,m)-spinor spaces constructed in this manner:%
$$
{\Cal S}(_{\circ \circ ,\circ \circ })={\bold S}_{\circ }\oplus
{\bold S}_{\circ}^{\prime }\oplus {\bold S}_{|\circ }\oplus
{\bold S}_{|\circ }^{\prime },{\Cal S}(_{\circ
\circ ,\circ }\mid ^{\circ })={\bold S}_{\circ }\oplus
{\bold S}_{\circ }^{\prime}\oplus {\bold S}_{|\circ }\oplus
\widetilde{S}_{|\circ }^{\prime }, \tag 1.65
$$
$$
{\Cal S}(_{\circ \circ ,}\mid ^{\circ \circ })={\bold S}_{\circ }\oplus
{\bold S}_{\circ }^{\prime }\oplus \widetilde{S}_{|\circ }\oplus \widetilde{S}%
_{|\circ }^{\prime },{\Cal S}(_{\circ }\mid ^{\circ \circ \circ })=
{\bold S}_{\circ }\oplus \widetilde{S}_{\circ }^{\prime }\oplus \widetilde{S}%
_{|\circ }\oplus \widetilde{S}_{|\circ }^{\prime },
$$
$$
...............................................
$$
$$
{\Cal S}(_{\triangle }, _{\triangle })={\Cal S}_{\triangle
}^{(+)}\oplus S_{|\bigtriangleup }^{(+)},S ( _{\triangle
},^{\triangle })={\Cal S}_{\triangle }^{(+)}\oplus \widetilde{S}_{|\triangle
}^{(+)},
$$
$$
................................
$$
$$
{\Cal S}(_{\triangle }|^{\circ },_\diamondsuit )={\bold S}_{\triangle
}\oplus \widetilde{S_{\circ }}^{\prime }\oplus {\Cal S}_{|\diamondsuit },%
{\Cal S}(_{\triangle }|^{\circ },^\diamondsuit )={\bold S}_{\triangle
}\oplus \widetilde{{\Cal S}_{\circ }}^{\prime }\oplus { \widetilde{\Cal S}}%
_{|}^\diamondsuit , $$
$$
................................
$$
Considering the operation of dualization of prime components in (1.65) we
can generate different isomorphic variants of distinguished (n,m)-spinor
spaces.

We define a d-spinor space ${\Cal S}_{(n,m)}\ $ as a direct sum of a
horizontal and a vertical spinor spaces of type (1.64), for instance,
$$
{\Cal S}_{(8k,8k^{\prime })} =
{\bold S}_{\circ }\oplus {\bold S}_{\circ }^{\prime
}\oplus {\bold S}_{|\circ }\oplus {\bold S}_{|\circ }^{\prime },
 {\Cal S}_{(8k,8k^{\prime}+1)}\ =
{\bold S}_{\circ }\oplus{\bold S}_{\circ }^{\prime }\oplus {\Cal S}_{|\circ
}^{(-)}, ...,$$
 $$ {\Cal S}_{(8k+4,8k^{\prime }+5)} = {\bold S}_{\triangle
}\oplus {\bold S}_{\triangle }^{\prime }\oplus {\Cal S}_{|\triangle }^{(-)},...
$$
The scalar product on a ${\Cal S}_{(n,m)}\ $ is induced by (corresponding to
fixed values of $n$ and $m$ ) $\epsilon $-objects (1.60) considered for h-
and v-components.

Having introduced d-spinors for dimensions $\left( n,m\right) $ we can
write out the generalization for la-spaces of twistor equations
\cite{Penrose and Rindler 1986} by
using the distinguished $\sigma $-objects (1.55):%
$$
(\sigma _{(\widehat{\alpha }})_{|\underline{\beta }|}^{..\underline{\gamma }}%
\quad \frac{\delta \omega ^{\underline{\beta }}}{\delta u^{\widehat{\beta })}}%
=\frac 1{n+m} \quad G_{\widehat{\alpha }\widehat{\beta }}(\sigma ^{\widehat{%
\epsilon }})_{\underline{\beta }}^{..\underline{\gamma }} \quad
\frac{\delta \omega^{\underline{\beta }}}
{\delta u^{\widehat{\epsilon }}}, \tag1.66
$$
where $\left| \underline{\beta }\right| $ denotes that we do not consider
symmetrization on this index. The general solution of (1.66) on the d-vector
space ${\Cal F}$ looks like as
$$
\omega ^{\underline{\beta }}=\Omega ^{\underline{\beta }}+u^{\widehat{\alpha
}}(\sigma _{\widehat{\alpha }})_{\underline{\epsilon }}^{..\underline{\beta }%
}\Pi ^{\underline{\epsilon }}, \tag1.67
$$
where $\Omega ^{\underline{\beta }}$ and $\Pi ^{\underline{\epsilon }}$ are
constant d-spinors. For fixed values of dimensions $n$ and $m$ we mast
analyze the reduced and irreducible components of h- and v-parts of
equations (1.66) and their solutions (1.67) in order to find the symmetry
properties of a d-twistor ${\bold Z}^\alpha $ defined as a pair of d-spinors%
$$
{\bold Z}^\alpha = (\omega ^{\underline{\alpha }},\pi _{\underline{%
\beta }}^{\prime }),
$$
where $\pi _{\underline{\beta }^{\prime }}=\pi _{\underline{\beta }^{\prime
}}^{(0)}\in {\widetilde{\Cal S}}_{(n,m)}$ is a constant dual d-spinor. The
problem of definition of spinors and twistors on la-spaces was firstly
considered in \cite{Vacaru and Ostaf 1994}
(see also \cite{Vacaru 1987} and \cite{Vacaru and Ostaf 1996b})
in connection with the possibility to
extend the equations (1.66) and theirs solutions (1.67), by using nearly
autoparallel maps, on curved, locally isotropic or anisotropic, spaces.
\vskip15pt

\subhead{I.5.2 Mutual transforms of d-tensors and d-spinors}\endsubhead

The spinor algebra for spaces of higher dimensions can not be considered as
a real alternative to the tensor algebra as for locally isotropic spaces of
dimensions $n=3,4$ \cite{Penrose and Rindler 1984, 1986}.
 The same holds true for la-spaces and we
emphasize that it is not quite convenient to perform a spinor calculus for
dimensions $n,m>>4$. Nevertheless, the concept of spinors is important for
every type of spaces,  we can deeply understand the fundamental
properties of geometical objects on la-spaces, and we shall consider in this
subsection some questions concerning transforms of d-tensor objects into
d-spinor ones.
\vskip10pt

\subsubhead{I.5.2.1 Transformation of d-tensors into d-spinors}
\endsubsubhead

In order to pass from d-tensors to d-spinors we must use $\sigma $-objects
(1.55) written in reduced or irreduced form  \quad (in dependence of fixed values of
dimensions $n$ and $m$ ):
$$
(\sigma _{\widehat{\alpha }})_{\underline{\beta }}^{\cdot \underline{\gamma }%
},~(\sigma ^{\widehat{\alpha }})^{\underline{\beta }\underline{\gamma }%
},~(\sigma ^{\widehat{\alpha }})_{\underline{\beta }\underline{\gamma }%
},...,(\sigma _{\widehat{a}})^{\underline{b}\underline{c}},...,(\sigma _{%
\widehat{i}})_{\underline{j}\underline{k}},...,(\sigma _{\widehat{a}%
})^{AA^{\prime }},...,(\sigma ^{\widehat{i}})_{II^{\prime }},....
\tag1.68
$$
It is obvious that contracting with corresponding $\sigma $-objects (1.68)
we can introduce instead of d-tensors indices the d-spinor ones, for
instance,%
$$
\omega ^{\underline{\beta }\underline{\gamma }}=(\sigma ^{\widehat{\alpha }%
})^{\underline{\beta }\underline{\gamma }}\omega _{\widehat{\alpha }},\quad
\omega _{AB^{\prime }}=(\sigma ^{\widehat{a}})_{AB^{\prime }}\omega _{%
\widehat{a}},\quad ...,\zeta _{\cdot \underline{j}}^{\underline{i}}=(\sigma
^{\widehat{k}})_{\cdot \underline{j}}^{\underline{i}}\zeta _{\widehat{k}%
},....
$$
For d-tensors containing groups of antisymmetric indices there is a more
simple procedure of theirs transforming into d-spinors because the objects
$$
(\sigma _{\widehat{\alpha }\widehat{\beta }...\widehat{\gamma }})^{%
\underline{\delta }\underline{\nu }},\quad (\sigma ^{\widehat{a}\widehat{b}%
...\widehat{c}})^{\underline{d}\underline{e}},\quad ...,(\sigma ^{\widehat{i}%
\widehat{j}...\widehat{k}})_{II^{\prime }},\quad ... \tag1.69
$$
can be used for sets of such indices into pairs of d-spinor indices. Let us
enumerate some properties of $\sigma $-objects of type (1.69) (for
simplicity we consider only h-components having q indices $\widehat{i},%
\widehat{j},\widehat{k},...$ taking values from 1 to $n;$ the properties of
v-components can be written in a similar manner with respect to indices $%
\widehat{a},\widehat{b},\widehat{c}...$ taking values from 1 to $m$):%
$$
(\sigma _{\widehat{i}...\widehat{j}})^{\underline{k}\underline{l}}\text{
 is\ }\cases
\text{symmetric on }\underline{k},\underline{l}&\text{ for }n-2q\equiv
1,7~(mod~8); \\ \text{antisymmetric on }\underline{k},\underline{l}&\text{
for }n-2q\equiv 3,5~(mod~8)
\endcases
 \tag1.70$$
for odd values of $n,$ and an object
$$
(\sigma _{\widehat{i}...\widehat{j}})^{IJ}~\left( (\sigma _{\widehat{i}...%
\widehat{j}})^{I^{\prime }J^{\prime }}\right) \tag1.71$$  $$
  \text{ is\ }\cases
\text{symmetric on }I,J~(I^{\prime },J^{\prime })&\text{ for }n-2q\equiv
0~(mod~8); \\ \text{antisymmetric on }I,J~(I^{\prime },J^{\prime })&\text{
for }n-2q\equiv 4~(mod~8)
\endcases
$$
or%
$$
(\sigma _{\widehat{i}...\widehat{j}})^{IJ^{\prime }}=\pm (\sigma _{\widehat{i%
}...\widehat{j}})^{J^{\prime }I}\cases
n+2q\equiv 6(mod8); \\
n+2q\equiv 2(mod8), \endcases
\tag1.72
$$
with vanishing of the rest of reduced components of the d-tensor $(\sigma _{%
\widehat{i}...\widehat{j}})^{\underline{k}\underline{l}}$ with prime/unprime
sets of indices.

\subsubhead{I.5.2.2 Transformation of d-spinors into d-tensors; fundamental
d--spinors}\endsubsubhead

We can transform every d-spinor $\xi ^{\underline{\alpha }}=\left( \xi ^{%
\underline{i}},\xi ^{\underline{a}}\right) $ into a corresponding d-tensor.
For simplicity, we consider this construction only for a h-component $\xi ^{%
\underline{i}}$ on a h-space being of dimension $n$. The values%
$$
\xi ^{\underline{\alpha }}\xi ^{\underline{\beta }}(\sigma ^{\widehat{i}...%
\widehat{j}})_{\underline{\alpha }\underline{\beta }}\quad \left( n\text{ is
odd}\right) \tag1.73
$$
or
$$
\xi ^I\xi ^J(\sigma ^{\widehat{i}...\widehat{j}})_{IJ}~\left( \text{or }\xi
^{I^{\prime }}\xi ^{J^{\prime }}(\sigma ^{\widehat{i}...\widehat{j}%
})_{I^{\prime }J^{\prime }}\right) ~\left( n\text{ is even}\right)
\tag1.74
$$
with a different number of indices $\widehat{i}...\widehat{j},$ taken
together, defines the h-spinor $\xi ^{\underline{i}}\,$ to an accuracy to the
sign. We emphasize that it is necessary to choose only those h-components
of d-tensors (1.73) (or (1.74)) which are symmetric on pairs of indices $%
\underline{\alpha }\underline{\beta }$ (or $IJ\,$ (or $I^{\prime }J^{\prime }
$ )) and the number $q$ of indices $\widehat{i}...\widehat{j}$ satisfies the
condition (as a respective consequence of the properties (1.70) and/or
(1.71), (1.72))%
$$
n-2q\equiv 0,1,7~(mod~8). \tag1.75
$$
Of special interest is the case when
$$
q=\frac 12\left( n\pm 1\right) ~\left( n\text{ is odd}\right)
\tag1.76
$$
or
$$
q=\frac 12n~\left( n\text{ is even}\right) . \tag1.77
$$
If all expressions (1.73) and/or (1.74) are zero for all values of $q\,$
with the exception of one or two ones defined by the condition (1.76) (or
(1.77)), the value $\xi ^{\widehat{i}}$ (or $\xi ^I$ ($\xi ^{I^{\prime }}))$
is called a fundamental h-spinor. Defining in a similar manner the
fundamental v-spinors we can introduce fundamental d-spinors as pairs of
fundamental h- and v-spinors. Here we remark that a h(v)-spinor $\xi ^{%
\widehat{i}}~(\xi ^{\widehat{a}})\,$ (we can also consider reduced
components) is always a fundamental one for $n(m)<7,$ which is a consequence
of (1.75)).

Finally, in this section, we note that the geometry of fundamental h- and
v-spinors is similar to that of usual fundamental spinors (see  Appendix
to the monograph \cite{Penrose and Rindler 1986}).
We omit such details in this work, but emphasize
that constructions with fundamental d-spinors, for a la-space, must be
adapted to the corresponding global splitting by N-connection of the space.
\vskip25pt

\head{I.6\ The Differential Geometry of Locally Anisotropic Spinors}
\endhead

The goal of the section is to formulate the differential geometry of
d-spinors for la-spaces.

We shall use denotations of type
$$
v^\alpha =(v^i,v^a)\in {{\Cal \sigma}^\alpha } =({{\Cal \sigma}
^i,{\Cal \sigma}^a})\,\text{ and }\zeta ^{\underline{\alpha }}=(\zeta ^{%
\underline{i}},\zeta ^{\underline{a}})\in {{\Cal \sigma}^{\underline{\alpha }}%
}=({{\Cal \sigma} ^{\underline{i}}, {\Cal \sigma}^{\underline{a}}})\,
$$
for, respectively, elements of modules of d-vector and irreduced d-spinor
fields (see details in \cite{Vacaru 1996}). D--tensors and d--spinor tensors
 (irreduced or
reduced) will be interpreted as elements of corresponding ${\Cal \sigma }$
-modules, for instance,
$$
q_{~\beta ...}^\alpha \in {\Cal \sigma ^\alpha ~_{\beta ....}}, \psi
_{~\underline{\beta }\quad ...}^{\underline{\alpha }\quad \underline{\gamma }%
}\in {\Cal \sigma }_{~\underline{\beta }\quad ...}^{\underline{\alpha }\quad
\underline{\gamma }}~, \xi _{\quad JK^{\prime }N^{\prime
}}^{II^{\prime }}\in {\Cal \sigma }_{\quad JK^{\prime }N^{\prime
}}^{II^{\prime }}~,...
$$

We can establish a correspondence between the la-adapted metric $g_{\alpha
\beta }$ (1.12) and d-spinor metric $\epsilon _{\underline{\alpha }%
\underline{\beta }}$ ( $\epsilon $-objects (1.60) for both h- and
v-subspaces of ${\Cal E},$ ) of a la-space ${\Cal E}$ by using the relation%
$$
g_{\alpha \beta }=-\frac 1{N(n)+N(m)}((\sigma _{(\alpha }(u))^{\underline{%
\alpha }_1\underline{\beta }_1}(\sigma _{\beta )}(u))^{\underline{\beta }_2%
\underline{\alpha }_2})\epsilon _{\underline{\alpha }_1\underline{\alpha }%
_2}\epsilon _{\underline{\beta }_1\underline{\beta }_2}, \tag1.78
$$
where%
$$
(\sigma _\alpha (u))^{\underline{\nu }\underline{\gamma }}=l_\alpha ^{%
\widehat{\alpha }}(u)(\sigma _{\widehat{\alpha }})^{\underline{\nu }%
\underline{\gamma }}, \tag1.79
$$
which is a consequence of formulas (1.52)-(1.57). In brief we can write (1.78)
as
$$
g_{\alpha \beta }=\epsilon _{\underline{\alpha }_1\underline{\alpha }%
_2}\epsilon _{\underline{\beta }_1\underline{\beta }_2} \tag1.80
$$
if the $\sigma $-objects are considered as a fixed structure, whereas $%
\epsilon $-objects are treated as caring the metric ''dynamics '' , on
la-space. This variant is used, for instance, in the so-called 2-spinor
geometry \cite{Penrose and Rindler 1984, 1986} and should be preferred if we
 have to make explicit the
algebraic symmetry properties of d-spinor objects. An alternative way is to
considered as fixed the algebraic structure of $\epsilon $-objects and to
use variable components of $\sigma $-objects of type (1.79) for developing a
variational d-spinor approach to gravitational and matter field interactions
on la-spaces ( the spinor Ashtekar variables \cite{Ashtekar, Romano and
 Ranjet 1989} are introduced in this manner).

We note that a d-spinor metric
$$
\epsilon _{\underline{\nu }\underline{\tau }}=\left(
\matrix
\epsilon _{\underline{i}\underline{j}} & 0 \\
0 & \epsilon _{\underline{a}\underline{b}}
\endmatrix
\right)
$$
on the d-spinor space ${\Cal S}=({\Cal S}_{(h)},{\Cal S}_{(v)})$
can have symmetric or
antisymmetric h (v) -components $\epsilon _{\underline{i}\underline{j}}$ ($%
\epsilon _{\underline{a}\underline{b}})$ , see $\epsilon $-objects (1.60).
For simplicity, in this section (in order to avoid cumbersome calculations
connected with eight-fold periodicity on dimensions $n$ and $m$ of a la-space
${\Cal E}$) we shall develop a general d-spinor formalism only by using
irreduced spinor spaces ${\Cal S}_{(h)}$ and ${\Cal S}_{(v)}.$
\vskip15pt

\subhead{I.6.1 D-covariant derivation on la-spaces}\endsubhead

Let ${\Cal E}$ be a la--space. We define the action on a d-spinor of a
d-covariant operator%
$$
\nabla _\alpha =\left( \nabla _i,\nabla _a\right) =(\sigma _\alpha )^{%
\underline{\alpha }_1\underline{\alpha }_2}\nabla _{^{\underline{\alpha }_1%
\underline{\alpha }_2}}=\left( (\sigma _i)^{\underline{i}_1\underline{i}%
_2}\nabla _{^{\underline{i}_1\underline{i}_2}},~(\sigma _a)^{\underline{a}_1%
\underline{a}_2}\nabla _{^{\underline{a}_1\underline{a}_2}}\right)
$$
(in brief, we shall write
$$
\nabla _\alpha =\nabla _{^{\underline{\alpha }_1\underline{\alpha }%
_2}}=\left( \nabla _{^{\underline{i}_1\underline{i}_2}},~\nabla _{^{%
\underline{a}_1\underline{a}_2}}\right)  )
$$
as a map
$$
\nabla _{{\underline{\alpha }}_1 {\underline{\alpha }}_2}\ :\ {\Cal %
\sigma}^{\underline{\beta }}\rightarrow \sigma _\alpha ^{\underline{\beta }}=
\sigma _{{\underline{\alpha }}_1 {\underline{\alpha }}_2}^{\underline{\beta }}
$$
satisfying conditions%
$$
\nabla _\alpha (\xi ^{\underline{\beta }}+\eta ^{\underline{\beta }})=\nabla
_\alpha \xi ^{\underline{\beta }}+\nabla _\alpha \eta ^{\underline{\beta }},
$$
and%
$$
\nabla _\alpha (f\xi ^{\underline{\beta }})=f\nabla _\alpha \xi ^{\underline{%
\beta }}+\xi ^{\underline{\beta }}\nabla _\alpha f
$$
for every $\xi ^{\underline{\beta }},\eta ^{\underline{\beta }}\in {\Cal %
\sigma ^{\underline{\beta }}}$ and $f$ being a scalar field on ${\Cal E}.$
It is also required that one holds the Leibnitz rule%
$$
(\nabla _\alpha \zeta _{\underline{\beta }})\eta ^{\underline{\beta }%
}=\nabla _\alpha (\zeta _{\underline{\beta }}\eta ^{\underline{\beta }%
})-\zeta _{\underline{\beta }}\nabla _\alpha \eta ^{\underline{\beta }}
$$
and that $\nabla _\alpha \,$ is a real operator, i.e. it commuters with the
operation of complex conjugation:%
$$
\overline{\nabla _\alpha \psi _{\underline{\alpha }\underline{\beta }%
\underline{\gamma }...}}=\nabla _\alpha (\overline{\psi }_{\underline{\alpha
}\underline{\beta }\underline{\gamma }...}).
$$

Let now analyze the question on uniqueness of action on d-spinors of an
operator $\nabla _\alpha $ satisfying necessary conditions . Denoting by $%
\nabla _\alpha ^{(1)}$ and $\nabla _\alpha $ two such d-covariant operators
we consider the map%
$$
(\nabla _\alpha ^{(1)}-\nabla _\alpha ):{\Cal \sigma} ^{\underline{\beta }}
\rightarrow {\Cal \sigma}_{\underline{\alpha }_1\underline{\alpha }_2}^{%
\underline{\beta }}. \tag1.81
$$
Because the action on a scalar $f$ of both operators $\nabla _\alpha ^{(1)}$
and $\nabla _\alpha $ must be identical, i.e.%
$$
\nabla _\alpha ^{(1)}f=\nabla _\alpha f, \tag1.82
$$
the action (1.81) on $f=\omega _{\underline{\beta }}\xi ^{\underline{\beta }}$
must be written as
$$
(\nabla _\alpha ^{(1)}-\nabla _\alpha )(\omega _{\underline{\beta }}\xi ^{%
\underline{\beta }})=0.
$$
In consequence we conclude that there is an element $\Theta _{\underline{%
\alpha }_1\underline{\alpha }_2\underline{\beta }}^{\quad \quad \underline{%
\gamma }}\in {\Cal \sigma }_{\underline{\alpha }_1\underline{\alpha }_2%
\underline{\beta }}^{\quad \quad \underline{\gamma }}$ for which%
$$
\nabla _{\underline{\alpha }_1\underline{\alpha }_2}^{(1)}\xi ^{\underline{%
\gamma }}=\nabla _{\underline{\alpha }_1\underline{\alpha }_2}\xi ^{%
\underline{\gamma }}+\Theta _{\underline{\alpha }_1\underline{\alpha }_2%
\underline{\beta }}^{\quad \quad \underline{\gamma }}\xi ^{\underline{\beta }%
} \text{ and }
\nabla _{\underline{\alpha }_1\underline{\alpha }_2}^{(1)}\omega _{%
\underline{\beta }}=\nabla _{\underline{\alpha }_1\underline{\alpha }%
_2}\omega _{\underline{\beta }}-\Theta _{\underline{\alpha }_1\underline{%
\alpha }_2\underline{\beta }}^{\quad \quad \underline{\gamma }}\omega _{%
\underline{\gamma }}~. \tag1.83
$$
The action of the operator (1.81) on a d-vector $v^\beta =v^{\underline{\beta
}_1\underline{\beta }_2}$ can be written by using formula (1.83) for both
indices $\underline{\beta }_1$ and $\underline{\beta }_2$ :%
$$
(\nabla _\alpha ^{(1)}-\nabla _\alpha )v^{\underline{\beta }_1\underline{%
\beta }_2}=\Theta _{\alpha \underline{\gamma }}^{\quad \underline{\beta }%
_1}v^{\underline{\gamma }\underline{\beta }_2}+\Theta _{\alpha \underline{%
\gamma }}^{\quad \underline{\beta }_2}v^{\underline{\beta }_1\underline{%
\gamma }}=
$$
$$
(\Theta _{\alpha \underline{\gamma }_1}^{\quad \underline{\beta }_1}\delta _{%
\underline{\gamma }_2}^{\quad \underline{\beta }_2}+\Theta _{\alpha
\underline{\gamma }_1}^{\quad \underline{\beta }_2}\delta _{\underline{%
\gamma }_2}^{\quad \underline{\beta }_1})v^{\underline{\gamma }_1\underline{%
\gamma }_2}=Q_{\ \alpha \gamma }^\beta v^\gamma ,
$$
where%
$$
Q_{\ \alpha \gamma }^\beta =Q_{\qquad \underline{\alpha }_1\underline{\alpha
}_2~\underline{\gamma }_1\underline{\gamma }_2}^{\underline{\beta }_1%
\underline{\beta }_2}=\Theta _{\alpha \underline{\gamma }_1}^{\quad
\underline{\beta }_1}\delta _{\underline{\gamma }_2}^{\quad \underline{\beta
}_2}+\Theta _{\alpha \underline{\gamma }_1}^{\quad \underline{\beta }%
_2}\delta _{\underline{\gamma }_2}^{\quad \underline{\beta }_1}. \tag1.84
$$
The d-commutator $\nabla _{[\alpha }\nabla _{\beta ]}$ defines the d-torsion
(see (1.27),(1.28) and (1.29)). So, applying operators $\nabla _{[\alpha
}^{(1)}\nabla _{\beta ]}^{(1)}$ and $\nabla _{[\alpha }\nabla _{\beta ]}$ on
$f=\omega _{\underline{\beta }}\xi ^{\underline{\beta }}$ we can write
$$
T_{\quad \alpha \beta }^{(1)\gamma }-T_{~\alpha \beta }^\gamma =Q_{~\beta
\alpha }^\gamma -Q_{~\alpha \beta }^\gamma
$$
with $Q_{~\alpha \beta }^\gamma $ from (1.84).

The action of operator $\nabla _\alpha ^{(1)}$ on d-spinor tensors of type $%
\chi _{\underline{\alpha }_1\underline{\alpha }_2\underline{\alpha }%
_3...}^{\qquad \quad \underline{\beta }_1\underline{\beta }_2...}$ must be
constructed by using formula (1.83) for every upper index $\underline{\beta }%
_1\underline{\beta }_2...$ and formula (1.84) for every lower index $%
\underline{\alpha }_1\underline{\alpha }_2\underline{\alpha }_3...$ .
\vskip15pt

\subhead{I.6.2 Infeld - van der Waerden coefficients and d-con\-nec\-ti\-ons}
\endsubhead

Let $$\delta _{\underline{{\bold \alpha }}}^{\quad \underline{\alpha }}=\left(
\delta _{\underline{{\bold 1}}}^{\quad \underline{i}},
\delta _{\underline{{\bold %
2}}}^{\quad \underline{i}},...,\delta _{\underline{{\bold N}{(n)}}}^{\quad
\underline{i}},\delta _{\underline{{\bold 1}}}^{\quad \underline{a}},\delta _{%
\underline{{\bold 2}}}^{\quad \underline{a}},...,
\delta _{\underline{{\bold N}{(m)}}%
}^{\quad \underline{i}}\right) $$
be a d-spinor basis. The dual to it basis
is denoted as
$$\delta _{\underline{\alpha }}^{\quad \underline{{\bold \alpha }}%
}=\left( \delta _{\underline{i}}^{\quad \underline{{\bold 1}}},\delta _{%
\underline{i}}^{\quad \underline{{\bold 2}}},...,\delta _{\underline{i}%
}^{\quad \underline{{\bold N}{(n)}}},
\delta _{\underline{i}}^{\quad \underline{%
{\bold 1}}},\delta _{\underline{i}}^{\quad \underline{{\bold 2}}},...,
\delta _{%
\underline{i}}^{\quad \underline{{\bold N}{(m)}}}\right) .$$
 A d-spinor $\kappa ^{%
\underline{\alpha }}\in {\Cal \sigma }$ $^{\underline{\alpha }}$ has
components $\kappa ^{\underline{{\bold \alpha }}}=\kappa ^{\underline{\alpha }%
}\delta _{\underline{\alpha }}^{\quad \underline{{\bold \alpha }}}.$ Taking
into account that
$$
\delta _{\underline{{\bold \alpha }}}^{\quad \underline{\alpha }}\delta _{%
\underline{{\bold \beta }}}^{\quad \underline{\beta }}\nabla _{\underline{%
\alpha }\underline{\beta }}=\nabla _{\underline{{\bold \alpha }}\underline{%
{\bold \beta }}},
$$
we write out the components $\nabla _{\underline{\alpha }\underline{\beta }}$
$\kappa ^{\underline{\gamma }}$ as%
$$ \multline
\delta _{\underline{{\bold \alpha }}}^{\quad \underline{\alpha }}~\delta _{%
\underline{{\bold \beta }}}^{\quad \underline{\beta }}~\delta _{\underline{%
\gamma }}^{\quad \underline{{\bold \gamma }}}~\nabla _{\underline{\alpha }%
\underline{\beta }}\kappa ^{\underline{\gamma }}=\\
\delta _{\underline{{\bold %
\epsilon }}}^{\quad \underline{\tau }}~\delta _{\underline{\tau }%
}^{\quad \underline{{\bold \gamma }}}~\nabla _{\underline{{\bold \alpha }}%
\underline{{\bold \beta }}}\kappa ^{\underline{{\bold \epsilon }}}+\kappa ^{%
\underline{{\bold \epsilon }}}~\delta _{\underline{\epsilon }}^{\quad
\underline{{\bold \gamma }}}~\nabla _{\underline{{\bold \alpha }}
\underline{{\bold %
\beta }}}\delta _{\underline{{\bold \epsilon }}}^{\quad \underline{\epsilon }%
}=
\nabla _{\underline{{\bold \alpha }}\underline{{\bold \beta }}}\kappa ^{%
\underline{{\bold \gamma }}}+\kappa ^{\underline{{\bold \epsilon }}}\gamma _{~%
\underline{{\bold \alpha }}\underline{{\bold \beta }}
\underline{{\bold \epsilon }}%
}^{\underline{{\bold \gamma }}}, \endmultline \tag1.85
$$
where the coordinate components of the d-spinor connection $\gamma _{~%
\underline{{\bold \alpha }}\underline{{\bold \beta }}
\underline{{\bold \epsilon }}%
}^{\underline{{\bold \gamma }}}$ are defined as
$$
\gamma _{~\underline{{\bold \alpha }}\underline{{\bold \beta }}
\underline{{\bold %
\epsilon }}}^{\underline{{\bold \gamma }}}\doteq \delta _{\underline{\tau }%
}^{\quad \underline{{\bold \gamma }}}~\nabla _{\underline{{\bold \alpha }}%
\underline{{\bold \beta }}}\delta _{\underline{{\bold \epsilon }}}^{\quad
\underline{\tau }}. \tag1.86
$$
We call the Infeld - van der Waerden d-symbols a set of $\sigma $--objects ($%
\sigma _{{\bold \alpha }})^{\underline{{\bold \alpha }}
\underline{{\bold \beta }}}$
paramet\-ri\-zed with respect to a coordinate d-spinor basis. Defining
$
\nabla _{{\bold \alpha }}=
(\sigma _{{\bold \alpha }})^{\underline{{\bold \alpha }}%
\underline{{\bold \beta }}}~\nabla _{\underline{{\bold \alpha }}
\underline{{\bold\beta }}},
$\
introducing denotations
$
\gamma ^{\underline{{\bold \gamma }}}
{}_{{\bold \alpha \underline{\tau }}}\doteq
\gamma ^{\underline{{\bold \gamma }}}
{}_{\underline{\bold \alpha }\underline{\beta }\underline{\tau }}
(\sigma _{{\bold \alpha }})^{\underline{{\bold \alpha
}}\underline{{\bold \beta }}}
$
and using properties (1.85) we can write relations%
$$
l_{{\bold \alpha }}^\alpha ~\delta _{\underline{\beta }}^{\quad \underline{%
{\bold \beta }}}~\nabla _\alpha \kappa ^{\underline{\beta }}=\nabla _{{\bold %
\alpha }}\kappa ^{\underline{{\bold \beta }}}+\kappa ^{\underline{{\bold \delta }%
}}~\gamma _{~{\bold \alpha }\underline{{\bold \delta }}}^{\underline{{\bold \beta }%
}} \tag1.87
$$
and%
$$
l_{{\bold \alpha }}^\alpha ~\delta _{\underline{{\bold \beta }}}^{\quad
\underline{\beta }}~\nabla _\alpha ~\mu _{\underline{\beta }}=\nabla _{{\bold %
\alpha }}~\mu _{\underline{{\bold \beta }}}-\mu _{\underline{{\bold \delta }}%
}\gamma _{~{\bold \alpha }\underline{{\bold \beta }}}^{\underline{{\bold \delta }}%
} \tag1.88
$$
for d-covariant derivations $~\nabla _{\underline{\alpha }}\kappa ^{%
\underline{\beta }}$ and $\nabla _{\underline{\alpha }}~\mu _{\underline{%
\beta }}.$

We can consider expressions similar to (1.87) and (1.88) for values having
both types of d-spinor and d-tensor indices, for instance,%
$$
l_{{\bold \alpha }}^\alpha ~l_\gamma ^{{\bold \gamma }}
~\delta _{\underline{{\bold %
\delta }}}^{\quad \underline{\delta }}~\nabla _\alpha \theta _{\underline{%
\delta }}^{~\gamma }=\nabla _{{\bold \alpha }}
\theta _{\underline{{\bold \delta }%
}}^{~{\bold \gamma }}-\theta _{\underline{{\bold \epsilon }}}^{~{\bold \gamma }%
}\gamma _{~\underline{{\bold \alpha }}\underline{{\bold \delta }}}^{\underline{%
{\bold \epsilon }}}+
\theta _{\underline{{\bold \delta }}}^{~{\bold \tau }}~\Gamma
_{\quad {\bold \alpha \tau }}^{~{\bold \gamma }}
$$
(we can prove this by a straightforward calculation of the derivation
\newline
 $\nabla _\alpha (\theta _{\underline{{\bold \epsilon }}}^{~{\bold \tau }}$
 $%
~\delta _{\underline{\delta }}^{\quad \underline{{\bold \epsilon }}}~l_{{\bold %
\tau }}^\gamma )).$

Now we shall consider some possible relations between components of
d-connec\-ti\-ons
$\gamma _{~\underline{{\bold \alpha }}\underline{{\bold \delta }}%
}^{\underline{{\bold \epsilon }}}$ and
$\Gamma _{\quad {{\bold \alpha} {\bold \tau }}}^{~%
{\bold \gamma }}$ and derivations of $(\sigma _{{\bold \alpha }})^{\underline{%
{\bold \alpha }}\underline{{\bold \beta }}}$ . According to definitions (1.12)
we can write
$$
\Gamma _{~{\bold \beta}{\bold \gamma }}^{{\bold \alpha }}=
l_\alpha ^{{\bold \alpha }%
}\nabla _{{\bold \gamma }}l_{{\bold \beta }}^\alpha =
l_\alpha ^{{\bold \alpha }%
}\nabla _{{\bold \gamma }}(\sigma _{{\bold \beta }})^{\underline{\epsilon }
\underline{\tau }}=l_\alpha ^{{\bold \alpha }}\nabla _{{\bold \gamma }}((\sigma
_{{\bold \beta }})^{\underline{{\bold \epsilon }}
\underline{{\bold \tau }}}\delta
_{\underline{{\bold \epsilon }}}^{~\underline{\epsilon }}\delta _{\underline{%
{\bold \tau }}}^{~\underline{\tau }})=
$$
$$
l_\alpha ^{{\bold \alpha }}\delta _{\underline{{\bold \alpha }}}^{~\underline{%
\alpha }}\delta _{\underline{{\bold \epsilon }}}^{~\underline{\epsilon }%
}\nabla _{{\bold \gamma }}
(\sigma _{{\bold \beta }})^{\underline{{\bold \alpha }}
\underline{{\bold \epsilon }}}+
l_\alpha ^{{\bold \alpha }}(\sigma _{{\bold \beta }%
})^{\underline{{\bold \epsilon }}\underline{{\bold \tau }}}
(\delta _{\underline{%
{\bold \tau }}}^{~\underline{\tau }}\nabla _{{\bold \gamma }}
\delta _{\underline{%
{\bold \epsilon }}}^{~\underline{\epsilon }}+
\delta _{\underline{{\bold \epsilon
}}}^{~\underline{\epsilon }}\nabla _{{\bold \gamma }}
\delta _{\underline{{\bold\tau }}}^{~\underline{\tau }})=
$$
$$
l_{\underline{{\bold \epsilon }}\underline{{\bold \tau }}}^{{\bold \alpha }%
}~\nabla _{{\bold \gamma }}
(\sigma _{{\bold \beta }})^{\underline{{\bold \epsilon }%
}\underline{{\bold \tau }}}+l_{\underline{{\bold \mu }}
\underline{{\bold \nu }}}^{%
{\bold \alpha }}\delta _{\underline{\epsilon }}^{~\underline{{\bold \mu }}%
}\delta _{\underline{\tau }}^{~\underline{{\bold \nu }}}
(\sigma _{{\bold \beta }%
})^{\underline{\epsilon }\underline{\tau }}(\delta _{\underline{{\bold \tau }}%
}^{~\underline{\tau }}\nabla _{{\bold \gamma }}\delta _{\underline{{\bold %
\epsilon }}}^{~\underline{\epsilon }}+\delta _{\underline{{\bold \epsilon }}%
}^{~\underline{\epsilon }}\nabla _{{\bold \gamma }}\delta _{\underline{{\bold %
\tau }}}^{~\underline{\tau }}),
$$
where $l_{\alpha}^{{\bold \alpha }%
}=(\sigma _{\underline{{\bold \epsilon }}\underline{{\bold \tau }}})^{{\bold %
\alpha }}$ , from which it follows%
$$
(\sigma _{{\bold \alpha }})^{\underline{{\bold \mu }}\underline{{\bold \nu }}%
}(\sigma _{\underline{{\bold \alpha }}
\underline{{\bold \beta }}})^{{\bold \beta }%
}\Gamma _{~{\bold \gamma \beta }}^{{\bold \alpha }}=
(\sigma _{\underline{{\bold %
\alpha }}\underline{{\bold \beta }}})^{{\bold \beta }}\nabla _{{\bold \gamma }%
}(\sigma _{{\bold \alpha }})^{\underline{{\bold \mu }}\underline{{\bold \nu }}%
}+\delta _{\underline{{\bold \beta }}}^{~\underline{{\bold \nu }}}
\gamma _{~{\bold %
\gamma \underline{\alpha }}}^{\underline{{\bold \mu }}}+\delta _{\underline{%
{\bold \alpha }}}^{~\underline{{\bold \mu }}}
\gamma _{~{\bold \gamma \underline{%
\beta }}}^{\underline{{\bold \nu }}}.
$$
Connecting the last expression on $\underline{{\bold \beta }}$ and
 $\underline{
{\bold \nu }}$ and using an orthonormalized d--spinor basis when $\gamma _{~%
{{\bold \gamma} \underline{\bold\beta }}}^{\underline{{\bold \beta }}}=0$ (a
consequence from (1.86)) we have
$$
\gamma _{~{{\bold \gamma}
{ \underline{\bold\alpha }}}}^{\underline{{\bold \mu }}}=\frac
1{N(n)+N(m)}(\Gamma _{\quad {{\bold \gamma}
~\underline{\bold\alpha }\underline{\bold\beta
}}}^{\underline{{\bold \mu }}\underline{{\bold \beta }}}-
(\sigma _{\underline{{\bold \alpha }}
\underline{{\bold \beta }}})^{{\bold \beta }}\nabla _{{\bold \gamma }%
}(\sigma _{{\bold \beta }})^{\underline{{\bold \mu }}\underline{{\bold \beta }}%
}), \tag1.89
$$
where
$$
\Gamma _{\quad {{\bold \gamma} ~\underline{\bold\alpha }
\underline{\bold\beta }}}^{%
\underline{{\bold \mu }}\underline{{\bold \beta }}}=
(\sigma _{{\bold \alpha }})^{%
\underline{{\bold \mu }}\underline{{\bold \beta }}}(\sigma _{\underline{{\bold %
\alpha }}\underline{{\bold \beta }}})^{{\bold \beta }}\Gamma _{~{\bold \gamma
\beta }}^{{\bold \alpha }}. \tag1.90
$$
We also note here that, for instance, for the canonical and Berwald
connections, Christoffel d--symbols we can express the d-spinor connection
(1.90)
through corresponding locally adapted derivations of components of metric
and N-connection by introducing respectively coefficients (1.22) and (1.20),
or (1.23) instead of
$\Gamma _{{\bold \gamma}{\bold \beta }}^{{\bold \alpha }}$ in
(1.90) and than in (1.89).
\vskip15pt

\subhead{I.6.3 D-spinors of la-space curvature and torsion}\endsubhead

The d-tensor indices of the commutator (1.35), $\Delta _{\alpha \beta },$
can be transformed into d-spinor ones:%
$$
\square _{\underline{\alpha }\underline{\beta }}=(\sigma ^{\alpha \beta })_{%
\underline{\alpha }\underline{\beta }}\Delta _{\alpha \beta }=(\square _{%
\underline{i}\underline{j}},\square _{\underline{a}\underline{b}}), \tag1.91
$$
with h- and v-components,
$$
\square _{\underline{i}\underline{j}}=(\sigma ^{\alpha \beta })_{\underline{i}%
\underline{j}}\Delta _{\alpha \beta }\text{ and }\square _{\underline{a}%
\underline{b}}=(\sigma ^{\alpha \beta })_{\underline{a}\underline{b}}\Delta
_{\alpha \beta },
$$
being symmetric or antisymmetric in dependence of corresponding values of
dimensions $n\,$ and $m$ (see eight-fold parametizations (1.69),\ (1.70) and
(1.71)). Considering the actions of operator (1.91) on d-spinors $\pi ^{%
\underline{\gamma }}$ and $\mu _{\underline{\gamma }}$ we introduce
the d-spinor curvature $X_{\underline{\delta }\quad \underline{\alpha }%
\underline{\beta }}^{\quad \underline{\gamma }}\,$ as to satisfy equations%
$$
\square _{\underline{\alpha }\underline{\beta }}\ \pi ^{\underline{\gamma }}=X_{%
\underline{\delta }\quad \underline{\alpha }\underline{\beta }}^{\quad
\underline{\gamma }}\pi ^{\underline{\delta }}
\text{ and }
\square _{\underline{\alpha }\underline{\beta }}\ \mu _{\underline{\gamma }%
}=X_{\underline{\gamma }\quad \underline{\alpha }\underline{\beta }%
}^{\quad \underline{\delta }}\mu _{\underline{\delta }}.\tag1.92
$$
The gravitational d-spinor $\Psi _{\underline{\alpha }\underline{\beta }%
\underline{\gamma }\underline{\delta }}$ is defined by a corresponding
symmetrization of d-spinor indices:%
$$
\Psi _{\underline{\alpha }\underline{\beta }\underline{\gamma }\underline{%
\delta }}=X_{(\underline{\alpha }|\underline{\beta }|\underline{\gamma }%
\underline{\delta })}.
$$
We note that d-spinor tensors $X_{\underline{\delta }\quad \underline{\alpha
}\underline{\beta }}^{\quad \underline{\gamma }}$ and $\Psi _{\underline{%
\alpha }\underline{\beta }\underline{\gamma }\underline{\delta }}\,$ are
transformed into similar 2-spinor objects on locally isotropic spaces
\cite{Penrose and Rindler 1984, 1986}
if we consider vanishing of the N--connection structure and a limit
to a locally isotropic space.

Putting $\delta _{\underline{\gamma }}^{\quad {\underline{\bold\gamma }}}$
instead of $\mu _{\underline{\gamma }}$ in (1.92) and using (1.93)
we can express respectively the curvature and gravitational d-spinors as
$$
X_{\underline{\gamma }\underline{\delta }\underline{\alpha }\underline{\beta
}}=\delta _{\underline{\delta }\underline{{\bold \tau }}}\square _{\underline{%
\alpha }\underline{\beta }}\delta _{\underline{\gamma }}^{\quad {
\underline{\bold\tau }}}
\text{ and }
\Psi _{\underline{\gamma }\underline{\delta }\underline{\alpha }\underline{%
\beta }}=\delta _{\underline{\delta }\underline{{\bold \tau }}}\square _{(%
\underline{\alpha }\underline{\beta }}\delta _{\underline{\gamma })}^{\quad
{ \underline{\bold\tau }}}.
$$

The d-spinor torsion $T_{\qquad \underline{\alpha }\underline{\beta }}^{%
\underline{\gamma }_1\underline{\gamma }_2}$ is defined similarly as for
d-tensors (see (1.36)) by using the d-spinor commutator (1.91) and equations
$$
\square _{\underline{\alpha }\underline{\beta }}f=T_{\qquad \underline{\alpha }%
\underline{\beta }}^{\underline{\gamma }_1\underline{\gamma }_2}\nabla _{%
\underline{\gamma }_1\underline{\gamma }_2}f. \tag1.93
$$

The d-spinor components $R_{\underline{\gamma }_1\underline{\gamma }_2\qquad
\underline{\alpha }\underline{\beta }}^{\qquad \underline{\delta }_1%
\underline{\delta }_2}$ of the curvature d-tensor $R_{\gamma \quad \alpha
\beta }^{\quad \delta }$ can be computed by using relations (1.90), and
(1.91) and (1.93) as to satisfy the equations (the d-spinor analogous of
equations (1.37) )%
$$
(\square _{\underline{\alpha }\underline{\beta }}-T_{\qquad \underline{\alpha }%
\underline{\beta }}^{\underline{\gamma }_1\underline{\gamma }_2}\nabla _{%
\underline{\gamma }_1\underline{\gamma }_2})V^{\underline{\delta }_1%
\underline{\delta }_2}=R_{\underline{\gamma }_1\underline{\gamma }_2\qquad
\underline{\alpha }\underline{\beta }}^{\qquad \underline{\delta }_1%
\underline{\delta }_2}V^{\underline{\gamma }_1\underline{\gamma }_2},
 \tag1.94$$
here d-vector $V^{\underline{\gamma }_1\underline{\gamma }_2}$ is
considered as a product of d-spinors, i.e. $V^{\underline{\gamma }_1%
\underline{\gamma }_2}=\nu ^{\underline{\gamma }_1}\mu ^{\underline{\gamma }%
_2}$. We find
$$ \multline
R_{\underline{\gamma }_1\underline{\gamma }_2\qquad \underline{\alpha }%
\underline{\beta }}^{\qquad \underline{\delta }_1\underline{\delta }%
_2}=\left( X_{\underline{\gamma }_1~\underline{\alpha }\underline{\beta }%
}^{\quad \underline{\delta }_1}+T_{\qquad \underline{\alpha }\underline{%
\beta }}^{\underline{\tau }_1\underline{\tau }_2}\quad \gamma _{\quad
\underline{\tau }_1\underline{\tau }_2\underline{\gamma }_1}^{\underline{%
\delta }_1}\right) \delta _{\underline{\gamma }_2}^{\quad \underline{\delta }%
_2}+\\
\left( X_{\underline{\gamma }_2~\underline{\alpha }\underline{\beta }%
}^{\quad \underline{\delta }_2}+T_{\qquad \underline{\alpha }\underline{%
\beta }}^{\underline{\tau }_1\underline{\tau }_2}\quad \gamma _{\quad
\underline{\tau }_1\underline{\tau }_2\underline{\gamma }_2}^{\underline{%
\delta }_2}\right) \delta _{\underline{\gamma }_1}^{\quad \underline{\delta }%
_1}.\endmultline \tag1.95 $$

It is convenient to use this d-spinor expression for the curvature d-tensor
$$\multline
R_{\underline{\gamma }_1\underline{\gamma }_2\qquad \underline{\alpha }_1%
\underline{\alpha }_2\underline{\beta }_1\underline{\beta }_2}^{\qquad
\underline{\delta }_1\underline{\delta }_2}=\left( X_{\underline{\gamma }_1~%
\underline{\alpha }_1\underline{\alpha }_2\underline{\beta }_1\underline{%
\beta }_2}^{\quad \underline{\delta }_1}+T_{\qquad \underline{\alpha }_1%
\underline{\alpha }_2\underline{\beta }_1\underline{\beta }_2}^{\underline{%
\tau }_1\underline{\tau }_2}~\gamma _{\quad \underline{\tau }_1\underline{%
\tau }_2\underline{\gamma }_1}^{\underline{\delta }_1}\right) \delta _{%
\underline{\gamma }_2}^{\quad \underline{\delta }_2}+\\
\left( X_{\underline{%
\gamma }_2~\underline{\alpha }_1\underline{\alpha }_2\underline{\beta }_1%
\underline{\beta }_2}^{\quad \underline{\delta }_2}+T_{\qquad \underline{%
\alpha }_1\underline{\alpha }_2\underline{\beta }_1\underline{\beta }_2~}^{%
\underline{\tau }_1\underline{\tau }_2}\gamma _{\quad \underline{\tau }_1%
\underline{\tau }_2\underline{\gamma }_2}^{\underline{\delta }_2}\right)
\delta _{\underline{\gamma }_1}^{\quad \underline{\delta }_1}
\endmultline$$
in order to get the d-spinor components of the Ricci d-tensor%
$$ \multline
R_{\underline{\gamma }_1\underline{\gamma }_2\underline{\alpha }_1\underline{%
\alpha }_2}=R_{\underline{\gamma }_1\underline{\gamma }_2\qquad \underline{%
\alpha }_1\underline{\alpha }_2\underline{\delta }_1\underline{\delta }%
_2}^{\qquad \underline{\delta }_1\underline{\delta }_2}= \\
X_{\underline{\gamma }_1~\underline{\alpha }_1\underline{\alpha }_2%
\underline{\delta }_1\underline{\gamma }_2}^{\quad \underline{\delta }%
_1}+T_{\qquad \underline{\alpha }_1\underline{\alpha }_2\underline{\delta }_1%
\underline{\gamma }_2}^{\underline{\tau }_1\underline{\tau }_2}~\gamma
_{\quad \underline{\tau }_1\underline{\tau }_2\underline{\gamma }_1}^{%
\underline{\delta }_1}+X_{\underline{\gamma }_2~\underline{\alpha }_1%
\underline{\alpha }_2\underline{\delta }_1\underline{\gamma }_2}^{\quad
\underline{\delta }_2}+T_{\qquad \underline{\alpha }_1\underline{\alpha }_2%
\underline{\gamma }_1\underline{\delta }_2~}^{\underline{\tau }_1\underline{%
\tau }_2}\gamma _{\quad \underline{\tau }_1\underline{\tau }_2\underline{%
\gamma }_2}^{\underline{\delta }_2} \endmultline \tag1.96
$$
and this d-spinor decomposition of the scalar curvature:%
$$\multline                                                                          q
\overleftarrow{R}=R_{\qquad \underline{\alpha }_1\underline{\alpha }_2}^{%
\underline{\alpha }_1\underline{\alpha }_2}=X_{\quad ~\underline{~\alpha }%
_1\quad \underline{\delta }_1\underline{\alpha }_2}^{\underline{\alpha }_1%
\underline{\delta }_1~~\underline{\alpha }_2}+T_{\qquad ~~\underline{\alpha }%
_2\underline{\delta }_1}^{\underline{\tau }_1\underline{\tau }_2\underline{%
\alpha }_1\quad ~\underline{\alpha }_2}~\gamma _{\quad \underline{\tau }_1%
\underline{\tau }_2\underline{\alpha }_1}^{\underline{\delta }_1}+\\
 X_{\qquad
\quad \underline{\alpha }_2\underline{\delta }_2\underline{\alpha }_1}^{%
\underline{\alpha }_2\underline{\delta }_2\underline{\alpha }_1}+T_{\qquad
\underline{\alpha }_1\quad ~\underline{\delta }_2~}^{\underline{\tau }_1%
\underline{\tau }_2~~\underline{\alpha }_2\underline{\alpha }_1}\gamma
_{\quad \underline{\tau }_1\underline{\tau }_2\underline{\alpha }_2}^{%
\underline{\delta }_2}. \endmultline \tag1.97$$

Putting (1.96) and (1.97) into (1.43) and, correspondingly, (1.41) we find
the d-spinor components of the Einstein and $\Phi _{\alpha \beta }$
d-tensors:%
$$\multline
\overleftarrow{G}_{\gamma \alpha }=\overleftarrow{G}_{\underline{\gamma }_1%
\underline{\gamma }_2\underline{\alpha }_1\underline{\alpha }_2}=  X_{%
\underline{\gamma }_1~\underline{\alpha }_1\underline{\alpha }_2\underline{%
\delta }_1\underline{\gamma }_2}^{\quad \underline{\delta }_1}+T_{\qquad
\underline{\alpha }_1\underline{\alpha }_2\underline{\delta }_1\underline{%
\gamma }_2}^{\underline{\tau }_1\underline{\tau }_2}~\gamma _{\quad
\underline{\tau }_1\underline{\tau }_2\underline{\gamma }_1}^{\underline{%
\delta }_1}+ \\
X_{\underline{\gamma }_2~\underline{\alpha }_1\underline{\alpha }%
_2\underline{\delta }_1\underline{\gamma }_2}^{\quad \underline{\delta }%
_2}+
T_{\qquad \underline{\alpha }_1\underline{\alpha }_2\underline{\gamma }_1%
\underline{\delta }_2~}^{\underline{\tau }_1\underline{\tau }_2}\gamma
_{\quad \underline{\tau }_1\underline{\tau }_2\underline{\gamma }_2}^{%
\underline{\delta }_2}-
\frac 12\varepsilon _{\underline{\gamma }_1\underline{\alpha }_1}\varepsilon
_{\underline{\gamma }_2\underline{\alpha }_2} [ X_{\quad ~\underline{%
~\beta }_1\quad \underline{\mu }_1\underline{\beta }_2}^{\underline{\beta }_1%
\underline{\mu }_1~~\underline{\beta }_2}+\\ T_{\qquad ~~\underline{\beta }_2%
\underline{\mu }_1}^{\underline{\tau }_1\underline{\tau }_2\underline{\beta }%
_1\quad ~\underline{\beta }_2}~\gamma _{\quad \underline{\tau }_1\underline{%
\tau }_2\underline{\beta }_1}^{\underline{\mu }_1}+
X_{\qquad \quad
\underline{\beta }_2\underline{\mu }_2\underline{\nu }_1}^{\underline{\beta }%
_2\underline{\mu }_2\underline{\nu }_1}+T_{\qquad \underline{\beta }_1\quad ~%
\underline{\delta }_2~}^{\underline{\tau }_1\underline{\tau }_2~~\underline{%
\beta }_2\underline{\beta }_1}\gamma _{\quad \underline{\tau }_1\underline{%
\tau }_2\underline{\beta }_2}^{\underline{\delta }_2} ] \endmultline \tag1.98
$$
and%
$$\multline
\Phi _{\gamma \alpha }=\Phi _{\underline{\gamma }_1\underline{\gamma }_2%
\underline{\alpha }_1\underline{\alpha }_2}=
\frac 1{2(n+m)}\varepsilon _{\underline{\gamma }_1\underline{\alpha }%
_1}\varepsilon _{\underline{\gamma }_2\underline{\alpha }_2} [ X_{\quad ~%
\underline{~\beta }_1\quad \underline{\mu }_1\underline{\beta }_2}^{%
\underline{\beta }_1\underline{\mu }_1~~\underline{\beta }_2}+\\
T_{\qquad ~~%
\underline{\beta }_2\underline{\mu }_1}^{\underline{\tau }_1\underline{\tau }%
_2\underline{\beta }_1\quad ~\underline{\beta }_2}~\gamma _{\quad \underline{%
\tau }_1\underline{\tau }_2\underline{\beta }_1}^{\underline{\mu }%
_1}+
X_{\qquad \quad \underline{\beta }_2\underline{\mu }_2\underline{\nu }%
_1}^{\underline{\beta }_2\underline{\mu }_2\underline{\nu }_1}+
T_{\qquad
\underline{\beta }_1\quad ~\underline{\delta }_2~}^{\underline{\tau }_1%
\underline{\tau }_2~~\underline{\beta }_2\underline{\beta }_1}\gamma _{\quad
\underline{\tau }_1\underline{\tau }_2\underline{\beta }_2}^{\underline{%
\delta }_2} ] -
\frac 12 [ X_{\underline{\gamma }_1~\underline{\alpha }_1\underline{%
\alpha }_2\underline{\delta }_1\underline{\gamma }_2}^{\quad \underline{%
\delta }_1}+\\ T_{\qquad \underline{\alpha }_1\underline{\alpha }_2\underline{%
\delta }_1\underline{\gamma }_2}^{\underline{\tau }_1\underline{\tau }%
_2}~\gamma _{\quad \underline{\tau }_1\underline{\tau }_2\underline{\gamma }%
_1}^{\underline{\delta }_1} +
X_{\underline{\gamma }_2~\underline{\alpha }_1%
\underline{\alpha }_2\underline{\delta }_1\underline{\gamma }_2}^{\quad
\underline{\delta }_2}+T_{\qquad \underline{\alpha }_1\underline{\alpha }_2%
\underline{\gamma }_1\underline{\delta }_2~}^{\underline{\tau }_1\underline{%
\tau }_2}\gamma _{\quad \underline{\tau }_1\underline{\tau }_2\underline{%
\gamma }_2}^{\underline{\delta }_2} ]. \endmultline \tag1.99$$

The components of the conformal Weyl d-spinor can be computed by putting
d-spinor values of the curvature (1.95) and Ricci (1.96) d-tensors into
corresponding expression for the d-tensor (1.40). We omit this calculus in
this work.
\vskip25pt

\head{I.7\ Field  Equations on Locally Anisotropic Spaces}\endhead

The problem of formulation gravitational and gauge field equations on
different types of la--spaces is considered, for instance, in
\cite{Miron and Anastasiei 1994}, \cite{Bejancu 1990},
 \cite{Asanov and Ponomarenko 1988}  and \cite{Vacaru and Goncharenko 1995}.
 In this section we shall introduce the basic field equations for
gravitational and matter field la-interactions in a generalized form for
generic la-spaces.
\vskip15pt
\subhead{I.7.1\ Locally anisotropic scalar field interactions}\endsubhead

Let $\varphi \left( u\right) =(\varphi _1\left( u\right) ,\varphi _2\left(
u\right) \dot ,...,\varphi _k\left( u\right) )$ be a complex k-component
scalar field of mass $\mu $ on la-space ${\Cal E}.$ The d-covariant
generalization of the conformally invariant (in the massless case) scalar
field equation \cite{Penrose and Rindler 1984, 1986} can be defined by using
the d'Alambert locally anisotropic operator
 $\square =D^\alpha D_\alpha $, where $D_\alpha $
is a d-covariant derivation on ${\Cal E}$ satisfying conditions (1.14) and
(1.15):%
$$
(\square +\frac{n+m-2}{4(n+m-1)}\overleftarrow{R}+\mu ^2)\varphi \left(
u\right) =0. \tag1.100
$$
We must change d-covariant derivation $D_{\alpha }$ into $%
^{\diamond }D_\alpha =D_\alpha +ieA_\alpha $ and take into account the
d-vector current
$$
J_\alpha ^{(0)}\left( u\right) =i(\left( \overline \varphi  \left( u\right)
D_\alpha \varphi \left( u\right) -D_\alpha \overline \varphi \left( u\right)
)\varphi \left( u\right) \right)
$$
if interactions between locally anisotropic electromagnetic field ( d-vector
potential $A_\alpha $ ), where $e$ is the electromagnetic constant, and
charged scalar field $\varphi $ are considered. The equations (1.100) are
(locally adapted to the N-connection structure) Euler equations for the
Lagrangian%
$$
{\Cal L}^{(0)}\left( u\right) =\sqrt{|g|}\left[ g^{\alpha \beta }\delta
_\alpha \overline \varphi \left( u\right) \delta _\beta \varphi \left( u\right)
-\left( \mu ^2+\frac{n+m-2}{4(n+m-1)}\right) \overline \varphi \left( u\right)
\varphi \left( u\right) \right], \tag1.101
$$
where $|g|=det g_{\alpha \beta}.$

The locally adapted variations of the action with Lagrangian (1.101) on
variables $\varphi \left( u\right) $ and $\overline \varphi \left( u\right) $
leads to the locally anisotropic generalization of the energy-momentum
tensor,%
$$
E_{\alpha \beta }^{(0,can)}\left( u\right) =\delta _\alpha \overline \varphi
 \left( u\right) \delta _\beta \varphi \left( u\right) +\delta _\beta
\overline \varphi \left( u\right) \delta _\alpha \varphi \left( u\right) -
\frac 1{%
\sqrt{|g|}}g_{\alpha \beta }{\Cal L}^{(0)}\left( u\right) , \tag1.102
$$
and a similar variation on the components of a d-metric (1.12) leads to a
symmetric metric energy-momentum d-tensor,%
$$
E_{\alpha \beta }^{(0)}\left( u\right) =E_{(\alpha \beta )}^{(0,can)}\left(
u\right) -\frac{n+m-2}{2(n+m-1)}\left[ R_{(\alpha \beta )}+D_{(\alpha
}D_{\beta )}-g_{\alpha \beta }\square \right] \overline \varphi \left( u\right)
\varphi \left( u\right) . \tag1.103
$$
Here we note that we can obtain a nonsymmetric energy-momentum d-tensor if
we firstly vary on $G_{\alpha \beta }$ and than impose constraints of type
(1.10) in order to have a compatibility with the N-connection structure. We
also conclude that the existence of a N-connection in v-bundle ${\Cal E}$
results in a nonequivalence of energy-momentum d-tensors (1.102) and (1.103),
nonsymmetry of the Ricci tensor (see (1.33)), nonvanishing of the
d-covariant derivation of the Einstein d-tensor (1.43), $D_\alpha
\overleftarrow{G}^{\alpha \beta }\neq 0$ and, in consequence, a
corresponding breaking of conservation laws on la-spaces when $D_\alpha
E^{\alpha \beta }\neq 0\,$  \cite{Miron and Anastasiei 1987, 1994}.
The problem of formulation
of conservation laws on la-spaces is discussed in details and two variants
of its solution (by using nearly autoparallel maps and tensor integral
formalism on la-multispaces) are proposed in
\cite{Vacaru and Ostaf 1994, 1996a} (see Chapter 3).
In this section we shall present only straightforward generalizations of
  field equations and
 necessary formulas for energy-momentum d-tensors of matter fields on
${\Cal E}$ considering that it is naturally that the conservation laws (usually
being consequences of global, local and/or intrinsic symmetries of the
fundamental space-time and of the type of field interactions) have to be
broken on locally anisotropic spaces.
\vskip15pt

\subhead{I.7.2 Proca equations on la--spaces}\endsubhead

Let consider a d-vector field $\varphi _\alpha \left( u\right) $
with mass $\mu ^2$ (locally anisotropic Proca field) interacting with
exterior la-gravitational field. From the Lagrangian
$$
{\Cal L}^{(1)}\left( u\right) =\sqrt{\left| g\right| }\left[ -\frac
12 {\overline f}_{\alpha \beta }\left( u\right)
 f^{\alpha \beta }\left( u\right) +\mu
^2 {\overline \varphi}_\alpha \left( u\right) \varphi ^\alpha \left( u\right)
\right] , \tag1.104
$$
where $f_{\alpha \beta }=D_\alpha \varphi _\beta -D_\beta \varphi _\alpha ,$
one follows the Proca equations on la-spaces
$$
D_\alpha f^{\alpha \beta }\left( u\right) +\mu ^2\varphi ^\beta \left(
u\right) =0. \tag1.105
$$
Equations (1.105) are a first type constraints for $\beta =0.$ Acting with $%
D_\alpha $ on (1.105), for $\mu \neq 0$ we obtain second type constraints%
$$
D_\alpha \varphi ^\alpha \left( u\right) =0. \tag1.106
$$

Putting (1.106) into (1.105) we obtain second order field equations with respect
to $\varphi _{\alpha }$ :%
$$
\square \varphi _\alpha \left( u\right) +R_{\alpha \beta }\varphi ^\beta \left(
u\right) +\mu ^2\varphi _\alpha \left( u\right) =0. \tag1.107
$$
The energy-momentum d-tensor and d-vector current following from the (1.107)
can be written as
$$
E_{\alpha \beta }^{(1)}\left( u\right) =-g^{\varepsilon \tau }\left(
{\overline f}_{\beta \tau } f_{\alpha \varepsilon }+
{\overline f}_{\alpha \varepsilon} f_{\beta \tau }\right) +
\mu ^2\left({\overline \varphi}_\alpha \varphi _\beta
+{\overline \varphi}_\beta \varphi _\alpha \right) -
\frac{g_{\alpha \beta }}{\sqrt{%
\left| g\right| }}{\Cal L}^{(1)}\left( u\right)
$$
and%
$$
J_\alpha ^{\left( 1\right) }\left( u\right) =
i\left( {\overline  f}_{\alpha \beta}
\left( u\right) \varphi ^\beta \left( u\right) - {\overline \varphi}^{\beta}
\left( u\right) f_{\alpha \beta }\left( u\right) \right) .
$$

For $\mu =0$ the d-tensor $f_{\alpha \beta }$ and the Lagrangian (1.104) are
invariant with respect to locally anisotropic gauge transforms of type
$$
\varphi _\alpha \left( u\right) \rightarrow \varphi _\alpha \left( u\right)
+\delta _\alpha \Lambda \left( u\right) ,
$$
where $\Lambda \left( u\right) $ is a d-differentiable scalar function, and
we obtain a locally anisot\-rop\-ic variant of Maxwell theory.
\vskip15pt

\subhead{I.7.3 La-gravitons on la-backgrounds}\endsubhead

Let a massless d-tensor field $h_{\alpha \beta }\left( u\right) $ is
interpreted as a small perturbation of the locally anisotropic background
metric d-field $g_{\alpha \beta }\left( u\right) .$ Considering, for
simplicity, a torsionless background we have locally anisotropic Fierz--Pauli
equations%
$$
\square h_{\alpha \beta }\left( u\right) +2R_{\tau \alpha \beta \nu }\left(
u\right) ~h^{\tau \nu }\left( u\right) =0 \tag1.108
$$
and d-gauge conditions%
$$
D_\alpha h_\beta ^\alpha \left( u\right) =0,\quad h\left( u\right) \equiv
h_\beta ^\alpha (u)=0, \tag1.109
$$
where $R_{\tau \alpha \beta \nu }\left( u\right) $ is curvature d-tensor of
the la-background space (these formulae can be obtained by using a
perturbation formalism with respect to $h_{\alpha \beta }\left( u\right) $).

We note that we can rewrite d-tensor formulas (1.100)-(1.109) into similar
d-spinor ones by using formulas (1.78)-(1.80), (1.90), (1.92) and
(1.96)-(1.105) (for simplicity, we omit these
considerations in this work).
\vskip15pt

\subhead{I.7.4\ Locally anisotropic Dirac equations}\endsubhead

Let denote the Dirac d-spinor field on ${\Cal E}$ as $\psi \left( u\right)
=\left( \psi ^{\underline{\alpha }}\left( u\right) \right) $ and consider as
the generalized Lorentz transforms the group of automorphysm of the metric $%
G_{\widehat{\alpha }\widehat{\beta }}$ (see the la-frame decomposition of
d-metric (1.54)).The d-covariant derivation of field $\psi \left( u\right) $
is written as
$$
\overrightarrow{\nabla _\alpha }\psi =\left[ \delta _\alpha +\frac 14C_{%
\widehat{\alpha }\widehat{\beta }\widehat{\gamma }}\left( u\right) ~l_\alpha
^{\widehat{\alpha }}\left( u\right) \sigma ^{\widehat{\beta }}\sigma ^{%
\widehat{\gamma }}\right] \psi , \tag1.110
$$
where coefficients $C_{\widehat{\alpha }\widehat{\beta }\widehat{\gamma }%
}=\left( D_\gamma l_{\widehat{\alpha }}^\alpha \right) l_{\widehat{\beta }%
\alpha }l_{\widehat{\gamma }}^\gamma $ generalize for la-spaces the
corresponding Ricci coefficients on Riemannian spaces.
Using $\sigma $%
-objects $\sigma ^\alpha \left( u\right) $ (see (1.79) and (1.55))
we define the Dirac equations on la-spaces:
$$
(i\sigma ^\alpha \left( u\right) \overrightarrow{\nabla _\alpha }-\mu )\psi
=0, \tag1.111$$
which are the Euler equations for the Lagrangian%
$$\multline
{\Cal L}^{(1/2)}\left( u\right) =\sqrt{\left| g\right| } \{ [ \psi
^{+}\left( u\right) \sigma ^\alpha \left( u\right) \overrightarrow{\nabla
_\alpha }\psi \left( u\right) - \\
 (\overrightarrow{\nabla _\alpha }\psi
^{+}\left( u\right) )\sigma ^\alpha \left( u\right) \psi \left( u\right)
 ] -\mu \psi ^{+}\left( u\right) \psi \left( u\right) \} ,
\endmultline \tag1.112
$$
where $\psi ^{+}\left( u\right) $ is the complex conjugation and transposition
of the column$~\psi \left( u\right) .$

From (1.112) we obtain the d-metric energy-momentum d-tensor%
$$
E_{\alpha \beta }^{(1/2)}\left( u\right) =\frac i4 [ \psi ^{+}\left(
u\right) \sigma _\alpha \left( u\right) \overrightarrow{\nabla _\beta }\psi
\left( u\right) +\psi ^{+}\left( u\right) \sigma _\beta \left( u\right)
\overrightarrow{\nabla _\alpha }\psi \left( u\right) -$$
$$ (\overrightarrow{%
\nabla _\alpha }\psi ^{+}\left( u\right) )\sigma _\beta \left( u\right) \psi
\left( u\right) -(\overrightarrow{\nabla _\beta }\psi ^{+}\left( u\right)
)\sigma _\alpha \left( u\right) \psi \left( u\right) ]
$$
and the d-vector source%
$$
J_\alpha ^{(1/2)}\left( u\right) =\psi ^{+}\left( u\right) \sigma _\alpha
\left( u\right) \psi \left( u\right) .
$$
We emphasize that la-interactions with exterior gauge fields can be
introduced by changing the la-partial derivation from (1.110) in this manner:%
$$
\delta _\alpha \rightarrow \delta _\alpha +ie^{\star }B_\alpha ,
\tag1.113
$$
where $e^{\star }$ and $B_\alpha $ are respectively the constant d-vector
potential of la-gauge interactions on la-spaces (see
\cite{Vacaru and Goncharenko 1995} and the next
subsection).
\vskip15pt

\subhead{I.7.5 D-spinor Yang-Mills equations}\endsubhead

We consider a v-bundle ${\Cal B}_E,~\pi _B:{\Cal B\rightarrow E,}$ on
la-space ${\Cal E}.$ Additionally to d-tensor and d-spinor indices we
shall use capital Greek letters, $\Phi ,\Upsilon ,\Xi ,\Psi ,...$ for fibre
(of this bundle) indices (see details in
\cite{Penrose and Rindler 1984, 1986} for the case when the base
space of the v-bundle $\pi _B$ is a locally isotropic space--time). Let $%
\underline{\nabla }_\alpha $ be, for simplicity, a torsionless,
linear connection in ${\Cal B}_E$ satisfying conditions:
$$
\underline{\nabla }_\alpha :{\bold \Upsilon }^\Theta \rightarrow {\bold %
\Upsilon }_\alpha ^\Theta \quad \left[ \text{ or }{\bold \Xi }^\Theta
\rightarrow {\bold \Xi }_\alpha ^\Theta \right] ,
$$
$$
\underline{\nabla }_\alpha \left( \lambda ^\Theta +\nu ^\Theta
\right) = \underline{\nabla }_\alpha \lambda ^\Theta +%
\underline{\nabla }_\alpha \nu ^\Theta ,
$$
$$
\underline{\nabla }_\alpha ~(f\lambda ^\Theta )=\lambda ^\Theta
\underline{\nabla }_\alpha f+f\underline{\nabla }_\alpha
\lambda ^\Theta ,\quad f\in {\bold \Upsilon }^\Theta ~[\text{ or }{\bold \Xi }%
^\Theta ],
$$
where by ${\bold \Upsilon }^\Theta ~\left( {\bold \Xi }^\Theta \right) $ we
denote the module of sections of the real (complex) v-bundle ${\Cal B}_E$
provided with the abstract index $\Theta .$ The curvature of connection $%
\underline{\nabla }_\alpha $ is defined as
$$
K_{\alpha \beta \Omega }^{\qquad \Theta }\lambda ^\Omega =\left(
\underline{\nabla }_\alpha \underline{\nabla }_\beta -%
\underline{\nabla }_\beta \underline{\nabla }_\alpha
\right) \lambda ^\Theta .
$$

For Yang-Mills fields as a rule one considers that ${\Cal B}_E$ is enabled
with a unitary (complex) structure (complex conjugation changes mutually the
upper and lower Greek indices). It is useful to introduce instead of $%
K_{\alpha \beta \Omega }^{\qquad \Theta }$ a Hermitian matrix $F_{\alpha
\beta \Omega }^{\qquad \Theta }=i$ $K_{\alpha \beta \Omega }^{\qquad \Theta }
$ connected with components of the Yang-Mills d-vector potential $B_{\alpha
\Xi }^{\quad \Phi }$ according the formula:%
$$
\frac 12F_{\alpha \beta \Xi }^{\qquad \Phi }=\underline{\nabla }%
_{[\alpha }B_{\beta ]\Xi }^{\quad \Phi }-iB_{[\alpha |\Lambda |}^{\quad \Phi
}B_{\beta ]\Xi }^{\quad \Lambda }, \tag1.114
$$
where the la-space indices commute with capital Greek indices. The gauge
transforms are written in the form:%
$$
B_{\alpha \Theta }^{\quad \Phi }\mapsto B_{\alpha \widehat{\Theta }}^{\quad
\widehat{\Phi }}=B_{\alpha \Theta }^{\quad \Phi }~s_\Phi ^{\quad \widehat{%
\Phi }}~q_{\widehat{\Theta }}^{\quad \Theta }+is_\Theta ^{\quad \widehat{%
\Phi }}\underline{\nabla }_\alpha ~q_{\widehat{\Theta }}^{\quad
\Theta },
$$
$$
F_{\alpha \beta \Xi }^{\qquad \Phi }\mapsto F_{\alpha \beta \widehat{\Xi }%
}^{\qquad \widehat{\Phi }}=F_{\alpha \beta \Xi }^{\qquad \Phi }s_\Phi
^{\quad \widehat{\Phi }}q_{\widehat{\Xi }}^{\quad \Xi },
$$
where matrices $s_\Phi ^{\quad \widehat{\Phi }}$ and $q_{\widehat{\Xi }%
}^{\quad \Xi }$ are mutually inverse (Hermitian conjugated in the unitary
case). The Yang-Mills equations on  la--spaces
 (see details in the next Chapter) are written in this form:%
$$
\underline{\nabla }^\alpha F_{\alpha \beta \Theta }^{\qquad \Psi
}=J_{\beta \ \Theta}^{\qquad \Psi}   , \tag1.115
$$
$$
\underline{\nabla }_{[\alpha }F_{\beta \gamma ]\Theta }^{\qquad \Xi
}=0. \tag1.116
$$
We must introduce deformations of connection of type (1.14) and (1.15),
 $\underline{\nabla }_\alpha ^{\star }~\longrightarrow
\underline{\nabla }_\alpha +P_\alpha, $ (the deformation d-tensor $%
P_\alpha $ is induced by the torsion in v-bundle ${\Cal B}_E)$ into the
definition of the curvature of la-gauge fields (1.114) and motion equations
(1.115) and (1.116) if interactions are modeled on a generic la-space.

Now we can write out the field equations of the Einstein-Cartan theory in
the d-spinor form. So, for the Einstein equations (1.42) we have
$$
\overleftarrow{G}_{\underline{\gamma }_1\underline{\gamma }_2\underline{%
\alpha }_1\underline{\alpha }_2}+\lambda \varepsilon _{\underline{\gamma }_1%
\underline{\alpha }_1}\varepsilon _{\underline{\gamma }_2\underline{\alpha }%
_2}=\kappa E_{\underline{\gamma }_1\underline{\gamma }_2\underline{\alpha }_1%
\underline{\alpha }_2},
$$
with $\overleftarrow{G}_{\underline{\gamma }_1\underline{\gamma }_2%
\underline{\alpha }_1\underline{\alpha }_2}$ from (1.98), or, by using the
d-tensor (1.99),%
$$
\Phi _{\underline{\gamma }_1\underline{\gamma }_2\underline{\alpha }_1%
\underline{\alpha }_2}+(\frac{\overleftarrow{ R}}4-
\frac \lambda 2)\varepsilon _{\underline{%
\gamma }_1\underline{\alpha }_1}\varepsilon _{\underline{\gamma }_2%
\underline{\alpha }_2}=-\frac \kappa 2E_{\underline{\gamma }_1\underline{%
\gamma }_2\underline{\alpha }_1\underline{\alpha }_2},
$$
which are the d-spinor equivalent of the equations (1.44). These equations
 must be completed by the algebraic equations
(1.45) for the d-torsion and d-spin density with d-tensor indices changed
into corresponding d-spinor ones.
\vskip80pt
\newpage
\topmatter
\title\chapter{2} Gauge Fields and Locally Anisotropic
Gravity \endtitle \endtopmatter

The aim of this Chapter is twofold. The first objective is to develop some of
 our results \cite{Vacaru and Goncharenko 1995} on formulation of
a geometrical approach to
interactions of Yang-Mills fields on spaces with local anisotropy in the
framework of the theory of linear connections in vector bundles (with
semisimple structural groups) on la-spaces. The second objective is to extend
the geometrical formalism in a manner including theories with nonsemisimple
groups which permit a unique fiber bundle treatment for both locally
anisotropic Yang-Mills field and gravitational interactions. In general
lines, we shall follow the ideas and geometric methods proposed in
\cite{Bishop and Crittenden 1964}, \cite{Popov 1975},
\cite{Popov and Dikhin 1975}, \cite{Tseytlin 1982} and
\cite{Ponomariov, Barvinsky and Obukhov 1985},
but we shall apply them in a form convenient for
introducing into consideration geometrical constructions
 and physical theories on la-spaces.

There is a number of works on gauge models of interactions on Finsler spaces
and theirs extensions(see, for instance,
 \cite{Asanov 1985}, \cite{Bejancu 1990} and \cite{Ono and Takano 1993}).
 One has  introduced different
 variants of generalized gauge transforms, postulated corresponding
Lagrangians for gravitational, gauge and matter field interactions and
 formulated variational calculus). The main problem of such models is the
 dependence of the basic equations on chosen definition of gauge
"compensation" symmetries and on type of space and field interactions
anisotropy. In order to avoid the ambiguities connected with particular
characteristics of possible la-gauge theories we consider a "pure" geometric
 approach to gauge theories (on both locally isotropic and anisotropic spaces)
 in the framework of the theory of fiber bundles provided in general with
different types of nonlinear and linear multiconnection and metric
structures). This way, based on global geometric methods, holds also  good for
nonvariational, in the total spaces of bundles, gauge theories
(in the case of gauge gravity based on Poincare or affine gauge groups);
physical values and motion (field) equations have adequate geometric
interpretation and do not depend on the type of local anisotropy of
space-time background. It should be emphasized here that extensions for "higher
 order spaces" (see \cite{Yano and Ishihara 1973} and \cite{Miron and
 Atanasiu 1995}) can be realized in a straightforward manner.
\vskip25pt

\head{II.1 Gauge Fields on Locally Anisotropic Spaces}\endhead

This section is devoted to formulation of the geometrical background for
gauge field theories on spaces with local anisotropy.

Let $\left( P,\pi ,Gr,{\Cal E}\right) $ be a principal bundle on base ${\Cal %
E}$ (being a la--space) with structural group $Gr$ and surjective map $\pi
:P\rightarrow {\Cal E}.$ At every point $u=\left( x,y\right) \in {\Cal E}$
there is a vicinity ${\Cal U\subset E,}u\in {\Cal U,}$ with trivializing $P$
diffeomorphisms $f$ and $\varphi :$%
$$
f_{{\Cal U}}:\pi ^{-1}\left( {\Cal U}\right) \rightarrow {\Cal U\times }%
Gr,\qquad f\left( p\right) =\left( \pi \left( p\right) ,\varphi \left(
p\right) \right) ,
$$
$$
\varphi _{{\Cal U}}:\pi ^{-1}\left( {\Cal U}\right) \rightarrow Gr,%
\varphi (pq)=\varphi \left( p\right) q,\quad \forall q\in Gr,~p\in P.
$$
We remark that in the general case for two open regions%
$$
{\Cal{U,V}}\subset {\Cal{E,U\cap V}}\neq \emptyset ,f_{{\Cal{ U|}}_p}\neq f_{%
{\Cal{V|}_p}},\text{ even }p\in {\Cal{U\cap V}}.
$$

Transition functions $g_{{\Cal {UV}}}$ are defined as
$$
g_{{\Cal{UV}}}:{\Cal {U\cap V}\rightarrow }Gr,g_{{\Cal {UV}}%
}\left( u\right) =\varphi _{{\Cal U}}\left( p\right)
\left( \varphi _{{\Cal V}}\left( p\right) ^{-1}\right) ,\pi \left( p\right) =u.
$$

Hereafter we shall omit, for simplicity, the specification of trivializing
regions of maps and denote, for example, $f\equiv f_{{\Cal U}},\varphi
\equiv \varphi _{{\Cal U}},$ $s\equiv s_{{\Cal U}},$ if this will not give
rise to ambiguities.

Let $\theta \,$ be the canonical left invariant 1-form on $Gr$ with values
in algebra Lie ${\Cal G}$ of group $Gr$ uniquely defined from the relation $%
\theta \left( q\right) =q,\forall q\in {\Cal G},$ and consider a 1-form $%
\omega $ on ${\Cal{U\subset E}}$ with values in ${\Cal G}.$ Using $\theta $
and $\omega ,$ we can locally define the connection form $\Omega $ in $P$ as
a 1-form:%
$$
\Omega =\varphi ^{*}\theta +Ad~\varphi ^{-1}\left( \pi ^{*}\omega \right)
\tag2.1
$$
where $\varphi ^{*}\theta $ and $\pi ^{*}\omega $ are, respectively, forms
induced on $\pi ^{-1}\left( {\Cal U}\right) $ and $P$ by maps $\varphi $ and
$\pi $ and $\omega =s^{*}\Omega .$ The adjoint action on a form $\lambda $
with values in ${\Cal G}$ is defined as
$$
\left( Ad~\varphi ^{-1}\lambda \right) _p=\left( Ad~\varphi ^{-1}\left(
p\right) \right) \lambda _p
$$
where $\lambda _p$ is the value of form $\lambda $ at point $p\in P.$

Introducing a basis $\{\Delta _{\widehat{a}}\}$ in ${\Cal G}$ (index $%
\widehat{a}$ enumerates the generators making up this basis), we write the
1-form $\omega $ on ${\Cal E}$ as
$$
\omega =\Delta _{\widehat{a}}\omega ^{\widehat{a}}\left( u\right) ,~\omega ^{%
\widehat{a}}\left( u\right) =\omega _\mu ^{\widehat{a}}\left( u\right)
\delta u^\mu \tag2.2
$$
where $\delta u^\mu =\left( dx^i,\delta y^a\right) $ and the Einstein
summation rule on indices $\widehat{a}$ and $\mu $ is used. Functions $%
\omega _\mu ^{\widehat{a}}\left( u\right) =$ $\omega _\mu ^{\widehat{a}%
}\left( x,y\right) $ from (2.2) will be called the components of Yang-Mills
fields on la-space ${\Cal E}.$ Gauge transforms of $ \omega $ can be
geometrically interpreted as transition relations for $\omega _{{\Cal U}%
}$ and $\omega _{{\Cal V}},$ when $u\in {\Cal{U\cap V}},$%
$$
\left( \omega _{{\Cal U}}\right) _u=\left( g_{{\Cal {UV}}}^{*}\theta \right)
_u+Ad~g_{{\Cal{UV}}}\left( u\right) ^{-1}\left( \omega _{{\Cal V}}\right)
_u .\tag2.3
$$

To relate $\omega _\mu ^{\widehat{a}}$ with a covariant derivation we shall
consider a vector bundle $\Upsilon $ associated to $P.$ Let $\rho
:Gr\rightarrow GL\left( {\Bbb R}^m\right) $ and $\rho ^{\prime }:{\Cal G}%
\rightarrow End\left( E^m\right) $ be, respectively, linear representations
of group $Gr$ and Lie algebra ${\Cal G}$ (in a more general case we can
consider ${\Bbb C}^m$ instead of ${\Bbb R}^m).$ Map $\rho $ defines a left
action on $Gr$ and associated vector bundle $\Upsilon =P\times {\Bbb R}%
^m/Gr,~\pi _E:E\rightarrow {\Cal E}.$ Introducing the standard basis $\xi _{%
\underline{i}}=\{\xi _{\underline{1}},\xi _{\underline{2}},...,\xi _{%
\underline{m}}\}$ in ${\Bbb R}^m,$ we can define the right action on $%
P\times $ ${\Bbb R}^m,\left( \left( p,\xi \right) q=\left( pq,\rho \left(
q^{-1}\right) \xi \right) ,q\in Gr\right) ,$ the map induced from $P$%
$$
p:{\Bbb R}^m\rightarrow \pi _E^{-1}\left( u\right) ,\quad \left( p\left( \xi
\right) =\left( p\xi \right) Gr,\xi \in {\Bbb R}^m,\pi \left( p\right)
=u\right)
$$
and a basis of local sections $e_{\underline{i}}:U\rightarrow \pi
_E^{-1}\left( U\right) ,~e_{\underline{i}}\left( u\right) =s\left( u\right)
\xi _{\underline{i}}.$ Every section $\varsigma :{\Cal E \rightarrow }%
\Upsilon $ can be written locally as $\varsigma =\varsigma ^ie_i,\varsigma
^i\in C^\infty \left( {\Cal U}\right) .$ To every vector field $X$ on ${\Cal %
E}$ and Yang-Mills field $\omega ^{\widehat{a}}$ on $P$ we associate
operators of covariant derivations:%
$$
\nabla _X\zeta =e_{\underline{i}}\left[ X\zeta ^{\underline{i}}+B\left(
X\right) _{\underline{j}}^{\underline{i}}\zeta ^{\underline{j}}\right]
 \text{ and }
B\left( X\right) =\left( \rho ^{\prime }X\right) _{\widehat{a}}\omega ^{%
\widehat{a}}\left( X\right) .\tag2.4
$$
Transformation laws (2.3) and operators (2.4) are interrelated by these
transition transforms for values $e_{\underline{i}},\zeta ^{\underline{i}},$
and $B_\mu :$%
$$ \split
e_{\underline{i}}^{{\Cal V}}\left( u\right) =\left[ \rho g_{{\Cal{UV}}}\left(
u\right) \right] _{\underline{i}}^{\underline{j}}e_{\underline{i}}^{{\Cal U}%
},~\zeta _{{\Cal U}}^{\underline{i}}\left( u\right) =\left[ \rho g_{{\Cal{UV}}%
}\left( u\right) \right] _{\underline{i}}^{\underline{j}}\zeta _{{\Cal V}}^{%
\underline{i}},\\
B_\mu ^{{\Cal V}}\left( u\right) =\left[ \rho g_{{\Cal{UV}}}\left( u\right)
\right] ^{-1}\delta _\mu \left[ \rho g_{{\Cal{UV}}}\left( u\right) \right]
+\left[ \rho g_{{\Cal{UV}}}\left( u\right) \right] ^{-1}B_\mu ^{{\Cal U}%
}\left( u\right) \left[ \rho g_{{\Cal{UV}}}\left( u\right) \right],
\endsplit \tag2.5
$$
where $B_\mu ^{{\Cal U}}\left( u\right) =B^\mu \left( \delta /du^\mu \right)
\left( u\right) .$

Using (2.5), we can verify that the operator $\nabla _X^{{\Cal U}},$ acting
on sections of $\pi _\Upsilon :\Upsilon \rightarrow {\Cal E}$ according to
definition (2.4), satisfies the properties%
$$
\multline
\nabla _{f_1X+f_2Y}^{
{\Cal U}}=f_1\nabla _X^{{\Cal U}}+f_2\nabla _X^{{\Cal U}},~\nabla _X^{{\Cal U%
}}\left( f\zeta \right) =f\nabla _X^{{\Cal U}}\zeta +\left( Xf\right) \zeta
, \\ \nabla _X^{{\Cal U}}\zeta =\nabla _X^{{\Cal V}}\zeta ,\quad u\in {\Cal{ %
U\cap V},} f_1,f_2\in C^\infty \left( {\Cal U}\right) .
\endmultline
$$

So, we can conclude that the Yang\_Mills connection in the vector bundle $%
\pi _\Upsilon :\Upsilon \rightarrow {\Cal E}$ is not a general one, but is
induced from the principal bundle $\pi :P\rightarrow {\Cal E}$ with
structural group $Gr.$

The curvature ${\Cal K}$ of connection $\Omega $ from (2.1) is defined as%
$${\Cal K}=D\Omega ,~D=\widehat{H}\circ d \tag2.6
$$
where $d$ is the operator of exterior derivation acting on ${\Cal G}$--valued
forms as $$d\left( \Delta _{\widehat{a}}\otimes \chi ^{\widehat{a}}\right)
=\Delta _{\widehat{a}}\otimes d\chi ^{\widehat{a}}$$ and $\widehat{H}\,$ is
the horizontal projecting operator actin, for example, on the 1-form $%
\lambda $ as $\left( \widehat{H}\lambda \right) _P\left( X_p\right) =\lambda
_p\left( H_pX_p\right) ,$ where $H_p$ projects on the horizontal subspace $%
{\Cal H}_p\in P_p\left[ X_p\in {\Cal H}_p\text{ is equivalent to }\Omega
_p\left( X_p\right) =0\right] .$ We can express (2.6) locally as
$$
{\Cal K}=Ad~\varphi _{{\Cal U}}^{-1}\left( \pi ^{*}{\Cal K}_{{\Cal U}%
}\right) \tag2.7
$$
where
$$
{\Cal K}_{{\Cal U}}=d\omega _{{\Cal U}}+\frac 12\left[ \omega _{{\Cal U}%
},\omega _{{\Cal U}}\right] .\tag2.8
$$
The exterior product of ${\Cal G}$--valued form (2.8) is defined as
$$
\left[ \Delta _{\widehat{a}}\otimes \lambda ^{\widehat{a}},\Delta _{\widehat{%
b}}\otimes \xi ^{\widehat{b}}\right] =\left[ \Delta _{\widehat{a}},\Delta _{%
\widehat{b}}\right] \otimes \lambda ^{\widehat{a}}\bigwedge \xi ^{\widehat{b}%
},
$$
where
$
\lambda ^{
\widehat{a}}\bigwedge \xi ^{\widehat{b}}=\lambda ^{\widehat{a}}\xi ^{%
\widehat{b}}-\xi ^{\widehat{b}}\lambda ^{\widehat{a}}.
$\ is the antisymmetric tensor product.

Introducing structural coefficients $f_{\widehat{b}\widehat{c}}^{\quad
\widehat{a}}$ of ${\Cal G}$ satisfying relations
$
\left[ \Delta _{\widehat{b}},\Delta _{\widehat{c}}\right] =f_{\widehat{b}%
\widehat{c}}^{\quad \widehat{a}}\Delta _{\widehat{a}}
$\
we can rewrite (2.8) in a form more convenient for local considerations:%
$$
{\Cal K}_{{\Cal U}}=\Delta _{\widehat{a}}\otimes {\Cal K}_{\mu \nu }^{%
\widehat{a}}\delta u^\mu \bigwedge \delta u^\nu \tag2.9
$$
where
$$
{\Cal K}_{\mu \nu }^{\widehat{a}}=\frac{\delta \omega _\nu ^{\widehat{a}}}{%
\delta u^\mu }-\frac{\delta \omega _\mu ^{\widehat{a}}}{\delta u^\nu }+\frac
12f_{\widehat{b}\widehat{c}}^{\quad \widehat{a}}\left( \omega _\mu ^{%
\widehat{b}}\omega _\nu ^{\widehat{c}}-\omega _\nu ^{\widehat{b}}\omega _\mu
^{\widehat{c}}\right) .
$$

This section ends by considering the problem of reduction of the local
an\-i\-sot\-rop\-ic gauge symmetries and gauge fields to isotropic ones. For local
trivial considerations we can consider that the vanishing of dependencies on $%
y$ variables leads to isotropic Yang-Mills fields with the same gauge group
as in the anisotropic case, Global geometric constructions require a more
rigorous topological study of possible obstacles for reduction of total
spaces and structural groups on anisotropic bases to their analogous on
isotropic (for example, pseudo-Riemannian) base spaces.
\vskip25pt

\head{II.2\ Yang--Mills Equations on Locally Anisotropic Spaces}\endhead

Interior gauge (nongravitational) symmetries are associated to semisimple
structural groups. On the principal bundle $\left( P, \pi , Gr, {\Cal E}\right)
$ with nondegenerate Killing form for semisimple group $Gr$ we can define
the generalized Lagrange metric
$$
h_p\left( X_p,Y_p\right) =G_{\pi \left( p\right) }\left( d\pi _PX_P,d\pi
_PY_P\right) +K\left( \Omega _P\left( X_P\right) ,\Omega _P\left( X_P\right)
\right) ,\tag2.10
$$
where $d\pi _P$ is the differential of map $\pi :P\rightarrow {\Cal E},$ $%
G_{\pi \left( p\right) }$ is locally generated as the la-metric (1.12), and $K$
is the Killing form on ${\Cal G}:$%
$$
K\left( \Delta _{
\widehat{a}},\Delta _{\widehat{b}}\right) =f_{\widehat{b}\widehat{d}}^{\quad
\widehat{c}}f_{\widehat{a}\widehat{c}}^{\quad \widehat{d}}=K_{\widehat{a}%
\widehat{b}}.
$$

Using the metric $G_{\alpha \beta }$ on ${\Cal E}$ $\left[ h_P\left(
X_P,Y_P\right) \text{ on }P\right] ,$ we can introduce operators $*_G$ and $%
\widehat{\delta }_G$ acting in the space of forms on ${\Cal E}$ ($*_H$ and $%
\widehat{\delta }_H$ acting on forms on ${\Cal E).}$ Let $e_{\underline{\mu }%
}$ be orthonormalized frames on ${\Cal{U\subset E}}$ and $e^\mu $ the
adjoint coframes. Locally%
$$
G=\sum\limits_\mu \eta \left( \mu \right) e^\mu \otimes e^\mu ,
$$
where $\eta _{\mu \mu }=\eta \left( \mu \right) =\pm 1,$ $\mu
=1,2,...,n,n+1,...,n+m,$ and the Hodge operator $*_{G}$ can be
defined as $*_{G}:\Lambda ^{\prime }\left( {\Cal E}\right)
\rightarrow \Lambda ^{n+m}\left( {\Cal E}\right) ,$ or, in explicit form, as%
$$\multline
*_G\left( e^{\mu _1}\bigwedge ...\bigwedge e^{\mu _r}\right) =\eta \left(
\nu _1\right) ...\eta \left( \nu _{n+m-r}\right) \times\\
sign\left(
\matrix
1 & 2 & ...r & r+1 & ...n+m \\
\mu _1 & \mu _2 & ...\mu _r & \nu _1 & ...\nu _{n+m-r}
\endmatrix
\right) \times e^{\nu _1}\bigwedge ...\bigwedge e^{\nu _{n+m-r}}.
\endmultline \tag2.11
$$
Next, define the operator%
$$
*_{G^{}}^{-1}=\eta \left( 1\right) ...\eta \left( n+m\right) \left(
-1\right) ^{r\left( n+m-r\right) }*_G
$$
and introduce the scalar product on forms $\beta _1,\beta _2,...\subset
\Lambda ^r\left( {\Cal E}\right) $ with compact carrier:%
$$
\left( \beta _1,\beta _2\right) =\eta \left( 1\right) ...\eta \left(
n+m\right) \int \beta _1\bigwedge *_G\beta _2.
$$
The operator $\widehat{\delta }_G$ is defined as the adjoint to $d$
associated to the scalar product for forms, specified for $r$-forms as
$$
\widehat{\delta }_G=\left( -1\right) ^r*_{G^{}}^{-1}\circ d\circ
*_G.\tag2.12
$$

We remark that operators $*_H$ and $\delta _H$ acting in the total space of $%
P$ can be defined similarly to (2.11) and (2.12), but by using metric
(2.10). Both these operators also act in the space of ${\Cal G}$%
-valued forms:%
$$
*\left( \Delta _{\widehat{a}}\otimes \varphi ^{\widehat{a}}\right) =\Delta _{%
\widehat{a}}\otimes (*\varphi ^{\widehat{a}}),\
\widehat{\delta }\left( \Delta _{\widehat{a}}\otimes \varphi ^{\widehat{a}%
}\right) =\Delta _{\widehat{a}}\otimes (\widehat{\delta }\varphi ^{\widehat{a%
}}).
$$

The form $\lambda $ on $P$ with values in ${\Cal G}$ is called horizontal if
$\widehat{H}\lambda =\lambda $ and equivariant if $R^{*}\left( q\right)
\lambda =Ad~q^{-1}\varphi ,~\forall g\in Gr,R\left( q\right) $ being the
right shift on $P.$ We can verify that equivariant and horizontal forms also
satisfy the conditions%
$$
\lambda =Ad~\varphi _{{\Cal U}}^{-1}\left( \pi ^{*}\lambda \right) ,\qquad
\lambda _{{\Cal U}}=S_{{\Cal U}}^{*}\lambda ,\
\left( \lambda _{{\Cal V}}\right) _{{\Cal U}}=Ad~\left( g_{{\Cal{UV}}}\left(
u\right) \right) ^{-1}\left( \lambda _{{\Cal U}}\right) _u.
$$

\proclaim{\bf Definition 2.1}
The field equations for curvature (2.7) and connection (2.1) are defined
by using geometric operators (2.11) and (2.12):
$$
\Delta {\Cal K}=0, \tag2.13
$$
$$
\nabla {\Cal K}=0,\tag2.14$$
where $\Delta =\widehat{H}\circ \widehat{\delta }_H.$
\endproclaim

Equations (2.13) are
similar to the well-known Maxwell equations and for non-Abelian gauge fields
are called Yang-Mills equations. The structural equations (2.14) are called
Bianchi identities.

The field equations (2.13) do not have a physical meaning because they are
written in the total space of bundle $\Upsilon $ and not on the base
anisotropic space-time ${\Cal E}.$ But this difficulty may be obviated by
projecting the mentioned equations on the base. The 1-form $\Delta {\Cal K}$
is horizontal by definition and its equivariance follows from the right
invariance of metric (2.10). So, there is a unique form $(\Delta {\Cal K})_{%
{\Cal U}}$ satisfying
$$
\Delta {\Cal K=}Ad~\varphi _{{\Cal U}}^{-1}\pi ^{*}(\Delta {\Cal K})_{{\Cal U%
}}.
$$
Projection of (2.13) on the base can be written as $(\Delta {\Cal K})_{{\Cal %
U}}=0.$ To calculate $(\Delta {\Cal K})_{{\Cal U}},$ we use
the equality \cite{Bishop and Crittenden 1964} and
\cite{Popov and Dikhin 1975}
$$
d\left( Ad~\varphi _{{\Cal U}}^{-1}\lambda \right) =Ad~~\varphi _{{\Cal U}%
}^{-1}~d\lambda -\left[ \varphi _{{\Cal U}}^{*}\theta ,Ad~\varphi _{{\Cal U}%
}^{-1}\lambda \right] \tag2.15
$$
where $\lambda $ is a form on $P$ with values in ${\Cal G}.$ For r--forms we
have
$$
\widehat{\delta }\left( Ad~\varphi _{{\Cal U}}^{-1}\lambda \right)
=Ad~\varphi _{{\Cal U}}^{-1}\widehat{\delta }\lambda -\left( -1\right)
^r*_H\{\left[ \varphi _{{\Cal U}}^{*}\theta ,*_HAd~\varphi _{{\Cal U}%
}^{-1}\lambda \right]
$$
and, as a consequence,%
$$\multline
\widehat{\delta }{\Cal K}=Ad~\varphi _{{\Cal U}}^{-1}\{\widehat{\delta }%
_H\pi ^{*}{\Cal K}_{{\Cal U}}+*_H^{-1}[\pi ^{*}\omega _{{\Cal U}},*_H\pi ^{*}%
{\Cal K}_{{\Cal U}}]\}- \\
*_H^{-1}\left[ \Omega ,Ad~\varphi _{{\Cal U}}^{-1}*_H\left( \pi ^{*}{\Cal K}%
\right) \right] .\endmultline \tag2.16
$$
By using straightforward calculations in a locally adapted dual basis on $\pi
^{-1}\left( {\Cal U}\right) $ we can verify the equalities%
$$\split
\left[ \Omega ,Ad~\varphi _{{\Cal U}}^{-1}~*_H\left( \pi ^{*}{\Cal K}_{{\Cal %
U}}\right) \right] =0,\widehat{H}\delta _H\left( \pi ^{*}{\Cal K}_{{\Cal U}%
}\right) =\pi ^{*}\left( \widehat{\delta }_G{\Cal K}\right) ,\\
*_H^{-1}\left[ \pi ^{*}\omega _{{\Cal U}},*_H\left( \pi ^{*}{\Cal K}_{{\Cal U%
}}\right) \right] =\pi ^{*}\{*_G^{-1}\left[ \omega _{{\Cal U}},*_G{\Cal K}_{%
{\Cal U}}\right] \}. \endsplit \tag2.17
$$
From (2.16) and (2.17) it follows that
$$
\left( \Delta {\Cal K}\right) _{{\Cal U}}=\widehat{\delta }_G{\Cal K}_{{\Cal %
U}}+*_G^{-1}\left[ \omega _{{\Cal U}},*_G{\Cal K}_{{\Cal U}}\right].\tag2.18
$$

Taking into account (2.18) and (2.12), we prove that projection on ${\Cal E}$
of equations (2.13) and (2.14) can be expressed respectively as
$$
*_G^{-1}\circ d\circ *_G{\Cal K}_{{\Cal U}}+*_G^{-1}\left[ \omega _{{\Cal U}%
},*_G{\Cal K}_{{\Cal U}}\right] =0.\tag2.19
$$
$$
d{\Cal K}_{{\Cal U}}+\left[ \omega _{{\Cal U}},{\Cal K}_{{\Cal U}}\right]
=0.\tag2.20
$$

Equations (2.19) (see (2.18)) are gauge-invariant because%
$$
\left( \Delta {\Cal K}\right) _{{\Cal U}}=Ad~g_{{\Cal{UV}}}^{-1}\left( \Delta
{\Cal K}\right) _{{\Cal V}}.
$$

By using formulas (2.9)-(2.12) we can rewrite (2.19) in coordinate form%
$$
D_\nu \left( G^{\nu \lambda }{\Cal K}_{~\lambda \mu }^{\widehat{a}}\right)
+f_{\widehat{b}\widehat{c}}^{\quad \widehat{a}}G^{v\lambda }\omega _\lambda
^{~\widehat{b}}{\Cal K}_{~\nu \mu }^{\widehat{c}}=0,\tag2.21$$
where $D_\nu $ is, for simplicity, a compatible with metric covariant
derivation on la-space ${\Cal E}.$

It is possible to distinguish the $x$ and $y$ parts of equations (2.21) by
using formulas (1.4),(1.5),(1.12),(1.17),(1.18),(1.28)(1.32). We omit this
trivial calculus.

We point out that for our bundles with semisimple structural groups the
Yang-Mills equations (2.13) (and, as a consequence, their horizontal
projections (2.19) or (2.21)) can be obtained by variation of the action%
$$
I=\int {\Cal K}_{~\mu \nu }^{\widehat{a}}{\Cal K}_{~\alpha \beta }^{\widehat{%
b}}G^{\mu \alpha }G^{\nu \beta }K_{\widehat{a}\widehat{b}}\left| G_{\alpha
\beta }\right| ^{1/2}dx^1...dx^n\delta y^1...\delta y^m.\tag2.22
$$
Equations for extremals of (2.22) have the form
$$
K_{\widehat{r}\widehat{b}}G^{\lambda \alpha }G^{\kappa \beta }D_\alpha {\Cal %
K}_{~\lambda \beta }^{\widehat{b}}-K_{\widehat{a}\widehat{b}}G^{\kappa
\alpha }G^{\nu \beta }f_{\widehat{r}\widehat{l}}^{\quad \widehat{a}}\omega
_\nu ^{\widehat{l}}{\Cal K}_{~\alpha \beta }^{\widehat{b}}=0,
\tag2.23
$$
which are equivalent to ''pure'' geometric equations (2.21) (or (2.19)) due
to nondegeneration of the Killing form $K_{\widehat{r}\widehat{b}}$ for
semisimple groups.

To take into account gauge interactions with matter fields (sections of
vector bundle $\Upsilon $ on ${\Cal E}$ ) we have to introduce a source
1-form ${\Cal J}$ in equations (2.13) and to write them as%
$$
\Delta {\Cal K}={\Cal J}\tag2.24
$$

Explicit constructions of ${\Cal J}$ require concrete definitions of the
bundle $\Upsilon ;$ for example, for spinor fields an invariant formulation
of the Dirac equations on la-spaces is necessary. We omit spinor
considerations in this Chapter (see section 1.7),
 but we shall present the definition
of the source ${\Cal J}$ for gravitational interactions (by using the
energy--momentum tensor of matter on la--space) in the next section.
\vskip25pt

\head{II.3 Gauge Locally Anisotropic Gravity } \endhead

A considerable body of work on the formulation of gauge gravitational models
on isotropic spaces is based on using nonsemisimple groups, for example,
Poincare and affine groups, as structural gauge groups (see critical
analysis and original results in \cite{Tseytlin 1982} and
\cite{Ponomarev, Barvinsky and Obukhov 1985}). The main
impediment to developing such models is caused by the degeneration of
Killing forms for nonsemisimple groups, which make it impossible to
construct consistent variational gauge field theories (functional (2.22) and
extremal equations are degenarate in these cases). There are at least two
possibilities to get around the mentioned difficulty.\ The first is to
realize a minimal extension of the nonsemisimple group to a semisimple one,
similar to the extension of the Poincare group to the de Sitter group
(in the next section we
shall use this operation for the definition of locally anisotropic
gravitational instantons). The second possibility is to introduce into
consideration the bundle of adapted affine frames on la--space ${\Cal E},$ to
use an auxiliary nondegenerate bilinear form $a_{\widehat{a}\widehat{b}}$
instead of the degenerate Killing form $K_{\widehat{a}\widehat{b}}$ and to
consider a ''pure'' geometric method, illustrated in the previous section,
of defining gauge field equations. Projecting on the base ${\Cal E},$ we
shall obtain gauge gravitational field equations on la--space having a form
similar to Yang-Mills equations.

The goal of this section is to prove that a specific parametrization of
components of the Cartan connection in the bundle of adapted affine frames on
${\Cal E}$ establishes an equivalence between Yang-Mills equations (2.24)
and Einstein equations on la--spaces.
\vskip15pt

\subhead{\bf II.3.1 The bundle of linear locally adapted frames }
    \endsubhead

Let $\left( X_\alpha \right) _u=\left( X_i,X_a\right) _u$ be an adapted
frame (see (1.4)  at point $u\in {\Cal E}.$ We consider a local right
distinguished action of matrices
$$
A_{\alpha ^{\prime }}^{\quad \alpha }=\left(
\matrix
A_{i^{\prime }}^{\quad i} & 0 \\
0 & B_{a^{\prime }}^{\quad a}
\endmatrix
\right) \subset GL_{n+m}=GL\left( n,{\Bbb R}\right) \oplus GL\left( m,{\Bbb R%
}\right) .0
$$
Nondegenerate matrices $A_{i^{\prime }}^{\quad i}$ and $B_{j^{\prime
}}^{\quad j}$ respectively transforms linearly $X_{i|u}$ into $X_{i^{\prime
}|u}=A_{i^{\prime }}^{\quad i}X_{i|u}$ and $X_{a^{\prime }|u}$ into $%
X_{a^{\prime }|u}=B_{a^{\prime }}^{\quad a}X_{a|u},$ where $X_{\alpha
^{\prime }|u}=A_{\alpha ^{\prime }}^{\quad \alpha }X_\alpha $ is also an
adapted frame at the same point $u\in {\Cal E}.$ We denote by $La\left(
{\Cal E}\right) $ the set of all adapted frames $X_\alpha $ at all points of
${\Cal E}$ and consider the surjective map $\pi $ from $La\left( {\Cal E}%
\right) $ to ${\Cal E}$ transforming every adapted frame $X_{\alpha |u}$ and
point $u$ into point $u.$ Every $X_{\alpha ^{\prime }|u}$ has a unique
representation as $X_{\alpha ^{\prime }}=A_{\alpha ^{\prime }}^{\quad \alpha
}X_\alpha ^{\left( 0\right) },$ where $X_\alpha ^{\left( 0\right) }$ is a
fixed distinguished basis in tangent space $T\left( {\Cal E}\right) .$ It is
obvious that $\pi ^{-1}\left( {\Cal U}\right) ,{\Cal U}\subset {\Cal E},$ is
bijective to ${\Cal U}\times GL_{n+m}\left( {\Bbb R}\right) .$ We can
transform $La\left( {\Cal E}\right) $ in a differentiable manifold taking $%
\left( u^\beta ,A_{\alpha ^{\prime }}^{\quad \alpha }\right) $ as a local
coordinate system on $\pi ^{-1}\left( {\Cal U}\right) .$ Now, it is easy to
verify that ${\Cal{L}}a ({\Cal E}) = ( La ( {\Cal E}, {\Cal E},
GL_{n+m} ( {\Bbb R} ) ) ) $ is a principal bundle. We call
${\Cal{L}}a ( {\Cal E} ) $ the bundle of linear adapted frames
on ${\Cal E}.$

The next step is to identify the components of, for simplicity, compatible
d-connection $\Gamma _{\beta \gamma }^\alpha $ on ${\Cal E}:$%
$$
\Omega _{{\Cal U}}^{\widehat{a}}=\omega ^{\widehat{a}}=\{\omega _{\quad
\lambda }^{\widehat{\alpha }\widehat{\beta }}\doteq \Gamma _{\beta \gamma
}^\alpha \}.\tag2.25
$$
Introducing (2.25) in (2.18), we calculate the local 1-form%
$$
\left( \Delta {\Cal R}^{(\Gamma )}
\right) _{{\Cal U}}=\Delta _{\widehat{\alpha }\widehat{%
\alpha }_1}\otimes \left( G^{\nu \lambda }D_\lambda {\Cal R}_{\qquad \nu \mu
}^{\widehat{\alpha }\widehat{\alpha }_1}+f_{\qquad \widehat{\beta }\widehat{%
\beta }_1\widehat{\gamma }\widehat{\gamma }_1}^{\widehat{\alpha }\widehat{%
\alpha }_1}G^{\nu \lambda }\omega _{\qquad \lambda }^{\widehat{\beta }%
\widehat{\beta }_1}{\Cal R}_{\qquad \nu \mu }^{\widehat{\gamma }\widehat{%
\gamma }_1}\right) \delta u^\mu ,\tag2.26
$$
where
$$
\Delta _{\widehat{\alpha }\widehat{\alpha }_1}=\left(
\matrix
\Delta _{ii_1} & 0 \\
0 & \Delta _{aa_1}
\endmatrix
\right)
$$
is the standard distinguished basis in Lie algebra of matrices
${{\Cal{G}}l}_{n+m}\left( {\Bbb R}\right) $ with $\left( \Delta _{ii_1}\right)
_{jj_1}=\delta _{ij}\delta _{i_1j_1}$ and $\left( \Delta _{aa_1}\right)
_{bb_1}$ be\-ing res\-pec\-ti\-ve\-ly the stand\-ard bas\-es in \newline
  ${\Cal{G}}l\left({\Bbb R}^{n+m}\right) .$ We have denoted the curvature
of connection (2.25), considered in (2.26), as
$$
{\Cal R}_{{\Cal U}}^{(\Gamma )} =
\Delta _{\widehat{\alpha }\widehat{\alpha }_1}\otimes
{\Cal R}_{\qquad \nu \mu }^{\widehat{\alpha }\widehat{\alpha }_1}X^\nu
\bigwedge X^\mu ,\tag2.27
$$
where ${\Cal R}_{\qquad \nu \mu }^{\widehat{\alpha }\widehat{\alpha }%
_1}=R_{\alpha _1 \quad \nu \mu }^{\quad \alpha} $ (see curvatures (1.32)).
\vskip15pt

\subhead{\bf II.3.2 The bundle of  affine locally adapted frames }\endsubhead

Besides ${\Cal{L}}a\left( {\Cal E}\right) $ with la--space ${\Cal E},$
another bundle is naturally related, the bundle of adapted affine frames with
structural group $Af_{n+m}\left( {\Bbb R}\right) =GL_{n+m}\left( {\Cal E}%
\right) \otimes {\Bbb R}^{n+m}.$ Because as linear space the Lie Algebra $%
af_{n+m}\left( {\Bbb R}\right) $ is a direct sum of ${{\Cal{G}}l}_{n+m}\left(
{\Bbb R}\right) $ and ${\Bbb R}^{n+m},$ we can write forms on ${\Cal{A}}a%
\left( {\Cal E}\right) $ as $\Theta =\left( \Theta _1,\Theta _2\right) ,$
where $\Theta _1$ is the ${{\Cal{G}}l}_{n+m}\left( {\Bbb R}\right) $ component
and $\Theta _2$ is the ${\Bbb R}^{n+m}$ component of the form $\Theta .$
Connection (2.25), $\Omega $ in ${{\Cal{L}}a}\left( {\Cal E}\right) ,$ induces
the Cartan connection $\overline{\Omega }$ in ${{\Cal{A}}a}\left( {\Cal E}%
\right) ;$ see the isotropic case in \cite{Bishop and Crittenden 1964} and
 \cite{Popov and Dikhin 1975}. This
is the unique connection on ${{\Cal{A}}a}\left( {\Cal E}\right) $ represented
as $i^{*}\overline{\Omega }=\left( \Omega ,\chi \right) ,$ where $\chi $ is
the shifting form and $i:{{\Cal{A}}a}\rightarrow {{\Cal{L}}a}$
is the trivial reduction
of bundles. If $s_{{\Cal U}}^{(a)}$ is a local adapted frame in ${{\Cal{L}}a}%
\left( {\Cal E}\right) ,$ then $\overline{s}_{{\Cal U}}^{\left( 0\right)
}=i\circ s_{{\Cal U}}$ is a local section in
${{\Cal{A}}a}\left( {\Cal E}\right)$ and
$$
\left( \overline{\Omega }_{{\Cal U}}\right) =s_{{\Cal U}}\Omega =\left(
\Omega _{{\Cal U}},\chi _{{\Cal U}}\right) ,\tag2.28
$$
$$
\left( \overline{{\Cal R}}_{{\Cal U}}\right) =s_{{\Cal U}}\overline{{\Cal R}}%
=\left( {\Cal R}_{{\Cal U}}^{(\Gamma )},T_{{\Cal U}}\right) ,\tag2.29
$$
where $\chi =e_{\widehat{\alpha }}\otimes \chi _{\quad \mu }^{\widehat{%
\alpha }}X^\mu ,G_{\alpha \beta }=\chi _{\quad \alpha }^{\widehat{\alpha }%
}\chi _{\quad \beta }^{\widehat{\beta }}\eta _{\widehat{\alpha }\widehat{%
\beta }}\quad (\eta _{\widehat{\alpha }\widehat{\beta }}$ is diagonal with $%
\eta _{\widehat{\alpha }\widehat{\alpha }}=\pm 1)$ is a frame decomposition
of metric (1.12) on ${\Cal E},e_{\widehat{\alpha }}$ is the standard
distinguished basis on ${\Bbb R}^{n+m},$ and the projection of torsion , $T_{%
{\Cal U}},$ on base ${\Cal E}$ is defined as
$$
T_{{\Cal U}}=d\chi _{{\Cal U}}+\Omega _{{\Cal U}}\bigwedge \chi _{{\Cal U}%
}+\chi _{{\Cal U}}\bigwedge \Omega _{{\Cal U}}=e_{\widehat{\alpha }}\otimes
\sum\limits_{\mu <\nu }T_{\quad \mu \nu }^{\widehat{\alpha }}X^\mu \bigwedge
X^\nu .\tag2.30
$$
For a fixed local adapted basis on ${\Cal U}\subset {\Cal E}$ we can identify
components $T_{\quad \mu \nu }^{\widehat{a}}$ of torsion (2.30) with
components of torsion (1.28) on ${\Cal E},$ i.e.
 $T_{\quad \mu \nu }^{\widehat{\alpha }}=T_{\quad \mu \nu }^\alpha .$
By straightforward calculation we obtain
$$
{( \Delta \overline{{\Cal R}})}_{{\Cal U}}
= [ {( \Delta{\Cal R}^{(\Gamma )}) }_{{\Cal U}},\
 {( R\tau )}_{{\Cal U}} + {( Ri )}_{{\Cal U}} ] ,
\tag2.31
$$
where%
$$
\left( R\tau \right) _{{\Cal U}}=\widehat{\delta }_GT_{{\Cal U}%
}+*_G^{-1}\left[ \Omega _{{\Cal U}},*_GT_{{\Cal U}}\right] ,\quad \left(
Ri\right) _{{\Cal U}}=*_G^{-1}\left[ \chi _{{\Cal U}},*_G{\Cal R}_{{\Cal U}%
}^{(\Gamma )}\right] .
$$
Form $\left( Ri\right) _{{\Cal U}}$ from (2.31) is locally constructed by
using components of the Ricci tensor (see (1.33)) as follows from decomposition
on the local adapted basis $X^\mu =\delta u^\mu :$
$$
\left( Ri\right) _{{\Cal U}}=e_{\widehat{\alpha }}\otimes \left( -1\right)
^{n+m+1}R_{\lambda \nu }G^{\widehat{\alpha }\lambda }\delta u^\mu
\tag2.32
$$

We remark that for isotropic torsionless pseudo-Riemannian spaces the
requirement that $\left( \Delta \overline{{\Cal R}}\right) _{{\Cal U}}=0,$
i.e., imposing the connection (2.25) to satisfy Yang-Mills equations (2.13)
(equivalently (2.19) or (2.21) we obtain
the equivalence of the mentioned gauge gravitational equations with the
vacuum Einstein equations $R_{ij}=0.\,$ In the case of la--spaces with
arbitrary given torsion, even considering vacuum gravitational fields, we
have to introduce a source for gauge gravitational equations in order to
compensate for the contribution of torsion and to obtain equivalence with
the Einstein equations.

The above presented considerations  constitute the proof of
the following

\proclaim{\bf Theorem II.1}
The Einstein equations (1.42) for locally anisotropic gravity are equivalent to
Yang-Mills equations%
$$
\left( \Delta \overline{{\Cal R}}\right) =\overline{{\Cal J}}
\tag2.33
$$
for the induced Cartan connection $\overline{\Omega }$ (see (2.25), (2.28))
in the bundle of local adapted affine frames ${\Cal A}a\left( {\Cal E}\right)
$ with source $\overline{{\Cal J}}_{{\Cal U}}$ constructed locally by using
the same formulas (2.31) (for $\left( \Delta \overline{{\Cal R}}\right) $),
where $R_{\alpha \beta }$ is changed by the matter source
 ${\tilde E}_{\alpha \beta}-\frac 12G_{\alpha \beta }{\tilde E},$ where
${\tilde E}_{\alpha\beta}=kE_{\alpha\beta} - \lambda G_{\alpha\beta}.$
\endproclaim
\vskip25pt

\head{II.4 Nonlinear De Sitter Gauge Locally Anisotropic Gravity}\endhead

The equivalent reexpression of the Einstein theory as a gauge like theory
implies, for both locally isotropic and anisotropic space-times, the
nonsemisimplicity of the gauge group, which leads to a nonvariational theory
in the total space of the bundle of locally adapted affine frames. A
variational gauge gravitational theory can be formulated by using a minimal
extension of the affine structural group ${{\Cal A}f}_{n+m}\left( {\Bbb R}%
\right) $ to the de Sitter gauge group $S_{n+m}=SO\left( n+m+1\right) $
acting on distinguished ${\Bbb R}^{n+m+1}$ space.
\vskip15pt

\subhead{II.4.1 Nonlinear gauge theories of de Sitter group}\endsubhead

Let us consider the de Sitter space $\Sigma ^{n+m}$ as a
hypersurface given by the equations $\eta _{AB}u^Au^B=-l^2$ in the
(n+m)-dimensional spaces enabled with diagonal metric $\eta _{AB},\eta
_{AA}=\pm 1$ (in this section $A,B,C,...=1,2,...,n+m+1),$ where $\{u^A\}$
are global Cartesian coordinates in ${\Bbb R}^{n+m+1};l>0$ is the curvature
of de Sitter space. The de Sitter group $S_{\left( \eta \right) }=SO_{\left(
\eta \right) }\left( n+m+1\right) $ is defined as the isometry group of $%
\Sigma ^{n+m}$-space with $\frac{n+m}2\left( n+m+1\right) $ generators of
Lie algebra ${{\Cal s}o}_{\left( \eta \right) }\left( n+m+1\right) $
satisfying the commutation relations%
$$
\left[ M_{AB},M_{CD}\right] =\eta _{AC}M_{BD}-\eta _{BC}M_{AD}-\eta
_{AD}M_{BC}+\eta _{BD}M_{AC}.\tag2.34
$$

Decomposing indices $A,B,...$ as $A=\left( \widehat{\alpha },n+m+1\right)
, B=\left( \widehat{\beta },n+m+1\right) ,$  $ ...,$ the metric $\eta _{AB}$ as $%
\eta _{AB}=\left( \eta _{\widehat{\alpha }\widehat{\beta }},\eta _{\left(
n+m+1\right) \left( n+m+1\right) }\right) ,$ and operators $M_{AB}$ as $M_{%
\widehat{\alpha }\widehat{\beta }}={\Cal F}_{\widehat{\alpha }\widehat{\beta
}}$ and $P_{\widehat{\alpha }}=l^{-1}M_{n+m+1,\widehat{\alpha }},$ we can
write (2.34) as
$$
\left[ {\Cal F}_{\widehat{\alpha }\widehat{\beta }},{\Cal F}_{\widehat{%
\gamma }\widehat{\delta }}\right] =\eta _{\widehat{\alpha }\widehat{\gamma }}%
{\Cal F}_{\widehat{\beta }\widehat{\delta }}-\eta _{\widehat{\beta }\widehat{%
\gamma }}{\Cal F}_{\widehat{\alpha }\widehat{\delta }}+\eta _{\widehat{\beta
}\widehat{\delta }}{\Cal F}_{\widehat{\alpha }\widehat{\gamma }}-\eta _{%
\widehat{\alpha }\widehat{\delta }}{\Cal F}_{\widehat{\beta }\widehat{\gamma
}},
$$
$$
\left[ P_{\widehat{\alpha }},P_{\widehat{\beta }}\right] =-l^{-2}{\Cal F}_{%
\widehat{\alpha }\widehat{\beta }},\quad \left[ P_{\widehat{\alpha }},{\Cal F%
}_{\widehat{\beta }\widehat{\gamma }}\right] =\eta _{\widehat{\alpha }%
\widehat{\beta }}P_{\widehat{\gamma }}-\eta _{\widehat{\alpha }\widehat{%
\gamma }}P_{\widehat{\beta }},
$$
where we have indicated the possibility to decompose ${{\Cal s}o}_{\left( \eta
\right) }\left( n+m+1\right) $ into a direct sum, ${{\Cal s}o}_{\left( \eta
\right) }\left( n+m+1\right) ={{\Cal s}o}_{\left( \eta \right) }(n+m)\oplus
V_{n+m},$ where $V_{n+m}$ is the vector space stretched on vectors $P_{%
\widehat{\alpha }}.$ We remark that $\Sigma ^{n+m}=S_{\left( \eta \right)
}/L_{\left( \eta \right) },$ where $L_{\left( \eta \right) }=SO_{\left( \eta
\right) }\left( n+m\right) .$ For $\eta _{AB}=diag\left( 1,-1,-1,-1\right) $
and $S_{10}=SO\left( 1,4\right) ,L_6=SO\left( 1,3\right) $ is the group of
Lorentz rotations.

Let $W\left( {\Cal E},{\Bbb R}^{n+m+1},S_{\left( \eta \right) },P\right) $
be the
vector bundle associated with principal bundle $P\left( S_{\left( \eta
\right) },{\Cal E}\right) $ on la--spaces. The action of the structural group
$S_{\left( \eta \right) }$ on $E$ can be realized by using $\left(
n+m\right) \times \left( n+m\right) $ matrices with a parametrization
distinguishing subgroup $L_{\left( \eta \right) }:$%
$$
B=bB_L,\ \text{ where }
B_L=\left(
\matrix
L & 0 \\
0 & 1
\endmatrix
\right) , \tag2.35
$$
$L\in L_{\left( \eta \right) }$ is the de Sitter bust matrix transforming
the vector $\left( 0,0,...,\rho \right) \in {\Bbb R}^{n+m+1}$ into the
arbitrary point $\left( V^1,V^2,...,V^{n+m+1}\right) \in \Sigma _\rho
^{n+m}\subset {\Bbb R}^{n+m+1}$ with curvature $\rho \quad \left(
V_AV^A=-\rho ^2,V^A=t^A\rho \right) .$ Matrix $b$ can be expressed as
$$
b=\left(
\matrix
\delta _{\quad \widehat{\beta }}^{\widehat{\alpha }}+\frac{t^{\widehat{%
\alpha }}t_{\widehat{\beta }}}{\left( 1+t^{n+m+1}\right) } & t^{
\widehat{\alpha }} \\ t_{\widehat{\beta }} & t^{n+m+1}
\endmatrix
\right) .
$$

The de Sitter gauge field is associated with a linear connection in $W$,
i.e., with a ${{\Cal s}o}_{\left( \eta \right) }\left( n+m+1\right) $-valued
connection 1-form on ${\Cal E}:$

$$
\widetilde{\Omega }=\left(
\matrix
\omega _{\quad \widehat{\beta }}^{\widehat{\alpha }} & \widetilde{\theta }^{%
\widehat{\alpha }} \\ \widetilde{\theta }_{\widehat{\beta }} & 0
\endmatrix \right) , \tag2.36
$$
where
 $\omega _{\quad \widehat{\beta }}^{\widehat{\alpha }}\in  so(n+m)%
_{\left( \eta \right) },$ $\widetilde{\theta }^{\widehat{\alpha }}\in {\Bbb R%
}^{n+m},\widetilde{\theta }_{\widehat{\beta }}\in \eta _{\widehat{\beta }%
\widehat{\alpha }}\widetilde{\theta }^{\widehat{\alpha }}.$

Because $S_{\left( \eta \right) }$-transforms mix $\omega _{\quad \widehat{%
\beta }}^{\widehat{\alpha }}$ and $\widetilde{\theta }^{\widehat{\alpha }}$
fields in (2.36) (the introduced para\-met\-ri\-za\-ti\-on is
invariant on action on $%
SO_{\left( \eta \right) }\left( n+m\right) $ group we cannot identify $%
\omega _{\quad \widehat{\beta }}^{\widehat{\alpha }}$ and $\widetilde{\theta
}^{\widehat{\alpha }},$ respectively, with the connection $\Gamma _{~\beta
\gamma }^\alpha $ and the fundamental form $\chi ^\alpha $ in ${\Cal E}$ (as
we have for (2.25) and (2.28)). To avoid this difficulty we consider
\cite{Tseytlin 1982} a nonlinear gauge realization of the de
Sitter group $S_{\left( \eta \right) },$ namely, we introduce into
consideration the nonlinear gauge field%
$$
\Omega =b^{-1}\Omega b+b^{-1}db=\left(
\matrix
\Gamma _{~\widehat{\beta }}^{\widehat{\alpha }} & \theta ^{
\widehat{\alpha }} \\ \theta _{\widehat{\beta }} & 0
\endmatrix
\right) ,\tag2.37
$$
where
$$
\Gamma _{\quad \widehat{\beta }}^{\widehat{\alpha }}=\omega _{\quad \widehat{%
\beta }}^{\widehat{\alpha }}-\left( t^{\widehat{\alpha }}Dt_{\widehat{\beta }%
}-t_{\widehat{\beta }}Dt^{\widehat{\alpha }}\right) /\left(
1+t^{n+m+1}\right) ,
$$
$$
\theta ^{\widehat{\alpha }}=t^{n+m+1}\widetilde{\theta }^{\widehat{\alpha }%
}+Dt^{\widehat{\alpha }}-t^{\widehat{\alpha }}\left( dt^{n+m+1}+\widetilde{%
\theta }_{\widehat{\gamma }}t^{\widehat{\gamma }}\right) /\left(
1+t^{n+m+1}\right) ,
$$
$$
Dt^{\widehat{\alpha }}=dt^{\widehat{\alpha }}+\omega _{\quad \widehat{\beta }%
}^{\widehat{\alpha }}t^{\widehat{\beta }}.
$$

The action of the group $S\left( \eta \right) $ is nonlinear, yielding
transforms  $\Gamma ^{\prime }=L^{\prime }\Gamma \left( L^{\prime }\right)
^{-1}+L^{\prime }d\left( L^{\prime }\right) ^{-1},\theta ^{\prime }=L\theta
, $ where the nonlinear matrix-valued function  $L^{\prime }=L^{\prime
}\left( t^\alpha ,b,B_T\right) $ is defined from $B_b=b^{\prime
}B_{L^{\prime }}$ (see parametrization (2.35)).

Now, we can identify components of (2.37) with components of $\Gamma
_{~\beta \gamma }^\alpha $ and $\chi _{\quad \alpha }^{\widehat{\alpha }}$
on ${\Cal E}$ and induce in a consistent manner on the base of bundle  $%
W\left( {\Cal E},{\Bbb R}^{n+m+1},S_{\left( \eta \right) },P\right) $ the
la--geometry.
\vskip15pt

\subhead{II.4.2\ Dynamics of the nonlinear locally anisotropic De Sitter
gravity}\endsubhead

Instead of the gravitational potential (2.25), we introduce the
gravitational connection (similar to (2.37))
$$
\Gamma =\left(
\matrix
\Gamma _{\quad \widehat{\beta }}^{\widehat{\alpha }} & l_0^{-1}\chi ^{
\widehat{\alpha }} \\ l_0^{-1}\chi _{\widehat{\beta }} & 0
\endmatrix
\right) \tag2.38
$$
where $$\Gamma _{\quad \widehat{\beta }}^{\widehat{\alpha }}=\Gamma _{\quad
\widehat{\beta }\mu }^{\widehat{\alpha }}\delta u^\mu ,\Gamma _{\quad
\widehat{\beta }\mu }^{\widehat{\alpha }}=\chi _{\quad \alpha }^{\widehat{%
\alpha }}\chi _{\quad \beta }^{\widehat{\beta }}\Gamma _{\quad \beta \gamma
}^\alpha +\chi _{\quad \alpha }^{\widehat{\alpha }}\delta _\mu \chi _{\quad
\widehat{\beta }}^\alpha ,\chi ^{\widehat{\alpha }}=\chi _{\quad \mu }^{%
\widehat{\alpha }}\delta u^\mu ,$$
 $G_{\alpha \beta }=\chi _{\quad \alpha
}^{\widehat{\alpha }}\chi _{\quad \beta }^{\widehat{\beta }}\eta _{\widehat{%
\alpha }\widehat{\beta }},$  $\eta _{\widehat{\alpha }\widehat{\beta }}$
is parametrized as
$
\eta _{\widehat{\alpha }\widehat{\beta }}=\left(
\smallmatrix
\eta _{ij} & 0 \\
0 & \eta _{ab}
\endsmallmatrix
\right) , \eta _{ij}=\left( 1,-1,...,-1\right) ,
$ and $l_0$ is a dimensional constant.

The curvature of (2.39), ${\Cal R}^{(\Gamma )} =
d\Gamma +\Gamma \bigwedge \Gamma ,$ can
be written as%
$$
{\Cal R}^{(\Gamma )} = \left(
\matrix
{\Cal R}_{\quad \widehat{\beta }}^{\widehat{\alpha }}+l_0^{-1}\pi _{\widehat{%
\beta }}^{\widehat{\alpha }} & l_0^{-1}T^{
\widehat{\alpha }} \\ l_0^{-1}T^{\widehat{\beta }} & 0
\endmatrix
\right) ,\tag2.39
$$
where $$\pi _{\widehat{\beta }}^{\widehat{\alpha }} = \chi ^{\widehat{\alpha }%
}\bigwedge \chi _{\widehat{\beta }},{\Cal R}_{\quad \widehat{\beta }}^{%
\widehat{\alpha }}=\frac 12{\Cal R}_{\quad \widehat{\beta }\mu \nu }^{%
\widehat{\alpha }}\delta u^\mu \bigwedge \delta u^\nu ,$$ and $${\Cal R}%
_{\quad \widehat{\beta }\mu \nu }^{\widehat{\alpha }}=\chi _{\widehat{\beta }%
}^{\quad \beta }\chi _\alpha ^{\quad \widehat{\alpha }}R_{\quad \beta _{\mu
\nu }}^\alpha $$ (see (1.31) and (1.32),
the components of d-curvatures). The de Sitter
gauge group is semisimple and we are able to construct a variational gauge
gravitational locally anisotropic theory (bundle metric (2.10) is
nondegenerate). The Lagrangian of the theory is postulated as
$$
L=L_{\left( G\right) }+L_{\left( m\right) }
$$
where the gauge gravitational Lagrangian is defined as
$$
L_{\left( G\right) }=\frac 1{4\pi }Tr\left( {\Cal R}^{(\Gamma )}
\bigwedge *_G{\Cal R}^{(\Gamma )}%
\right) ={\Cal L}_{\left( G\right) }\left| G\right| ^{1/2}\delta ^{n+m}u,
$$
$$
{\Cal L}_{\left( G\right) }=\frac 1{2l^2}T_{\quad \mu \nu }^{\widehat{\alpha
}}T_{\widehat{\alpha }}^{\quad \mu \nu }+\frac 1{8\lambda }{\Cal R}_{\quad
\widehat{\beta }\mu \nu }^{\widehat{\alpha }}{\Cal R}_{\quad \widehat{\alpha
}}^{\widehat{\beta }\quad \mu \nu }-
\frac 1{l^2}\left({\overleftarrow{R}}\left( \Gamma
\right) -2\lambda _1\right) ,\tag2.40
$$
$T_{\quad \mu \nu }^{\widehat{\alpha }}=\chi _{\quad \alpha }^{\widehat{%
\alpha }}T_{\quad \mu \nu }^\alpha $ (the gravitational constant $l^{2}$ in
(2.40) satisfies the relations $l^2=2l_0^2\lambda ,\lambda
_1=-3/l_0]),\quad Tr$ denotes the trace on $\widehat{\alpha },\widehat{\beta }
$ indices, and the matter field Lagrangian is defined as%
$$
L_{\left( m\right) }=-1\frac 12Tr\left( \Gamma \bigwedge *_G{\Cal I}\right) =%
{\Cal L}_{\left( m\right) }\left| G\right| ^{1/2}\delta ^{n+m}u,
$$
$$
{\Cal L}_{\left( m\right) }=\frac 12\Gamma _{\quad \widehat{\beta }\mu }^{%
\widehat{\alpha }}S_{\quad \alpha }^{\widehat{\beta }\quad \mu }-t_{\quad
\widehat{\alpha }}^\mu l_{\quad \mu }^{\widehat{\alpha }}.
\tag2.41
$$
The matter field source ${\Cal I}$ is obtained as a variational derivation
of ${\Cal L}_{\left( m\right) }$ on $\Gamma $ and is parametrized as
$$
{\Cal I}=\left(
\matrix
S_{\quad \widehat{\beta }}^{\widehat{\alpha }} & -l_0t^{
\widehat{\alpha }} \\ -l_0t_{\widehat{\beta }} & 0
\endmatrix
\right) \tag2.42
$$
with $t^{\widehat{\alpha }}=t_{\quad \mu }^{\widehat{\alpha }}\delta u^\mu $
and $S_{\quad \widehat{\beta }}^{\widehat{\alpha }}=S_{\quad \widehat{\beta }%
\mu }^{\widehat{\alpha }}\delta u^\mu $ being respectively the canonical
tensors of energy-momentum and spin density. Because of the contraction of
the ''interior'' indices $\widehat{\alpha },\widehat{\beta }$ in (2.40) and
(2.41) we used the Hodge operator $*_G$ instead of $*_H$ (hereafter we
consider $*_G=*).$

Varying (by taking into account the distinguishing by N--connection) the action
$$
S=\int \left| G\right| ^{1/2}\delta ^{n+m}u\left( {\Cal L}_{\left( G\right)
}+{\Cal L}_{\left( m\right) }\right)
$$
on the $\Gamma $-variables (2.31), we obtain the gauge-gravitational field
equations:%
$$
d\left( *{\Cal R}^{(\Gamma )} \right) +
 \Gamma \bigwedge \left( *{\Cal R}^{(\Gamma )} \right) -\left( *%
{\Cal R}^{(\Gamma )}\right) \bigwedge \Gamma =-\lambda \left( *{\Cal I}\right).
\tag2.43
$$

Specifying the variations on $\Gamma _{\quad \widehat{\beta }}^{\widehat{%
\alpha }}$ and $l^{\widehat{\alpha }}$-variables, we rewrite (2.43) as
$$
\widehat{{\Cal D}}\left( *{\Cal R}^{(\Gamma )} \right) +
\frac{2\lambda }{l^2}\left(
\widehat{{\Cal D}}\left( *\pi \right) +\chi \bigwedge \left( *T^T\right)
-\left( *T\right) \bigwedge \chi ^T\right) =-\lambda \left( *S\right),
\tag2.44
$$
$$
\widehat{{\Cal D}}\left( *T\right) -\left( *{\Cal R}^{(\Gamma )}
\right) \bigwedge \chi -%
\frac{2\lambda }{l^2}\left( *\pi \right) \bigwedge \chi =\frac{l^2}2\left(
*t+\frac 1\lambda *\tau \right) ,\tag2.45
$$
where
$$
T^t=\{T_{\widehat{\alpha }}=\eta _{\widehat{\alpha }\widehat{\beta }}T^{%
\widehat{\beta }},~T^{\widehat{\beta }}=\frac 12T_{\quad \mu \nu }^{\widehat{%
\beta }}\delta u^\mu \bigwedge \delta u^\nu \},
$$
$$
\chi ^T=\{\chi _{\widehat{\alpha }}=\eta _{\widehat{\alpha }\widehat{\beta }%
}\chi ^{\widehat{\beta }},~\chi ^{\widehat{\beta }}=\chi _{\quad \mu }^{%
\widehat{\beta }}\delta u^\mu \},\qquad \widehat{{\Cal D}}=d+\widehat{\Gamma
}
$$
($\widehat{\Gamma }$ acts as $\Gamma _{\quad \widehat{\beta }\mu }^{\widehat{%
\alpha }}$ on indices $\widehat{\gamma },\widehat{\delta },...$ and as $%
\Gamma _{\quad \beta \mu }^\alpha $ on indices $\gamma ,\delta ,...).$ In
(2.45), $\tau $ defines the energy-momentum tensor of the $S_{\left( \eta
\right) }$-gauge gravitational field $\widehat{\Gamma }:$%
$$
\tau _{\mu \nu }\left( \widehat{\Gamma }\right) =\frac 12Tr\left( {\Cal R}%
_{\mu \alpha }{\Cal R}_{\quad \nu }^\alpha -\frac 14{\Cal R}_{\alpha \beta }%
{\Cal R}^{\alpha \beta }G_{\mu \nu }\right) .\tag2.46
$$

Equations (2.43) (or equivalently (2.44),(2.45)) make up the complete system
of variational field equations for nonlinear de Sitter gauge gravity with
local anisot\-ro\-py. They can be interpreted as a generalization of Miron and
Anastasiei equations for la--gravity \cite{Vacaru and Goncharenko 1995}
(equivalently, of gauge gravitational equations (2.33)] to a system of gauge
field equations with dynamical torsion and corresponding spin--density source.

Finally, we remark that we can obtain a nonvariational Poincare gauge
gravitational theory on la-spaces if we consider the contraction of the
gauge potential (2.38) to a potential with values in the Poincare Lie
algebra
$$
\Gamma =\left(
\matrix
\Gamma _{\quad \widehat{\beta }}^{\widehat{\alpha }} & l_0^{-1}\chi ^{
\widehat{\alpha }} \\ l_0^{-1}\chi _{\widehat{\beta }} & 0
\endmatrix
\right) \rightarrow \Gamma =\left(
\matrix
\Gamma _{\quad \widehat{\beta }}^{\widehat{\alpha }} & l_0^{-1}\chi ^{
\widehat{\alpha }} \\ 0 & 0
\endmatrix
\right) .
$$
Isotropic Poincare gauge gravitational theories are studied in a number of
papers (see, for example, references from
\cite{Tseytlin 1982} and \cite{Ponomarev, Barvinsky and Obukhov 1985}).
In a manner
similar to considerations presented in this work, we can generalize Poincare
gauge models for spaces with local anisotropy.
\vskip80pt
\newpage
\topmatter
\title
\chapter{3} NEARLY AUTOPARALLEL MAPS AND CONSERVATION LAWS \endtitle
\endtopmatter

The study of models of classical and quantum field interactions on la--spaces
is in order of the day. The development of this branch of theoretical and
mathematical physics entails great difficulties because of problematical
character of the possibility and manner of definition of conservation laws on
la--spaces. It will be recalled that, for instance, conservation laws of
energy--momentum type are a consequence of existence of a global group of
automorphisms of the fundamental Mikowski spaces
(the tangent space's automorphisms and particular cases when there are
symmetries generated by existence of Killing vectors are considered for
 (pseudo)Riemannian spaces). No
global or local automorphisms exist on generic la-spaces and in result of
this fact the formulation of la-conservation laws is sophisticate and full
of ambiguities. R. Miron and M. Anastasiei firstly pointed out the nonzero
divergence of the matter energy-momentum d--tensor, the source in Einstein
equations on la--spaces, and considered an original approach to the geometry
of time--dependent Lagrangians \cite{Miron and Anastasiei 1994}.
 Nevertheless, the
rigorous definition of energy-momentum values for la--gravitational and
matter fields and the form of conservation laws for such values have not
been considered in present--day studies of the mentioned problem.

The question of definition of tensor integration as the inverse operation of
covariant derivation was posed and studied by \cite{Mo\'or 1951}.
Tensor--integral and bitensor formalisms turned out to be very useful in
solving certain problems connected with conservation laws in general
relativity.  In order to extend tensor--integral constructions we have
proposed \cite{Gottlieb and Vacaru 1996} to take into consideration nearly
 autoparallel \cite{Vacaru 1992} and nearly geodesic  \cite{Sinyukov 1992}
maps, ng--maps, which forms a subclass of local 1--1 maps of  curved spaces
with deformation of the connection and metric structures.

 One of the main purposes of this Chapter is to synthesize the results on
nearly autoparallel maps and tensor integral  and to formulate them for a very
 general class of la--spaces. As the next step the la--gravity and an
 analysis of la--conservation laws are considered.

We note that proofs of our theorems are mechanical, but, in most cases, they
are rather tedious calculations similar to those presented in
\cite{Sinyukov 1979}, \cite{Vacaru and Ostaf 1996a, 1996b} and
\cite{Gottlieb and Vacaru 1996}. Some of them, on la--spaces, will be given
 in detail the rest, being similar, or consequences, will be only sketched or
 omitted.
\vskip25pt

\head{III.1 Nearly Autoparallel Maps of Locally Anisotropic Spaces}\endhead

In our  geometric constructions we shall use pairs of open regions
$( U, {\underline U})$ of la--spaces,
$U{\subset}{\xi},\, {\underline U}{\subset}{\underline {\xi}%
}$, and 1--1 local maps $f : U{\to}{\underline U}$ given by functions $%
f^{a}(u)$ of smoothly class $C^r( U) \, (r>2, $ or $r={\omega}$~ for
analytic functions) and their inverse functions $f^{\underline a}({%
\underline u})$ with corresponding non--zero Jacobians in every point $u{\in}
U$ and ${\underline u}{\in}{\underline U}.$

We consider that two open regions $U$~ and ${\underline U}$~ are attributed
to a common for f--map coordinate system if this map is realized on the
principle of coordinate equality $q(u^{\alpha}) {\to} {\underline q}
(u^{\alpha})$~ for every point $q {\in} U$~ and its f--image ${\underline q}
{\in} {\underline U}$ and note that all calculations
included in this Chapter will be local in nature and taken to refer to open
 subsets of mappings of
type $\frac{{\xi} {\supset} U \buildrel f }{\longrightarrow {\underline U} {%
\subset} {\underline {\xi}}.}$ For simplicity, we suppose that in a fixed
common coordinate system for $U$ and ${\underline U}$ spaces $\xi$ and ${%
\underline {\xi}}$ are characterized by a common N--connection structure (in
consequence of (1.10) by a corresponding concordance of d--metric structure),
i.e.
$$
N^{a}_{j}(u)={\underline N}^{a}_{j}(u)={\underline N}^{a}_{j} ({\underline u}%
),%
$$
which leads to the possibility to establish common local bases, adapted to a
given N--connection, on both regions $U$ and ${\underline U.}$ We consider
that on $\xi$ it is defined the linear d--connection structure with
components ${\Gamma}^{\alpha}_{{.}{\beta}{\gamma}}.$ On the space $%
\underline {\xi}$ the linear d--connection is considered to be a  general
one with torsion
$$
{\underline T}^{\alpha}_{{.}{\beta}{\gamma}}={\underline {\Gamma}}^{\alpha}_
{{.}{\beta}{\gamma}}-{\underline {\Gamma}}^{\alpha}_{{.}{\gamma}{\beta}}+
w^{\alpha}_{{.}{\beta}{\gamma}}
$$
and nonmetricity
$$
{\underline K}_{{\alpha}{\beta}{\gamma}}={{\underline D}_{\alpha}} {%
\underline G}_{{\beta}{\gamma}}. \tag3.1
$$

Geometrical objects on ${\underline {\xi}}$ are specified by underlined
symbols (for example, ${\underline A}^{\alpha}, {\underline B}^{{\alpha}{%
\beta}})$~ or underlined indices (for example, $A^{\underline a}, B^{{%
\underline a}{\underline b}}).$

For our purposes it is convenient to introduce auxiliary sym\-met\-ric
d--con\-nec\-ti\-ons, ${\gamma}^{\alpha}_{{.}{\beta}{\gamma}}={\gamma}%
^{\alpha}_{{.}{\gamma}{\beta}} $~ on $\xi$ and ${\underline {\gamma}}%
^{\alpha}_{.{\beta}{\gamma}}= {\underline {\gamma}}^{\alpha}_{{.}{\gamma}{%
\beta}}$ on ${\underline {\xi}}$ defined, correspondingly, as
$$
{\Gamma}^{\alpha}_{{.}{\beta}{\gamma}}= {\gamma}^{\alpha}_{{.}{\beta}{\gamma}%
}+  T^{\alpha}_{{.}{\beta}{\gamma}}\quad { and }\quad {\underline
{\Gamma}}^{\alpha}_{{.}{\beta}{\gamma}}= {\underline {\gamma}}^{\alpha}_{{.}{%
\beta}{\gamma}}+ {\underline T}^{\alpha}_{{.}{\beta}{\gamma}}.%
$$

We are interested in definition of local 1--1 maps from $U$ to ${\underline U%
}$ characterized by symmetric, $P^{\alpha}_{{.}{\beta}{\gamma}},$ and
antisymmetric, $Q^{\alpha}_{{.}{\beta}{\gamma}}$,~ deformations:
$$
{\underline {\gamma}}^{\alpha}_{{.}{\beta}{\gamma}} ={\gamma}^{\alpha}_{{.}{%
\beta}{\gamma}}+ P^{\alpha}_{{.}{\beta}{\gamma}} \tag3.2
$$
and
$$
{\underline T}^{\alpha}_{{.}{\beta}{\gamma}}= T^{\alpha}_{{.}{\beta}{\gamma}%
}+ Q^{\alpha}_{{.}{\beta}{\gamma}}. \tag3.3
$$
The auxiliary linear covariant derivations induced by ${\gamma}^{\alpha}_{{.}%
{\beta}{\gamma}}$ and ${\underline {\gamma}}^{\alpha}_{{.}{\beta}{\gamma}}$~
are denoted respectively as $^{({\gamma})}D$~ and $^{({\gamma})}{\underline D%
}.$~

Let introduce this local coordinate parametrization of curves on $U$~:
$$
u^{\alpha}=u^{\alpha}({\eta})=(x^{i}({\eta}), y^{i}({\eta}),~{\eta}_{1}<{\eta%
}<{\eta}_{2},%
$$
where corresponding tangent vector field is defined as
$$
v^{\alpha}={\frac{{du^{\alpha}} }{d{\eta}}}= ({\frac{{dx^{i}({\eta})} }{{d{%
\eta}}}}, {\frac{{dy^{j}({\eta})} }{d{\eta}}}).%
$$

\proclaim{\bf Definition 3.1}
A curve $l$~ is called auto parallel, a--parallel, on $\xi $ if its tangent
vector field $v^\alpha $~ satisfies a--parallel equations:
$$
vDv^\alpha =v^\beta {^{({\gamma })}D}_\beta v^\alpha ={\rho }({\eta }%
)v^\alpha ,\tag3.4
$$
where ${\rho }({\eta })$~ is a scalar function on $\xi $.
\endproclaim

Let curve ${\underline l} {\subset} {\underline {\xi}}$ is given in
parametric form as $u^{\alpha}=u^{\alpha}({\eta}),~{\eta}_1 < {\eta} <{\eta}%
_2$ with tangent vector field $v^{\alpha} = {\frac{{du^{\alpha}} }{{d{\eta}}}%
} {\ne} 0.$ We suppose that a 2--dimensional distribution $E_2({\underline l}%
)$ is defined along ${\underline l} ,$ i.e. in every point $u {\in} {%
\underline l}$ is fixed a 2-dimensional vector space $E_{2}({\underline l}) {%
\subset} {\underline {\xi}}.$ The introduced distribution $E_{2}({\underline
l})$~ is coplanar along ${\underline l}$~ if every vector ${\underline p}%
^{\alpha}(u^{b}_{(0)}) {\subset} E_{2}({\underline l}), u^{\beta}_{(0)} {%
\subset} {\underline l}$~ rests contained in the same distribution after
parallel transports along ${\underline l},$~ i.e. ${\underline p}%
^{\alpha}(u^{\beta}({\eta})) {\subset} E_{2} ({\underline l}).$

\proclaim{\bf Definition 3.2}
A curve ${\underline{l}}$~ is called nearly autoparallel, or in brief an
na--parallel, on space ${\underline{\xi }}$~ if a coplanar along ${%
\underline{l}}$~ distribution $E_2({\underline{l}})$ containing tangent to ${%
\underline{l}}$~ vector field $v^\alpha ({\eta })$,~ i.e. $v^\alpha ({\eta })%
{\subset }E_2({\underline{l}}),$~ is defined.
\endproclaim

The nearly autoparallel maps of la--spaces are introduced according the
 definition:

\proclaim{\bf Definition 3.3}
Nearly autoparallel maps, na--maps, of la--spaces are defined as local 1--1
mappings of v--bundles, $\xi {\to }{\underline{\xi }},$ changing every
a--parallel on $\xi $ into a na--parallel on ${\underline{\xi }}.$
\endproclaim

Now we formulate the general conditions when deformations (3.2) and (3.3)
charac\-ter\-ize na-maps : Let a-parallel $l{\subset }U$~ is given by
func\-ti\-ons
 $u^\alpha =u^{({\alpha })}({\eta }),v^\alpha ={\frac{{%
du^\alpha }}{d{\eta }}}$, ${\eta }_1<{\eta }<{\eta }_2$, satisfying
equations (3.4). We suppose that to this a--parallel corresponds a
na--parallel ${\underline{l}}\subset {\underline{U}}$ given by the same
parameterization in a common for a chosen na--map coordinate system on $U$~
and ${\underline{U}}.$ This condition holds for vectors ${\underline{v}}%
_{(1)}^\alpha =v{\underline{D}}v^\alpha $~ and $v_{(2)}^\alpha =v{\underline{%
D}}v_{(1)}^\alpha $ satisfying equality
$$
{\underline{v}}_{(2)}^\alpha ={\underline{a}}({\eta })v^\alpha +{\underline{b%
}}({\eta }){\underline{v}}_{(1)}^\alpha \tag3.5
$$
for some scalar functions ${\underline{a}}({\eta })$~ and ${\underline{b}}({%
\eta })$~ (see Definitions 3.2 and 3.3). Putting splittings (3.2) and (3.3)
into expressions for ${\underline{v}}_{(1)}^\alpha$  and
${\underline{v}}_{(2)}^\alpha $ in (3.5) we obtain:
$$
v^\beta v^\gamma v^\delta (D_\beta P_{{.}{\gamma }{\delta }}^\alpha +P_{{.}{%
\beta }{\tau }}^\alpha P_{{.}{\gamma }{\delta }}^\tau +Q_{{.}{\beta }{\tau }%
}^\alpha P_{{.}{\gamma }{\delta }}^\tau )=bv^\gamma v^\delta P_{{.}{\gamma }{%
\delta }}^\alpha +av^\alpha ,\tag3.6
$$
where
$$
b({\eta },v)={\underline{b}}-3{\rho },\qquad \text{and}\qquad a({\eta },v)={%
\underline{a}}+{\underline{b}}{\rho }-v^b{\partial }_b{\rho }-{\rho }^2%
\tag3.7
$$
are called the deformation parameters of na--maps.

The algebraic equations for the deformation of torsion $Q_{{.}{\beta }{\tau }%
}^\alpha $ should be written as the compatibility conditions for a given
nonmetricity tensor ${\underline{K}}_{{\alpha }{\beta }{\gamma }}$~ on ${%
\underline{\xi }}$ ( or as the metricity conditions if d--connection ${%
\underline{D}}_\alpha $~ on ${\underline{\xi }}$~ is required to be metric) :%
$$
D_\alpha G_{{\beta }{\gamma }}-P_{{.}{\alpha }({\beta }}^\delta G_{{{\gamma }%
)}{\delta }}-{\underline{K}}_{{\alpha }{\beta }{\gamma }}=Q_{{.}{\alpha }({%
\beta }}^\delta G_{{\gamma }){\delta }},\tag3.8
$$
where $({\quad })$ denotes the symmetric alternation.

So, we have proved this

\proclaim{\bf Theorem 3.1}
The na--maps from la--space $\xi $ to la--space ${\underline{\xi }}$~ with a
fixed common nonlinear connection $N_j^a(u)={\underline{N}}_j^a(u)$ and
given d--connections, ${\Gamma }_{{.}{\beta }{\gamma }}^\alpha $~ on $\xi $~
and ${\underline{\Gamma }}_{{.}{\beta }{\gamma }}^\alpha $~ on ${\underline{%
\xi }}$ are locally parametrized by the solutions of equations (3.6) and
(3.8) for every point $u^\alpha $~ and direction $v^\alpha $~ on $U{\subset }%
{\xi }.$
\endproclaim

We call (3.6) and (3.8) the basic equations for na--maps of la--spaces. They
generalize the corresponding  equations \cite{Sinyukov 1979} for isotropic
spaces provided with symmetric affine connection structure.
\vskip25pt

\head{III.2\ Classification of Nearly Autoparallel Maps of LA--Spaces}\endhead

Na--maps are classed on possible polynomial parametrizations on variables $%
v^{\alpha}$~ of deformations parameters $a$ and $b$ (see (3.6) and (3.7) ).

\proclaim{\bf Theorem 3.2}
There are four classes of na--maps characterized by corresponding
deformation parameters and tensors and basic equations:
\item 0.  for $na_{(0)}$--maps, ${\pi }_{(0)}$--maps,
$$
P_{{\beta }{\gamma }}^\alpha (u)={\psi }_{{(}{\beta }}{\delta }_{{\gamma }%
)}^\alpha
$$
(${\delta }_\beta ^\alpha $~ is Kronecker symbol and ${\psi }_\beta ={\psi }%
_\beta (u)$~ is a covariant vector d--field) ;

{\ \item 1.}  for $na_{(1)}$--maps
$$
a(u,v)=a_{{\alpha }{\beta }}(u)v^\alpha v^\beta ,\quad b(u,v)=b_\alpha
(u)v^\alpha
$$
and $P_{{.}{\beta }{\gamma }}^\alpha (u)$~ is the solution of equations
$$
D_{({\alpha }}P_{{.}{\beta }{\gamma })}^\delta +P_{({\alpha }{\beta }}^\tau
P_{{.}{\gamma }){\tau }}^\delta -P_{({\alpha }{\beta }}^\tau Q_{{.}{\gamma })%
{\tau }}^\delta =b_{({\alpha }}P_{{.}{\beta }{\gamma })}^{{\delta }}+a_{({%
\alpha }{\beta }}{\delta }_{{\gamma })}^\delta ;\tag3.9
$$
{\ \item 2.}  for $na_{(2)}$--maps
$$
a(u,v)=a_\beta (u)v^\beta ,\quad b(u,v)={\frac{{b_{{\alpha }{\beta }%
}v^\alpha v^\beta }}{{{\sigma }_\alpha (u)v^\alpha }}},\quad {\sigma }%
_\alpha v^\alpha {\neq }0,
$$
$$
P_{{.}{\alpha }{\beta }}^\tau (u)={{\psi }_{({\alpha }}}{\delta }_{{\beta }%
)}^\tau +{\sigma }_{({\alpha }}F_{{\beta })}^\tau
$$
and $F_\beta ^\alpha (u)$~ is the solution of equations
$$
{D}_{({\gamma }}F_{{\beta })}^\alpha +F_\delta ^\alpha F_{({\gamma }}^\delta
{\sigma }_{{\beta })}-Q_{{.}{\tau }({\beta }}^\alpha F_{{\gamma })}^\tau ={%
\mu }_{({\beta }}F_{{\gamma })}^\alpha +{\nu }_{({\beta }}{\delta }_{{\gamma
})}^\alpha \tag3.10
$$
$({\mu }_\beta (u),{\nu }_\beta (u),{\psi }_\alpha (u),{\sigma }_\alpha (u)$%
~ are covariant d--vectors) ;
{\ \item 3.}  for $na_{(3)}$--maps
$$
b(u,v)={\frac{{{\alpha }_{{\beta }{\gamma }{\delta }}v^\beta v^\gamma
v^\delta }}{{{\sigma }_{{\alpha }{\beta }}v^\alpha v^\gamma }}},
$$
$$
P_{{.}{\beta }{\gamma }}^\alpha (u)={\psi }_{({\beta }}{\delta }_{{\gamma }%
)}^\alpha +{\sigma }_{{\beta }{\gamma }}{\varphi }^\alpha ,
$$
where ${\varphi }^\alpha $~ is the solution of equations
$$
D_\beta {\varphi }^\alpha ={\nu }{\delta }_\beta ^\alpha +{\mu }_\beta {%
\varphi }^\alpha +{\varphi }^\gamma Q_{{.}{\gamma }{\delta }}^\alpha ,%
\tag3.11
$$
${\alpha }_{{\beta }{\gamma }{\delta }}(u),{\sigma }_{{\alpha }{\beta }}(u),{%
\psi }_\beta (u),{\nu }(u)$~ and ${\mu }_\beta (u)$~ are d--tensors.
\endproclaim
\vskip10pt

{\bf Proof.} We sketch the proof respectively for every point in the theorem:

{\ \item 0.}
  It is easy to verify that a--parallel equations (3.4) on $\xi $
transform into similar ones on $\underline{\xi }$ if and only if
deformations (3.2) with deformation d--tensors of type ${P^\alpha }_{\beta
\gamma }(u)={\psi }_{(\beta }{\delta }_{\gamma )}^\alpha $ are considered.

{\  \item 1.}
  Using corresponding to $na_{(1)}$--maps parametrizations of $a(u,v)$
and $b(u,v)$ (see conditions of the theorem) for arbitrary $v^\alpha \neq 0$
on $U\in \xi $ and after a redefinition of deformation parameters we obtain
that equations (3.6) hold if and only if ${P^\alpha }_{\beta \gamma }$
satisfies (3.3).

{\ \item 2.}  In a similar manner we obtain basic $na_{(2)}$--map equations (3.10)
from (3.6) by considering $na_{(2)}$--parametrizations of deformation
parameters and d--tensor.

{\ \item 3.}  For $na_{(3)}$--maps we mast take into consideration deformations of
torsion (3.3) and introduce $na_{(3)}$--parametrizations for $b(u,v)$ and ${%
P^\alpha }_{\beta \gamma }$ into the basic na--equa\-ti\-ons (3.6). The last
ones for $na_{(3)}$--maps are equivalent to equations (3.11) (with a
corresponding redefinition of deformation parameters). \qquad $\blacksquare$
\vskip10pt

We point out that for ${\pi}_{(0)}$-maps we do not have differential equations
on $P^{\alpha}_{{.}{\beta}{\gamma}}$ (in the  isotropic case one considers a
first order system of differential equations on metric \cite{Sinyukov 1979};
 we omit
constructions  with deformation of metric in this work).\

To formulate invariant conditions for reciprocal na--maps (when every
a-parallel on ${\underline {\xi}}$~ is also transformed into na--parallel on
$\xi$ ) it is convenient to introduce into consideration the curvature and
Ricci tensors defined for auxiliary connection ${\gamma}^{\alpha}_{{.}{\beta}%
{\gamma}}$~ :
$$
r^{{.}{\delta}} _{{\alpha}{.}{\beta}{\tau}}={\partial}_{[{\beta}}{\gamma}%
^{\delta}_ {{.}{\tau}]{\alpha}}+{\gamma}^{\delta}_{{.}{\rho}[{\beta}}{\gamma}%
^{\rho}_ {{.}{\tau}]{\alpha}} + {{\gamma}^{\delta}}_{\alpha \phi} {w^{\phi}}%
_{\beta \tau}%
$$
and, respectively, $r_{{\alpha}{\tau}}=r^{{.}{\gamma}} _{{\alpha}{.}{\gamma}{%
\tau}} $, where $[\quad ]$ denotes antisymmetric alternation of indices, and
to define values:
$$
^{(0)}T^{\mu}_{{.}{\alpha}{\beta}}= {\Gamma}^{\mu}_{{.}{\alpha}{\beta}} -
T^{\mu}_{{.}{\alpha}{\beta}}- {\frac{1 }{(n+m + 1)}}({\delta}^{\mu}_{({\alpha%
}}{\Gamma}^{\delta}_ {{.}{\beta}){\delta}}-{\delta}^{\mu}_{({\alpha}%
}T^{\delta}_ {{.}{\beta}){\gamma}}),
$$
$$
{}^{(0)}{W}^{\cdot \tau}_{\alpha \cdot \beta \gamma} = {r}^{\cdot
\tau}_{\alpha \cdot \beta \gamma} + {\frac{1}{n+m+1}} [ {\gamma}%
^{\tau}_{\cdot \varphi \tau} {\delta}^{\tau}_{( \alpha} {w^{\varphi}}_{\beta
) \gamma} - ( {\delta}^{\tau}_{\alpha}{r}_{[ \gamma \beta ]} + {\delta}%
^{\tau}_{\gamma} {r}_{[ \alpha \beta ]} - {\delta}^{\tau}_{\beta} {r}_{[
\alpha \gamma ]} )] -%
$$
$$
{\frac{1}{{(n+m+1)}^2}} [ {\delta}^{\tau}_{\alpha}  ( 2 {\gamma}%
^{\tau}_{\cdot \varphi \tau} {w^{\varphi}}_{[ \gamma \beta ] } - {\gamma}%
^{\tau}_{\cdot \tau [ \gamma } {w^{\varphi}}_{\beta ] \varphi} ) + {\delta}%
^{\tau}_{\gamma}  ( 2 {\gamma}^{\tau}_{\cdot \varphi \tau} {w^{\varphi}}%
_{\alpha \beta} -{\gamma}^{\tau}_{\cdot \alpha \tau} {w^{\varphi}}_{\beta
\varphi}) -
$$
$$
{\delta}^{\tau}_{\beta}  ( 2 {\gamma}^{\tau}_{\cdot \varphi \tau} {%
w^{\varphi}}_{\alpha \gamma} - {\gamma}^{\tau}_{\cdot \alpha \tau} {%
w^{\varphi}}_{\gamma \varphi} ) ],%
$$

$$
{^{(3)}T}^{\delta}_{{.}{\alpha}{\beta}}= {\gamma}^{\delta}_{{.}{\alpha}{\beta%
}}+ {\epsilon}{\varphi}^{\tau}{^{({\gamma})}D}_{\beta}q_{\tau}+ {\frac{1 }{%
n+m}}({\delta}^{\gamma}_{\alpha}- {\epsilon}{\varphi}^{\delta}q_{\alpha})[{%
\gamma}^{\tau}_{{.}{\beta}{\tau}}+ {\epsilon}{\varphi}^{\tau}{^{({\gamma})}D}%
_{\beta}q_{\tau}+%
$$
$$
{\frac{1 }{{n+m -1}}}q_{\beta}({\epsilon}{\varphi}^{\tau}{\gamma}^{\lambda}_
{{.}{\tau}{\lambda}}+ {\varphi}^{\lambda}{\varphi}^{\tau}{^{({\gamma})}D}%
_{\tau}q_{\lambda})]- {\frac{1 }{n+m}}({\delta}^{\delta}_{\beta}-{\epsilon}{%
\varphi}^{\delta} q_{\beta})[{\gamma}^{\tau}_{{.}{\alpha}{\tau}}+ {\epsilon}{%
\varphi}^{\tau} {^{({\gamma})}D}_{\alpha}q_{\tau}+%
$$
$$
{\frac{1 }{{n+m -1}}}q_{\alpha}({\epsilon}{\varphi}^{\tau}{\gamma}%
^{\lambda}_ {{.}{\tau}{\lambda}}+ {\varphi}^{\lambda}{\varphi}^{\tau} {%
^{(\gamma)}D} _{\tau}q_{\lambda})],%
$$

$$
{^{(3)}W}^\alpha {{.}{\beta }{\gamma }{\delta }}={\rho }_{{\beta }{.}{\gamma
}{\delta }}^{{.}{\alpha }}+{\epsilon }{\varphi }^\alpha q_\tau {\rho }_{{%
\beta }{.}{\gamma }{\delta }}^{{.}{\tau }}+({\delta }_\delta ^\alpha -
$$
$$
{\epsilon }{\varphi }^\alpha q_\delta )p_{{\beta }{\gamma }}-({\delta }%
_\gamma ^\alpha -{\epsilon }{\varphi }^\alpha q_\gamma )p_{{\beta }{\delta }%
}-({\delta }_\beta ^\alpha -{\epsilon }{\varphi }^\alpha q_\beta )p_{[{%
\gamma }{\delta }]},
$$
$$
(n+m-2)p_{{\alpha }{\beta }}=-{\rho }_{{\alpha }{\beta }}-{\epsilon }q_\tau {%
\varphi }^\gamma {\rho }_{{\alpha }{.}{\beta }{\gamma }}^{{.}{\tau }}+{\frac
1{n+m}}[{\rho }_{{\tau }{.}{\beta }{\alpha }}^{{.}{\tau }}-{\epsilon }q_\tau
{\varphi }^\gamma {\rho }_{{\gamma }{.}{\beta }{\alpha }}^{{.}{\tau }}+{%
\epsilon }q_\beta {\varphi }^\tau {\rho }_{{\alpha }{\tau }}+
$$
$$
{\epsilon }q_\alpha (-{\varphi }^\gamma {\rho }_{{\tau }{.}{\beta }{\gamma }%
}^{{.}{\tau }}+{\epsilon }q_\tau {\varphi }^\gamma {\varphi }^\delta {\rho }%
_{{\gamma }{.}{\beta }{\delta }}^{{.}{\tau }}]),
$$
where $q_\alpha {\varphi }^\alpha ={\epsilon }=\pm 1,$
$$
{{\rho }^\alpha }_{\beta \gamma \delta }=r_{\beta \cdot \gamma \delta
}^{\cdot \alpha }+{\frac 12}({\psi }_{(\beta }{\delta }_{\varphi )}^\alpha +{%
\sigma }_{\beta \varphi }{\varphi }^\tau ){w^\varphi }_{\gamma \delta }
$$
( for a similar value on $\underline{\xi }$ we write ${\quad }{\underline{%
\rho }}_{\cdot \beta \gamma \delta }^\alpha ={\underline{r}}_{\beta \cdot
\gamma \delta }^{\cdot \alpha }-{\frac 12}({\psi }_{(\beta }{\delta }_{{%
\varphi })}^\alpha -{\sigma }_{\beta \varphi }{\varphi }^\tau ){w^\varphi }%
_{\gamma \delta }{\quad })$ and ${\rho }_{\alpha \beta }={\rho }_{\cdot
\alpha \beta \tau }^\tau .$

Similar values, $
^{(0)}{\underline T}^{\alpha}_{{.}{\beta}{\gamma}}, ^{(0)}{\underline W}%
^{\nu}_{{.}{\alpha}{\beta}{\gamma}}, {\hat T}^{\alpha} _{{.}{\beta}{\gamma}%
}, {\check T}^{\alpha}_ {{.}{\beta}{\tau}}, {\hat W}^{\delta}_{{.}{\alpha}{%
\beta}{\gamma}}, {\check W}^{\delta}_ {{.}{\alpha}{\beta}{\gamma}}, ^{(3)}{%
\underline T}^{\delta} _{{.}{\alpha}{\beta}},$
and $^{(3)}{\underline W}^ {\alpha}_{{.}{\beta}{\gamma}{\delta}} $ are
given, correspondingly, by auxiliary connections ${\quad}{\underline {\Gamma}%
}^{\mu}_{{.}{\alpha}{\beta}},$~
$$
{\star {\gamma}}^{\alpha}_{{.}{\beta}{\lambda}}={\gamma}^{\alpha} _{{.}{\beta%
}{\lambda}} + {\epsilon}F^{\alpha}_{\tau}{^{({\gamma})}D}_{({\beta}}
F^{\tau}_{{\lambda})}, \quad {\check {\gamma}}^{\alpha}_{{.}{\beta}{\lambda}%
}= {\widetilde {\gamma}}^{\alpha}_{{.}{\beta}{\lambda}} + {\epsilon}%
F^{\lambda} _{\tau} {\widetilde D}_{({\beta}}F^{\tau}_{{\lambda})},%
$$
$$
{\widetilde {\gamma}}^{\alpha}_{{.}{\beta}{\tau}}={\gamma}^{\alpha} _{{.}{%
\beta}{\tau}}+ {\sigma}_{({\beta}}F^{\alpha}_{{\tau})}, \quad {\hat {\gamma}}%
^{\alpha}_{{.}{\beta}{\lambda}}={\star {\gamma}}^{\alpha}_ {{.}{\beta}{%
\lambda}} + {\widetilde {\sigma}}_{({\beta}}{\delta}^{\alpha}_ {{\lambda})},%
$$
where ${\widetilde {\sigma}}_{\beta}={\sigma}_{\alpha}F^{\alpha}_{\beta}.$

\proclaim{\bf Theorem 3.3}
Four classes of reciprocal na--maps of la--spaces are characterized by
corresponding invariant criterions:

{\ \item 0.}
  for a--maps $^{(0)}T_{{.}{\alpha }{\beta }}^\mu =^{(0)}{\underline{T}}%
_{{.}{\alpha }{\beta }}^\mu ,$
$$
{}^{(0)}W_{{.}{\alpha }{\beta }{\gamma }}^\delta =^{(0)}{\underline{W}}_{{.}{%
\alpha }{\beta }{\gamma }}^\delta ;\tag3.12
$$

{\ \item 1.}  for $na_{(1)}$--maps
$$
3({^{({\gamma })}D}_\lambda P_{{.}{\alpha }{\beta }}^\delta +P_{{.}{\tau }{%
\lambda }}^\delta P_{{.}{\alpha }{\beta }}^\tau )=r_{({\alpha }{.}{\beta }){%
\lambda }}^{{.}{\delta }}-{\underline{r}}_{({\alpha }{.}{\beta }){\lambda }%
}^{{.}{\delta }}+\tag3.13
$$
$$
[T_{{.}{\tau }({\alpha }}^\delta P_{{.}{\beta }{\lambda })}^\tau +Q_{{.}{%
\tau }({\alpha }}^\delta P_{{.}{\beta }{\lambda })}^\tau +b_{({\alpha }}P_{{.%
}{\beta }{\lambda })}^\delta +{\delta }_{({\alpha }}^\delta a_{{\beta }{%
\lambda })}];
$$

{\ \item 2.}
  for $na_{(2)}$--maps ${\hat T}_{{.}{\beta }{\tau }}^\alpha ={\star T}%
_{{.}{\beta }{\tau }}^\alpha ,$
$$
{\hat W}_{{.}{\alpha }{\beta }{\gamma }}^\delta ={\star W}_{{.}{\alpha }{%
\beta }{\gamma }}^\delta ;\tag3.14
$$

{\ \item 3.}
  for $na_{(3)}$--maps $^{(3)}T_{{.}{\beta }{\gamma }}^\alpha =^{(3)}{%
\underline{T}}_{{.}{\beta }{\gamma }}^\alpha ,$
$$
{}^{(3)}W_{{.}{\beta }{\gamma }{\delta }}^\alpha =^{(3)}{\underline{W}}_{{.}{%
\beta }{\gamma }{\delta }}^\alpha .\tag3.15
$$
\endproclaim
\vskip10pt

{\bf Proof. }

{\ \item 0.}
  Let us prove that a--invariant conditions (3.12) hold. Deformations
of d--connections of type
$$
{}^{(0)}{\underline{\gamma }}_{\cdot \alpha \beta }^\mu ={{\gamma }^\mu }%
_{\alpha \beta }+{\psi }_{(\alpha }{\delta }_{\beta )}^\mu \tag3.16
$$
define a--applications. Contracting indices $\mu $ and $\beta $ we can write
$$
{\psi }_\alpha ={\frac 1{m+n+1}}({{\underline{\gamma }}^\beta }_{\alpha
\beta }-{{\gamma }^\beta }_{\alpha \beta }).\tag3.17
$$
Introducing d--vector ${\psi }_\alpha $ into previous relation and
expressing
$$
{{\gamma }^\alpha }_{\beta \tau }=-{T^\alpha }_{\beta \tau }+{{\Gamma }%
^\alpha }_{\beta \tau }
$$
and similarly for underlined values we obtain the first invariant conditions
from (3.12).

Putting deformation (3.16) into the formula for
$$
{\underline{r}}_{\alpha \cdot \beta \gamma }^{\cdot \tau }\quad \text{and}%
\quad {\underline{r}}_{\alpha \beta }={\underline{r}}_{\alpha \tau \beta
\tau }^{\cdot \tau }
$$
we obtain respectively relations
$$
{\underline{r}}_{\alpha \cdot \beta \gamma }^{\cdot \tau }-r_{\alpha \cdot
\beta \gamma }^{\cdot \tau }={\delta }_\alpha ^\tau {\psi }_{[\gamma \beta
]}+{\psi }_{\alpha [\beta }{\delta }_{\gamma ]}^\tau +{\delta }_{(\alpha
}^\tau {\psi }_{\varphi )}{w^\varphi }_{\beta \gamma }\tag3.18
$$
and
$$
{\underline{r}}_{\alpha \beta }-r_{\alpha \beta }={\psi }_{[\alpha \beta
]}+(n+m-1){\psi }_{\alpha \beta }+{\psi }_\varphi {w^\varphi }_{\beta \alpha
}+{\psi }_\alpha {w^\varphi }_{\beta \varphi },\tag3.19
$$
where
$$
{\psi }_{\alpha \beta }={}^{({\gamma })}D_\beta {\psi }_\alpha -{\psi }%
_\alpha {\psi }_\beta .
$$
Putting (3.16) into (3.19) we can express ${\psi }_{[\alpha \beta ]}$ as
$$\multline
{\psi }_{[\alpha \beta ]}={\frac 1{n+m+1}}[{\underline{r}}_{[\alpha \beta ]}+%
{\frac 2{n+m+1}}{\underline{\gamma }}_{\cdot \varphi \tau }^\tau {w^\varphi }%
_{[\alpha \beta ]}-{\frac 1{n+m+1}}{\underline{\gamma }}_{\cdot \tau [\alpha
}^\tau {w^\varphi }_{\beta ]\varphi }]- \\
{\frac 1{n+m+1}}[r_{[\alpha \beta ]}+{\frac 2{n+m+1}}{{\gamma }^\tau }%
_{\varphi \tau }{w^\varphi }_{[\alpha \beta ]}-{\frac 1{n+m+1}}{{\gamma }%
^\tau }_{\tau [\alpha }{w^\varphi }_{\beta ]\varphi }].
 \endmultline \tag3.20$$
To simplify our consideration we can choose an a--transform, parametrized by
corresponding $\psi $--vector from (3.16), (or fix a local coordinate cart)
the antisymmetrized relations (3.20) to be satisfied by d--tensor
$$\multline
{\psi }_{\alpha \beta }={\frac 1{n+m+1}}[{\underline{r}}_{\alpha \beta }+{%
\frac 2{n+m+1}}{\underline{\gamma }}_{\cdot \varphi \tau }^\tau {w^\varphi }%
_{\alpha \beta }-{\frac 1{n+m+1}}{\underline{\gamma }}_{\cdot \alpha \tau
}^\tau {w^\varphi }_{\beta \varphi }-r_{\alpha \beta }-\\
{\frac 2{n+m+1}}{{\gamma }^\tau }_{\varphi \tau }{w^\varphi }_{\alpha \beta
}+{\frac 1{n+m+1}}{{\gamma }^\tau }_{\alpha \tau }{w^\varphi }_{\beta
\varphi }] \endmultline \tag3.21$$
Introducing expressions (3.16),(3.20) and (3.21) into deformation of
curvature (3.17) we obtain the second conditions (3.12) of a--map invariance:
$$
^{(0)}W_{\alpha \cdot \beta \gamma }^{\cdot \delta }={}^{(0)}{\underline{W}}%
_{\alpha \cdot \beta \gamma }^{\cdot \delta },
$$
where the Weyl d--tensor on $\underline{\xi }$ (the extension of the usual
one for geodesic maps on (pseudo)--Riemannian spaces to the case of
v--bundles provided with N--connection structure) is defined as
$$
{}^{(0)}{\underline{W}}_{\alpha \cdot \beta \gamma }^{\cdot \tau }={%
\underline{r}}_{\alpha \cdot \beta \gamma }^{\cdot \tau }+{\frac 1{n+m+1}}[{%
\underline{\gamma }}_{\cdot \varphi \tau }^\tau {\delta }_{(\alpha }^\tau {%
w^\varphi }_{\beta )\gamma }-({\delta }_\alpha ^\tau {\underline{r}}%
_{[\gamma \beta ]}+{\delta }_\gamma ^\tau {\underline{r}}_{[\alpha \beta ]}-{%
\delta }_\beta ^\tau {\underline{r}}_{[\alpha \gamma ]})]-
$$
$$
{\frac 1{{(n+m+1)}^2}}[{\delta }_\alpha ^\tau (2{\underline{\gamma }}_{\cdot
\varphi \tau }^\tau {w^\varphi }_{[\gamma \beta ]}-{\underline{\gamma }}%
_{\cdot \tau [\gamma }^\tau {w^\varphi }_{\beta ]\varphi })+{\delta }_\gamma
^\tau (2{\underline{\gamma }}_{\cdot \varphi \tau }^\tau {w^\varphi }%
_{\alpha \beta }-{\underline{\gamma }}_{\cdot \alpha \tau }^\tau {w^\varphi }%
_{\beta \varphi })-
$$
$$
{\delta }_\beta ^\tau (2{\underline{\gamma }}_{\cdot \varphi \tau }^\tau {%
w^\varphi }_{\alpha \gamma }-{\underline{\gamma }}_{\cdot \alpha \tau }^\tau
{w^\varphi }_{\gamma \varphi })]
$$
The formula for $^{(0)}W_{\alpha \cdot \beta \gamma }^{\cdot \tau }$ written
similarly with respect to non--underlined values is presented in subsection
1.1.2 .

{\ \item 1.}
  To obtain $na_{(1)}$--invariant conditions we rewrite $na_{(1)}$%
--equations (3.9) as to consider in explicit form covariant derivation $^{({%
\gamma })}D$ and deformations (3.2) and (3.3):
$$\multline
2({}^{({\gamma })}D_\alpha {P^\delta }_{\beta \gamma }+{}^{({\gamma }%
)}D_\beta {P^\delta }_{\alpha \gamma }+{}^{({\gamma })}D_\gamma {P^\delta }%
_{\alpha \beta }+{P^\delta }_{\tau \alpha }{P^\tau }_{\beta \gamma }+
{P^\delta }_{\tau \beta }{P^\tau }_{\alpha \gamma }+{P^\delta }_{\tau \gamma
}{P^\tau }_{\alpha \beta }) \\ =
{T^\delta }_{\tau (\alpha }{P^\tau }_{\beta\gamma )}+
{H^\delta }_{\tau (\alpha }{P^\tau }_{\beta \gamma )}+b_{(\alpha }{P^\delta }%
_{\beta \gamma )}+a_{(\alpha \beta }{\delta }_{\gamma )}^\delta .
\endmultline \tag3.22$$
Alternating the first two indices in (3.22) we have
$$
2({\underline{r}}_{(\alpha \cdot \beta )\gamma }^{\cdot \delta }-r_{(\alpha
\cdot \beta )\gamma }^{\cdot \delta })=2({}^{(\gamma )}D_\alpha {P^\delta }%
_{\beta \gamma }+
$$
$$
{}^{(\gamma )}D_\beta {P^\delta }_{\alpha \gamma }-2{}^{(\gamma )}D_\gamma {%
P^\delta }_{\alpha \beta }+{P^\delta }_{\tau \alpha }{P^\tau }_{\beta \gamma
}+{P^\delta }_{\tau \beta }{P^\tau }_{\alpha \gamma }-2{P^\delta }_{\tau
\gamma }{P^\tau }_{\alpha \beta }).
$$
Substituting the last expression from (3.22) and rescalling the deformation
parameters and d--tensors we obtain the conditions (3.9).

{\bf \item 2.}
  Now we prove the invariant conditions for $na_{(2)}$--maps satisfying
conditions
$$
\epsilon \neq 0\quad \text{and}\quad \epsilon -F_\beta ^\alpha F_\alpha
^\beta \neq 0
$$
Let define the auxiliary d--connection
$$
{\tilde \gamma }_{\cdot \beta \tau }^\alpha ={\underline{\gamma }}_{\cdot
\beta \tau }^\alpha -{\psi }_{(\beta }{\delta }_{\tau )}^\alpha ={{\gamma }%
^\alpha }_{\beta \tau }+{\sigma }_{(\beta }F_{\tau )}^\alpha \tag3.23
$$
and write
$$
{\tilde D}_\gamma ={}^{({\gamma })}D_\gamma F_\beta ^\alpha +{\tilde \sigma }%
_\gamma F_\beta ^\alpha -{\epsilon }{\sigma }_\beta {\delta }_\gamma ^\alpha
,
$$
where ${\tilde \sigma }_\beta ={\sigma }_\alpha F_\beta ^\alpha ,$ or, as a
consequence from the last equality,
$$
{\sigma }_{(\alpha }F_{\beta )}^\tau ={\epsilon }F_\lambda ^\tau ({}^{({%
\gamma })}D_{(\alpha }F_{\beta )}^\alpha -{\tilde D}_{(\alpha }F_{\beta
)}^\lambda )+{\tilde \sigma }_{(}{\alpha }{\delta }_{\beta )}^\tau .
$$
Introducing auxiliary connections
$$
{\star {\gamma }}_{\cdot \beta \lambda }^\alpha ={\gamma }_{\cdot \beta
\lambda }^\alpha +{\epsilon }F_\tau ^\alpha {}^{({\gamma })}D_{(\beta
}F_{\lambda )}^\tau
$$
and
$$
{\check \gamma }_{\cdot \beta \lambda }^\alpha ={\tilde \gamma }_{\cdot
\beta \lambda }^\alpha +{\epsilon }F_\tau ^\alpha {\tilde D}_{(\beta
}F_{\lambda )}^\tau
$$
we can express deformation (3.23) in a form characteristic for a--maps:
$$
{\hat \gamma }_{\cdot \beta \gamma }^\alpha ={\star {\gamma }}_{\cdot \beta
\gamma }^\alpha +{\tilde \sigma }_{(\beta }{\delta }_{\lambda )}^\alpha .%
\tag3.24
$$
Now it's obvious that $na_{(2)}$--invariant conditions (3.24) are equivalent
with a--invariant conditions (3.12) written for d--connection (3.24). As a
matter of principle we can write formulas for such $na_{(2)}$--invariants in
terms of ''underlined'' and ''non--underlined'' values by expressing
consequently all used auxiliary connections as deformations of ''prime''
connections on $\xi $ and ''final'' connections on $\underline{\xi }.$ We
omit such tedious calculations in this work.

{\ \item 3.}
  Finally, we prove the last statement, for $na_{(3)}$--maps, of the
theorem. Let
$$
q_\alpha {\varphi }^\alpha =e=\pm 1,\tag3.25
$$
where ${\varphi }^\alpha $ is contained in
$$
{\underline{\gamma }}_{\cdot \beta \gamma }^\alpha ={{\gamma }^\alpha }%
_{\beta \gamma }+{\psi }_{(\beta }{\delta }_{\gamma )}^\alpha +{\sigma }%
_{\beta \gamma }{\varphi }^\alpha .\tag3.26
$$
Acting with operator $^{({\gamma })}{\underline{D}}_\beta $ on (3.25) we
write
$$
{}^{({\gamma })}{\underline{D}}_\beta q_\alpha ={}^{({\gamma })}D_\beta
q_\alpha -{\psi }_{(\alpha }q_{\beta )}-e{\sigma }_{\alpha \beta }.%
\tag3.27
$$
Contracting (3.27) with ${\varphi }^\alpha $ we can express
$$
e{\varphi }^\alpha {\sigma }_{\alpha \beta }={\varphi }^\alpha ({}^{({\gamma
})}D_\beta q_\alpha -{}^{({\gamma })}{\underline{D}}_\beta q_\alpha )-{%
\varphi }_\alpha q^\alpha q_\beta -e{\psi }_\beta .
$$
Putting the last formula in (3.26) contracted on indices $\alpha $ and $%
\gamma $ we obtain
$$
(n+m){\psi }_\beta ={\underline{\gamma }}_{\cdot \alpha \beta }^\alpha -{{%
\gamma }^\alpha }_{\alpha \beta }+e{\psi }_\alpha {\varphi }^\alpha q_\beta
+e{\varphi }^\alpha {\varphi }^\beta ({}^{({\gamma })}{\underline{D}}_\beta
-{}^{({\gamma })}D_\beta ).\tag3.28
$$
From these relations, taking into consideration (3.25), we have%
$$
(n+m-1){\psi }_\alpha {\varphi }^\alpha =
{\varphi }^\alpha ({\underline{\gamma }}_{\cdot \alpha \beta }^\alpha -{{%
\gamma }^\alpha }_{\alpha \beta })+e{\varphi }^\alpha {\varphi }^\beta ({}^{(%
{\gamma })}{\underline{D}}_\beta q_\alpha -{}^{({\gamma })}D_\beta q_\alpha)
$$

Using the equalities and identities (3.27) and (3.28) we can express
deformations (3.26) as the first $na_{(3)}$--invariant conditions from
(3.15).

To prove the second class of $na_{(3)}$--invariant conditions we introduce
two additional d--tensors:
$$
{{\rho }^\alpha }_{\beta \gamma \delta }=r_{\beta \cdot \gamma \delta
}^{\cdot \alpha }+{\frac 12}({\psi }_{(\beta }{\delta }_{\varphi )}^\alpha +{%
\sigma }_{\beta \varphi }{\varphi }^\tau ){w^\varphi }_{\gamma \delta }{%
\quad }
$$
and
$$
{\underline{\rho }}_{\cdot \beta \gamma \delta }^\alpha ={\underline{r}}%
_{\beta \cdot \gamma \delta }^{\cdot \alpha }-{\frac 12}({\psi }_{(\beta }{%
\delta }_{{\varphi })}^\alpha -{\sigma }_{\beta \varphi }{\varphi }^\tau ){%
w^\varphi }_{\gamma \delta }.\tag3.29
$$
Using deformation (3.26) and (3.29) we write relation
$$
{\tilde \sigma }_{\cdot \beta \gamma \delta }^\alpha ={\underline{\rho }}%
_{\cdot \beta \gamma \delta }^\alpha -{\rho }_{\cdot \beta \gamma \delta
}^\alpha ={\psi }_{\beta [\delta }{\delta }_{\gamma ]}^\alpha -{\psi }_{[{%
\gamma }{\delta }]}{\delta }_\beta ^\alpha -{\sigma }_{\beta \gamma \delta }{%
\varphi }^\alpha ,\tag3.30
$$
where
$$
{\psi }_{\alpha \beta }={}^{({\gamma })}D_\beta {\psi }_\alpha +{\psi }%
_\alpha {\psi }_\beta -({\nu }+{\varphi }^\tau {\psi }_\tau ){\sigma }%
_{\alpha \beta },
$$
and
$$
{\sigma }_{\alpha \beta \gamma }={}^{({\gamma })}D_{[\gamma }{\sigma }_{{%
\beta }]{\alpha }}+{\mu }_{[\gamma }{\sigma }_{{\beta }]{\alpha }}-{\sigma }%
_{{\alpha }[{\gamma }}{\sigma }_{{\beta }]{\tau }}{\varphi }^\tau .
$$
Let multiply (3.30) on $q_\alpha $ and write (taking into account relations
(3.25)) the relation
$$
e{\sigma }_{\alpha \beta \gamma }=-q_\tau {\tilde \sigma }_{\cdot \alpha
\beta \delta }^\tau +{\psi }_{\alpha [\beta }q_{\gamma ]}-{\psi }_{[\beta
\gamma ]}q_\alpha .\tag3.31
$$
The next step is to express ${\psi }_{\alpha \beta }$ trough d--objects on ${%
\xi }.$ To do this we contract indices $\alpha $ and $\beta $ in (3.30) and
obtain
$$
(n+m){\psi }_{[\alpha \beta ]}=-{\sigma }_{\cdot \tau \alpha \beta }^\tau
+eq_\tau {\varphi }^\lambda {\sigma }_{\cdot \lambda \alpha \beta }^\tau -e{%
\tilde \psi }_{[\alpha }{\tilde \psi }_{\beta ]}.
$$
Then contracting indices $\alpha $ and $\delta $ in (3.30) and using (3.31)
we write
$$
(n+m-2){\psi }_{\alpha \beta }={\tilde \sigma }_{\cdot \alpha \beta \tau
}^\tau -eq_\tau {\varphi }^\lambda {\tilde \sigma }_{\cdot \alpha \beta
\lambda }^\tau +{\psi }_{[\beta \alpha ]}+e({\tilde \psi }_\beta q_\alpha -{%
\hat \psi }_{(\alpha }q_{\beta )},\tag3.32
$$
where ${\hat \psi }_\alpha ={\varphi }^\tau {\psi }_{\alpha \tau }.$ If the
both parts of (3.32) are contracted with ${\varphi }^\alpha ,$ it results
that
$$
(n+m-2){\tilde \psi }_\alpha ={\varphi }^\tau {\sigma }_{\cdot \tau \alpha
\lambda }^\lambda -eq_\tau {\varphi }^\lambda {\varphi }^\delta {\sigma }%
_{\lambda \alpha \delta }^\tau -eq_\alpha ,
$$
and, in consequence of ${\sigma }_{\beta (\gamma \delta )}^\alpha =0,$ we
have
$$
(n+m-1){\varphi }={\varphi }^\beta {\varphi }^\gamma {\sigma }_{\cdot \beta
\gamma \alpha }^\alpha .
$$
By using the last expressions we can write
$$
(n+m-2){\underline{\psi }}_\alpha ={\varphi }^\tau {\sigma }_{\cdot \tau
\alpha \lambda }^\lambda -eq_\tau {\varphi }^\lambda {\varphi }^\delta {%
\sigma }_{\cdot \lambda \alpha \delta }^\tau -e{(n+m-1)}^{-1}q_\alpha {%
\varphi }^\tau {\varphi }^\lambda {\sigma }_{\cdot \tau \lambda \delta
}^\delta .\tag3.33
$$
Contracting (3.32) with ${\varphi }^\beta $ we have
$$
(n+m){\hat \psi }_\alpha ={\varphi }^\tau {\sigma }_{\cdot \alpha \tau
\lambda }^\lambda +{\tilde \psi }_\alpha
$$
and taking into consideration (3.33) we can express ${\hat \psi }_\alpha $
through ${\sigma }_{\cdot \beta \gamma \delta }^\alpha .$

As a consequence of (3.31)--(3.33) we obtain this formulas for d--tensor ${%
\psi }_{\alpha \beta }:$
$$
(n+m-2){\psi }_{\alpha \beta }={\sigma }_{\cdot \alpha \beta \tau }^\tau
-eq_\tau {\varphi }^\lambda {\sigma }_{\cdot \alpha \beta \lambda }^\tau +
$$
$$
{\frac 1{n+m}}\{-{\sigma }_{\cdot \tau \beta \alpha }^\tau +eq_\tau {\varphi
}^\lambda {\sigma }_{\cdot \lambda \beta \alpha }^\tau -q_\beta (e{\varphi }%
^\tau {\sigma }_{\cdot \alpha \tau \lambda }^\lambda -q_\tau {\varphi }%
^\lambda {\varphi }^\delta {\sigma }_{\cdot \alpha \lambda \delta }^\tau )+
$$
$$
eq_\alpha [{\varphi }^\lambda {\sigma }_{\cdot \tau \beta \lambda }^\tau
-eq_\tau {\varphi }^\lambda {\varphi }^\delta {\sigma }_{\cdot \lambda \beta
\delta }^\tau -{\frac e{n+m-1}}q_\beta ({\varphi }^\tau {\varphi }^\lambda {%
\sigma }_{\cdot \tau \gamma \delta }^\delta -eq_\tau {\varphi }^\lambda {%
\varphi }^\delta {\varphi }^\varepsilon {\sigma }_{\cdot \lambda \delta
\varepsilon }^\tau )]\}.
$$

Finally, putting the last formula and (3.31) into (3.30) and after a
rearrangement of terms we obtain the second group of $na_{(3)}$-invariant
conditions (3.15). If necessary we can rewrite these conditions in terms of
geometrical objects on $\xi $ and $\underline{\xi }.$ To do this we mast
introduce splittings (3.29) into (3.15). \qquad $\blacksquare$
\vskip15pt

For the particular case of $na_{(3)}$--maps when
$$
{\psi}_{\alpha}=0 , {\varphi}_{\alpha} = g_{\alpha \beta} {\varphi}^{\beta}
= {\frac{\delta }{\delta u^{\alpha}}} ( \ln {\Omega} ) , {\Omega}(u) > 0
$$
and
$$
{\sigma}_{\alpha \beta} = g_{\alpha \beta}%
$$
we define a subclass of conformal transforms ${\underline g}_{\alpha \beta}
(u) = {\Omega}^2 (u)  g_{\alpha \beta}$ which, in consequence of the fact
that d--vector ${\varphi}_{\alpha}$ must satisfy equations (3.11),
generalizes the class of concircular transforms (see \cite{Sinyukov 1979}  for
references and details on concircular mappings of Riemannaian spaces) .

We emphasize that basic na--equations (3.9)--(3.11) are systems of first
order partial differential equations. The study of their geometrical
properties and definition of integral varieties, general and particular
solutions are possible by using the formalism of Pffaf systems
\cite{Cartan 1945}.
Here we point out that by using algebraic methods we can always verify if
systems of na--equations of type (3.9)--(3.11) are, or not, involute, even
to find their explicit solutions it is a difficult task (see more detailed
considerations for isotropic ng--maps in \cite{Sinyukov 1979} and, on
language of
Pffaf systems for na--maps, in \cite{Vacaru 1994}). We can also formulate the
Cauchy problem for na--equations on $\xi $~ and choose deformation
parameters (3.7) as to make involute mentioned equations for the case of
maps to a given background space ${\underline{\xi }}$. If a solution, for
example, of $na_{(1)}$--map equations exists, we say that space $\xi $ is $%
na_{(1)}$--projective to space ${\underline{\xi }}.$ In general, we have to
introduce chains of na--maps in order to obtain involute systems of
equations for maps (superpositions of na-maps) from $\xi $ to ${\underline{%
\xi }}:$
%%%%%%%%%%%%%%%%%%%%%%%%%%%%%%%%%%%%%%%%%%%%%%%%%%%
$$U \buildrel {ng<i_{1}>} \over \longrightarrow  {U_{\underline 1}}
\buildrel ng<i_2> \over \longrightarrow \cdots
\buildrel ng<i_{k-1}> \over \longrightarrow U_{\underline {k-1}}
\buildrel ng<i_k> \over \longrightarrow {\underline U} $$
where $U\subset {\xi },U_{\underline{1}}\subset {\xi }_{\underline{1}%
},\ldots ,U_{k-1}\subset {\xi }_{k-1},{\underline{U}}\subset {\xi }_k$ with
corresponding splittings of auxiliary symmetric connections
$$
{\underline{\gamma }}_{.{\beta }{\gamma }}^\alpha =_{<i_1>}P_{.{\beta }{%
\gamma }}^\alpha +_{<i_2>}P_{.{\beta }{\gamma }}^\alpha +\cdots +_{<i_k>}P_{.%
{\beta }{\gamma }}^\alpha
$$
and torsion
$$
{\underline{T}}_{.{\beta }{\gamma }}^\alpha =T_{.{\beta }{\gamma }}^\alpha
+_{<i_1>}Q_{.{\beta }{\gamma }}^\alpha +_{<i_2>}Q_{.{\beta }{\gamma }%
}^\alpha +\cdots +_{<i_k>}Q_{.{\beta }{\gamma }}^\alpha
$$
where cumulative indices $<i_1>=0,1,2,3,$ denote possible types of na--maps.

\proclaim{\bf Definition 3.4}
Space $\xi $~ is nearly conformally projective to space ${\underline{\xi }},{%
\quad }nc:{\xi }{\to }{\underline{\xi }},$~ if there is a finite chain of
na--maps from $\xi $~ to ${\underline{\xi }}.$
\endproclaim

For nearly conformal maps we formulate :

\proclaim{\bf Theorem 3.4}
For every fixed triples $(N_j^a,{\Gamma }_{{.}{\beta }{\gamma }}^\alpha
,U\subset {\xi }$ and $(N_j^a,{\underline{\Gamma }}_{{.}{\beta }{\gamma }%
}^\alpha $, ${\underline{U}}\subset {\underline{\xi }})$, components of
nonlinear connection, d--connection and d--metric being of class $C^r(U),C^r(%
{\underline{U}})$, $r>3,$ there is a finite chain of na--maps $nc:U\to {%
\underline{U}}.$
\endproclaim

Proof is similar to that for isotropic maps \cite{Vacaru 1994} (we have to
introduce a finite number of na-maps with corresponding components of
deformation parameters and deformation tensors in order to transform step
by step coefficients of d-connection ${\Gamma}^{\alpha}_{\gamma \delta}$
into ${\underline {\Gamma}}^{\alpha}_{\beta \gamma} ).$

Now we introduce the concept of the Category of la--spaces, ${\Cal C}({\xi}).
$ The elements of ${\Cal C}({\xi})$ consist from $Ob{\Cal C}({\xi})=\{{\xi},
{\xi}_{<i_{1}>}, {\xi}_{<i_{2}>},{\ldots}, \}$ being la--spaces, for
simplicity in this work, having common N--connection structures, and $Mor
{\Cal C}({\xi})=\{ nc ({\xi}_{<i_{1}>}, {\xi}_{<i_{2}>})\}$ being chains of
na--maps interrelating la--spaces. We point out that we can consider
equivalent models of physical theories on every object of ${\Cal C}({\xi})$
(see details for isotropic gravitational models in \cite{Petrov 1970} and
 \cite{Vacaru 1994} and anisotropic gravity in \cite{Vacau and Ostaf 1994,
 1996a, 1996b}). One of the main purposes of this chapter is to develop a
d--tensor and variational formalism on ${\Cal C}({\xi}),$ i.e. on
la--multispaces, interrelated with nc--maps. Taking into account the
distinguished character of geometrical objects on la--spaces we call tensors
on ${\Cal C}({\xi})$ as distinguished tensors on la--space Category, or
dc--tensors.

Finally, we emphasize that presented in that section definitions and
theorems can be generalized for v--bundles with arbitrary given structures
of nonlinear connection, linear d--connection and metric structures.
We omit such combersome calculations connected with deformation of all
 basic N--connection, d--connection and d--metric structures.
\vskip25pt

\head{III.3 Nearly Autoparallel Tensor--Integral on LA--Spaces}\endhead

The aim of this section is to define tensor integration not only for
bitensors, objects defined on the same curved space, but for dc--tensors,
defined on two spaces, $\xi$ and ${\underline {\xi}}$, even it is necessary
on la--multispaces.

Let $T_{u}{\xi}$~ and $T_{\underline u}{\underline {\xi}}$ be tangent spaces
in corresponding points $u {\in} U {\subset} {\xi}$ and ${\underline u} {\in}
{\underline U} {\subset} {\underline {\xi}}$ and, respectively, $T^{\ast}_{u}%
{\xi}$ and $T^{\ast}_{\underline u}{\underline {\xi}} $ be their duals (in
general, in this section we shall not consider that a common
coordinatization is introduced for open regions $U$ and ${\underline U}$ ).
We call as the dc--tensors on the pair of spaces (${\xi}, {\underline {\xi}}$
) the elements of distinguished tensor algebra
$$
( {\otimes}_{\alpha} T_{u}{\xi}) {\otimes} ({\otimes}_{\beta} T^{\ast}_{u} {%
\xi}) {\otimes}({\otimes}_{\gamma} T_{\underline u}{\underline {\xi}}) {%
\otimes} ({\otimes}_{\delta} T^{\ast}_{\underline u} {\underline {\xi}})%
$$
defined over the space ${\xi}{\otimes} {\underline {\xi}}, $ for a given $nc
: {\xi} {\to} {\underline {\xi}} $.

We admit the convention that underlined and non--underlined indices refer,
respectively, to the points ${\underline u}$ and $u$. Thus $Q_{{.}{%
\underline {\alpha}}}^{\beta}, $ for instance, are the components of
dc--tensor $Q{\in} T_{u}{\xi} {\otimes} T_{\underline u}{\underline {\xi}}.$

Now, we define the transport dc--tensors. Let open regions $U$ and ${%
\underline U}$ be homeomorphic to sphere ${\Bbb R}^{2n}$ and introduce
isomorphism ${\mu}_{{u},{\underline u}}$ between $T_{u}{\xi}$ and $%
T_{\underline u}{\underline {\xi}}$ (given by map $nc : U {\to} {\underline U%
}).$ We consider that for every d--vector $v^{\alpha} {\in} T_{u}{\xi}$
corresponds the vector ${\mu}_{{u},{\underline u}}
(v^{\alpha})=v^{\underline {\alpha}} {\in} T_{\underline u}{\underline {\xi}}%
,$ with components $v^{\underline {\alpha}}$ being linear functions of $%
v^{\alpha}$:
$$
v^{\underline {\alpha}}=h^{\underline {\alpha}}_{\alpha}(u, {\underline u})
v^{\alpha}, \quad v_{\underline {\alpha}}= h^{\alpha}_{\underline {\alpha}}({%
\underline u}, u)v_{\alpha},%
$$
where $h^{\alpha}_{\underline {\alpha}}({\underline u}, u)$ are the
components of dc--tensor associated with ${\mu}^{-1}_{u,{\underline u}}$. In
a similar manner we have
$$
v^{\alpha}=h^{\alpha}_{\underline {\alpha}}({\underline u}, u) v^{\underline
{\alpha}}, \quad v_{\alpha}=h^{\underline {\alpha}}_{\alpha} (u, {\underline
u})v_{\underline {\alpha}}.%
$$

In order to reconcile just presented definitions and to assure the identity
for trivial maps ${\xi }{\to }{\xi },u={\underline{u}},$ the transport
dc--tensors must satisfy conditions:
$$
h_\alpha ^{\underline{\alpha }}(u,{\underline{u}})h_{\underline{\alpha }%
}^\beta ({\underline{u}},u)={\delta }_\alpha ^\beta ,h_\alpha ^{\underline{%
\alpha }}(u,{\underline{u}})h_{\underline{\beta }}^\alpha ({\underline{u}}%
,u)={\delta }_{\underline{\beta }}^{\underline{\alpha }}
$$
and ${\lim }_{{({\underline{u}}{\to }u})}h_\alpha ^{\underline{\alpha }}(u,{%
\underline{u}})={\delta }_\alpha ^{\underline{\alpha }},\quad {\lim }_{{({%
\underline{u}}{\to }u})}h_{\underline{\alpha }}^\alpha ({\underline{u}},u)={%
\delta }_{\underline{\alpha }}^\alpha .$

Let ${\overline S}_{p} {\subset} {\overline U} {\subset} {\overline {\xi}}$
is a homeomorphic to $p$-dimensional sphere and suggest that chains of
na--maps are used to connect regions :
$$ U \buildrel nc_{(1)} \over \longrightarrow {\overline S}_p
     \buildrel nc_{(2)} \over \longrightarrow {\underline U}.$$

\proclaim{\bf Definition 3.5}
The tensor integral in ${\overline{u}}{\in }{\overline{S}}_p$ of a
dc--tensor $N_{{\varphi }{.}{\underline{\tau }}{.}{\overline{\alpha }}_1{%
\cdots }{\overline{\alpha }}_p}^{{.}{\gamma }{.}{\underline{\kappa }}}$ $({%
\overline{u}},u),$ completely antisymmetric on the indices ${{\overline{%
\alpha }}_1},{\ldots },{\overline{\alpha }}_p,$ over domain ${\overline{S}}%
_p,$ is defined as
$$\multline
N_{{\varphi }{.}{\underline{\tau }}}^{{.}{\gamma }{.}{\underline{\kappa }}}({%
\underline{u}},u)=I_{({\overline{S}}_p)}^{\underline{U}}N_{{\varphi }{.}{%
\overline{\tau }}{.}{\overline{\alpha }}_1{\ldots }{\overline{\alpha }}_p}^{{%
.}{\gamma }{.}{\overline{\kappa }}}({\overline{u}},{\underline{u}})dS^{{%
\overline{\alpha }}_1{\ldots }{\overline{\alpha }}_p}=\\
{\int }_{({\overline{S}}_p)}h_{\underline{\tau }}^{\overline{\tau }}({%
\underline{u}},{\overline{u}})h_{\overline{\kappa }}^{\underline{\kappa }}({%
\overline{u}},{\underline{u}})N_{{\varphi }{.}{\overline{\tau }}{.}{%
\overline{\alpha }}_1{\cdots }{\overline{\alpha }}_p}^{{.}{\gamma }{.}{%
\overline{\kappa }}}({\overline{u}},u)d{\overline{S}}^{{\overline{\alpha }}_1%
{\cdots }{\overline{\alpha }}_p}, \endmultline \tag3.34
$$
where $dS^{{\overline{\alpha }}_1{\cdots }{\overline{\alpha }}_p}={\delta }%
u^{{\overline{\alpha }}_1}{\land }{\cdots }{\land }{\delta }u_p^{\overline{%
\alpha }}$.
\endproclaim

Let suppose that transport dc--tensors $h_\alpha ^{\underline{\alpha }}$~
and $h_{\underline{\alpha }}^\alpha $~ admit covariant derivations of
or\-der two and pos\-tu\-la\-te ex\-is\-ten\-ce of de\-for\-ma\-ti\-on
dc--ten\-sor $B_{{\alpha }{\beta }}^{{..}{\gamma }}(u,{\underline{u}})$~
satisfying relations
$$
D_\alpha h_\beta ^{\underline{\beta }}(u,{\underline{u}})=B_{{\alpha }{\beta
}}^{{..}{\gamma }}(u,{\underline{u}})h_\gamma ^{\underline{\beta }}(u,{%
\underline{u}})\tag3.35
$$
and, taking into account that $D_\alpha {\delta }_\gamma ^\beta =0,$
$$
D_\alpha h_{\underline{\beta }}^\beta ({\underline{u}},u)=-B_{{\alpha }{%
\gamma }}^{{..}{\beta }}(u,{\underline{u}})h_{\underline{\beta }}^\gamma ({%
\underline{u}},u).
$$
By using formulas  for torsion and, respectively,
curvature of connection ${\Gamma }_{{\beta }{\gamma }}^\alpha $~ we can
calculate next commutators:
$$
D_{[{\alpha }}D_{{\beta }]}h_\gamma ^{\underline{\gamma }}=-(R_{{\gamma }{.}{%
\alpha }{\beta }}^{{.}{\lambda }}+T_{{.}{\alpha }{\beta }}^\tau B_{{\tau }{%
\gamma }}^{{..}{\lambda }})h_\lambda ^{\underline{\gamma }}.\tag3.36
$$
On the other hand from (3.35) one follows that
$$
D_{[{\alpha }}D_{{\beta }]}h_\gamma ^{\underline{\gamma }}=(D_{[{\alpha }}B_{%
{\beta }]{\gamma }}^{{..}{\lambda }}+B_{[{\alpha }{|}{\tau }{|}{.}}^{{..}{%
\lambda }}B_{{\beta }]{\gamma }{.}}^{{..}{\tau }})h_\lambda ^{\underline{%
\gamma }},\tag3.37
$$
where ${|}{\tau }{|}$~ denotes that index ${\tau }$~ is excluded from the
action of antisymmetrization $[{\quad }]$. From (3.36) and (3.37) we obtain
$$
D_{[{\alpha }}B_{{\beta }]{\gamma }{.}}^{{..}{\lambda }}+B_{[{\beta }{|}{%
\gamma }{|}}B_{{\alpha }]{\tau }}^{{..}{\lambda }}=(R_{{\gamma }{.}{\alpha }{%
\beta }}^{{.}{\lambda }}+T_{{.}{\alpha }{\beta }}^\tau B_{{\tau }{\gamma }}^{%
{..}{\lambda }}).\tag3.38
$$

Let ${\overline{S}}_p$~ be the boundary of ${\overline{S}}_{p-1}$. The
Stoke's type formula for tensor--integral (3.34) is defined as
$$
I_{{\overline{S}}_p}N_{{\varphi }{.}{\overline{\tau }}{.}{\overline{\alpha }}%
_1{\ldots }{\overline{\alpha }}_p}^{{.}{\gamma }{.}{\overline{\kappa }}}dS^{{%
\overline{\alpha }}_1{\ldots }{\overline{\alpha }}_p}=
I_{{\overline{S}}_{p+1}}{^{{\star }{(p)}}{\overline{D}}}_{[{\overline{\gamma
}}{|}}N_{{\varphi }{.}{\overline{\tau }}{.}{|}{\overline{\alpha }}_1{\ldots }%
{{\overline{\alpha }}_p]}}^{{.}{\gamma }{.}{\overline{\kappa }}}dS^{{%
\overline{\gamma }}{\overline{\alpha }}_1{\ldots }{\overline{\alpha }}_p},%
\tag3.39
$$
where
$$\multline
{^{{\star }{(p)}}D}_{[{\overline{\gamma }}{|}}N_{{\varphi }{.}{\overline{%
\tau }}{.}{|}{\overline{\alpha }}_1{\ldots }{\overline{\alpha }}_p]}^{{.}{%
\gamma }{.}{\overline{\kappa }}}=
D_{[{\overline{\gamma }}{|}}N_{{\varphi }{.}{\overline{\tau }}{.}{|}{%
\overline{\alpha }}_1{\ldots }{\overline{\alpha }}_p]}^{{.}{\gamma }{.}{%
\overline{\kappa }}}+\\
 pT_{{.}[{\overline{\gamma }}{\overline{\alpha }}_1{|}}^{%
\underline{\epsilon }}N_{{\varphi }{.}{\overline{\tau }}{.}{\overline{%
\epsilon }}{|}{\overline{\alpha }}_2{\ldots }{\overline{\alpha }}_p]}^{{.}{%
\gamma }{.}{\overline{\kappa }}}-B_{[{\overline{\gamma }}{|}{\overline{\tau }%
}}^{{..}{\overline{\epsilon }}}N_{{\varphi }{.}{\overline{\epsilon }}{.}{|}{%
\overline{\alpha }}_1{\ldots }{\overline{\alpha }}_p]}^{{.}{\gamma }{.}{%
\overline{\kappa }}}+B_{[{\overline{\gamma }}{|}{\overline{\epsilon }}}^{..{%
\overline{\kappa }}}N_{{\varphi }{.}{\overline{\tau }}{.}{|}{\overline{%
\alpha }}_1{\ldots }{\overline{\alpha }}_p]}^{{.}{\gamma }{.}{\overline{%
\epsilon }}}. \endmultline \tag3.40
$$
We define the dual element of the hypersurfaces element $dS^{{j}_1{\ldots }{j%
}_p}$ as
$$
d{\Cal S}_{{\beta }_1{\ldots }{\beta }_{q-p}}={\frac 1{{p!}}}{\epsilon }_{{%
\beta }_1{\ldots }{\beta }_{k-p}{\alpha }_1{\ldots }{\alpha }_p}dS^{{\alpha }%
_1{\ldots }{\alpha }_p},\tag3.41
$$
where ${\epsilon }_{{\gamma }_1{\ldots }{\gamma }_q}$ is completely
antisymmetric on its indices and
$$
{\epsilon }_{12{\ldots }(n+m)}=\sqrt{{|}G{|}},G=det{|}G_{{\alpha }{\beta }{|}%
},
$$
$G_{{\alpha }{\beta }}$ is taken from (1.12). The dual of dc--tensor $N_{{%
\varphi }{.}{\overline{\tau }}{.}{\overline{\alpha }}_1{\ldots }{\overline{%
\alpha }}_p}^{{.}{\gamma }{\overline{\kappa }}}$ is defined as the
dc--tensor  ${\Cal N}_{{\varphi }{.}{\overline{\tau }}}^{{.}{\gamma }{.}{%
\overline{\kappa }}{\overline{\beta }}_1{\ldots }{\overline{\beta }}_{n+m-p}}
$ satisfying
$$
N_{{\varphi }{.}{\overline{\tau }}{.}{\overline{\alpha }}_1{\ldots }{%
\overline{\alpha }}_p}^{{.}{\gamma }{.}{\overline{\kappa }}}={\frac 1{{p!}}}%
{\Cal N}_{{\varphi }{.}{\overline{\tau }}}^{{.}{\gamma }{.}{\overline{\kappa
}}{\overline{\beta }}_1{\ldots }{\overline{\beta }}_{n+m-p}}{\epsilon }_{{%
\overline{\beta }}_1{\ldots }{\overline{\beta }}_{n+m-p}{\overline{\alpha }}%
_1{\ldots }{\overline{\alpha }}_p}.\tag3.42
$$
Using (3.16), (3.41) and (3.42) we can write
$$
I_{{\overline{S}}_p}N_{{\varphi }{.}{\overline{\tau }}{.}{\overline{\alpha }}%
_1{\ldots }{\overline{\alpha }}_p}^{{.}{\gamma }{.}{\overline{\kappa }}}dS^{{%
\overline{\alpha }}_1{\ldots }{\overline{\alpha }}_p}={\int }_{{\overline{S}}%
_{p+1}}{^{\overline{p}}D}_{\overline{\gamma }}{\Cal N}_{{\varphi }{.}{%
\overline{\tau }}}^{{.}{\gamma }{.}{\overline{\kappa }}{\overline{\beta }}_1{%
\ldots }{\overline{\beta }}_{n+m-p-1}{\overline{\gamma }}}d{\Cal S}_{{%
\overline{\beta }}_1{\ldots }{\overline{\beta }}_{n+m-p-1}},\tag3.43
$$
where
$$
{^{\overline{p}}D}_{\overline{\gamma }}{\Cal N}_{{\varphi }{.}{\overline{%
\tau }}}^{{.}{\gamma }{.}{\overline{\kappa }}{\overline{\beta }}_1{\ldots }{%
\overline{\beta }}_{n+m-p-1}{\overline{\gamma }}}=
$$
$$
{\overline{D}}_{\overline{\gamma }}{\Cal N}_{{\varphi }{.}{\overline{\tau }}%
}^{{.}{\gamma }{.}{\overline{\kappa }}{\overline{\beta }}_1{\ldots }{%
\overline{\beta }}_{n+m-p-1}{\overline{\gamma }}}+(-1)^{(n+m-p)}(n+m-p+1)T_{{%
.}{\overline{\gamma }}{\overline{\epsilon }}}^{[{\overline{\epsilon }}}{\Cal %
N}_{{\varphi }{.}{{\overline{\tau }}}}^{{.}{|}{\gamma }{.}{\overline{\kappa }%
}{|}{\overline{\beta }}_1{\ldots }{\overline{\beta }}_{n+m-p-1}]{\overline{%
\gamma }}}-
$$
$$
B_{{\overline{\gamma }}{\overline{\tau }}}^{{..}{\overline{\epsilon }}}{\Cal %
N}_{{\varphi }{.}{\overline{\epsilon }}}^{{.}{\gamma }{.}{\overline{\kappa }}%
{\overline{\beta }}_1{\ldots }{\overline{\beta }}_{n+m-p-1}{\overline{\gamma
}}}+B_{{\overline{\gamma }}{\overline{\epsilon }}}^{{..}{\overline{\kappa }}}%
{\Cal N}_{{\varphi }{.}{\overline{\tau }}}^{{.}{\gamma }{.}{\overline{%
\epsilon }}{\overline{\beta }}_1{\ldots }{\overline{\beta }}_{n+m-p-1}{%
\overline{\gamma }}}.
$$
To verify the equivalence of (3.42) and (3.43) we must take in consideration
that
$$
D_\gamma {\epsilon }_{{\alpha }_1{\ldots }{\alpha }_k}=0\ \text{and}\ {%
\epsilon }_{{\beta }_1{\ldots }{\beta }_{n+m-p}{\alpha }_1{\ldots }{\alpha }%
_p}{\epsilon }^{{\beta }_1{\ldots }{\beta }_{n+m-p}{\gamma }_1{\ldots }{%
\gamma }_p}=p!(n+m-p)!{\delta }_{{\alpha }_1}^{[{\gamma }_1}{\cdots }{\delta
}_{{\alpha }_p}^{{\gamma }_p]}.
$$
The developed in this section tensor integration formalism will be used in
the next section for definition of conservation laws on spaces with local
anisotropy.
\vskip25pt

\head{III.4\ On Conservation Laws on La--Spaces}\endhead

To define conservation laws on locally anisotropic spaces is a challenging
task because of absence of global and local groups of automorphisms of such
spaces. Our main idea is to use chains of na--maps from a given, called
hereafter as the fundamental la--space to an auxiliary one with trivial
curvatures and torsions admitting a global group of automorphisms. The aim
of this section is to  formulate conservation laws for la-gravitational
fields by using dc--objects and tensor--integral values, na--maps and
variational calculus on the Category of la--spaces.
\vskip15pt

\subhead{III.4.1\ Nonzero divergence of the energy--momentum d--tensor}
\endsubhead

In work \cite{Miron and Anastasiei 1994} it is pointed to this specific form of
conservation laws of matter on la--spaces: They calculated the divergence of
the energy--momentum d--tensor on la--space $\xi ,$%
$$
D_\alpha {E}_\beta ^\alpha ={\frac 1{\ {\kappa }_1}}U_\alpha ,\tag3.44
$$
and concluded that d--vector
$$
U_\alpha ={\frac 12}(G^{\beta \delta }{{R_\delta }^\gamma }_{\phi \beta }{T}%
_{\cdot \alpha \gamma }^\phi -G^{\beta \delta }{{R_\delta }^\gamma }_{\phi
\alpha }{T}_{\cdot \beta \gamma }^\phi +{R_\phi ^\beta }{T}_{\cdot \beta
\alpha }^\phi )
$$
vanishes if and only if d--connection $D$ is without torsion.

No wonder that conservation laws, in usual physical theories being a
consequence of global (for usual gravity of local) automorphisms of the
fundamental space--time, are more sophisticate on the spaces with local
anisotropy. Here it is important to emphasize the multiconnection character
of la--spaces. For example, for a d--metric (1.12) on $\xi $ we can
equivalently introduce another (see (1.23)) metric linear connection $%
\tilde D .$ The Einstein equations
$$
{\tilde R}_{\alpha \beta }-{\frac 12}G_{\alpha \beta }{\tilde R}={\kappa }_1{%
\tilde E}_{\alpha \beta }\tag3.45
$$
constructed by using connection (3.23) have vanishing divergences
$$
{\tilde D}^\alpha ({{\tilde R}_{\alpha \beta }}-{\frac 12}G_{\alpha \beta }{%
\tilde R})=0\text{ and }{\tilde D}^\alpha {\tilde E}_{\alpha \beta }=0,
$$
similarly as those on (pseudo)Riemannian spaces. We conclude that by using
the connection (1.23) we construct a model of la--gravity which looks like
locally isotropic on the total space $E.$  More general gravitational models
with local anisot\-ro\-py can be obtained by using deformations of connection
${\tilde \Gamma }_{\cdot \beta \gamma }^\alpha ,$
$$
{{\Gamma }^\alpha }_{\beta \gamma }={\tilde \Gamma }_{\cdot \beta \gamma
}^\alpha +{P^\alpha }_{\beta \gamma }+{Q^\alpha }_{\beta \gamma },
$$
were, for simplicity, ${{\Gamma }^\alpha }_{\beta \gamma }$ is chosen to be
also metric and satisfy Einstein equations (3.45). We can consider
deformation d--tensors ${P^\alpha }_{\beta \gamma }$ generated (or not) by
deformations of type (3.9)--(3.11) for na--maps. In
this case d--vector $U_\alpha $ can be interpreted as a generic source of
local anisotropy on $\xi $ satisfying generalized conservation laws (3.44).
\vskip15pt
\subhead{III.4.2\ Deformation d--tensors and tensor--integral conservation laws}
\endsubhead

From (3.34) we obtain a tensor integral on ${\Cal C}({\xi})$ of a d--tensor:
$$
N_{{\underline {\tau}}}^{{.}{\underline {\kappa}}}(\underline u)= I_{{%
\overline S}_{p}}N^{{..}{\overline {\kappa}}}_ {{\overline {\tau}}{..}{%
\overline {\alpha}}_{1}{\ldots}{\overline {\alpha}}_{p}} ({\overline u})h^{{%
\overline {\tau}}}_{{\underline {\tau}}}({\underline u}, {\overline u})h^{{%
\underline {\kappa}}}_{{\overline {\kappa}}} ({\overline u}, {\underline u}%
)dS^{{\overline {\alpha}}_{1}{\ldots} {\overline {\alpha}}_{p}}.%
$$

We point out that tensor--integral can be defined not only for dc--tensors
but and for d--tensors on $\xi $. Really, suppressing indices ${\varphi }$~
and ${\gamma }$~ in (3.42) and (3.43), considering instead of a deformation
dc--tensor a deformation tensor
$$
B_{{\alpha }{\beta }}^{{..}{\gamma }}(u,{\underline{u}})=B_{{\alpha }{\beta }%
}^{{..}{\gamma }}(u)=P_{{.}{\alpha }{\beta }}^\gamma (u)\tag3.46
$$
(we consider deformations induced by a nc--transform) and integration\newline
 $I_{S_p}{\ldots }dS^{{\alpha }_1{\ldots }{\alpha }_p}$ in la--space $\xi $ we
obtain from (3.34) a tensor--integral on ${\Cal C}({\xi })$~ of a
d--tensor:
$$
N_{{\underline{\tau }}}^{{.}{\underline{\kappa }}}({\underline{u}}%
)=I_{S_p}N_{{\tau }{.}{\alpha }_1{\ldots }{\alpha }_p}^{.{\kappa }}(u)h_{{%
\underline{\tau }}}^\tau ({\underline{u}},u)h_\kappa ^{\underline{\kappa }%
}(u,{\underline{u}})dS^{{\alpha }_1{\ldots }{\alpha }_p}.
$$
Taking into account (3.38) we can calculate that curvature
$$
{\underline{R}}_{{\gamma }{.}{\alpha }{\beta }}^{.{\lambda }}=D_{[{\beta }%
}B_{{\alpha }]{\gamma }}^{{..}{\lambda }}+B_{[{\alpha }{|}{\gamma }{|}}^{{..}%
{\tau }}B_{{\beta }]{\tau }}^{{..}{\lambda }}+T_{{.}{\alpha }{\beta }}^{{%
\tau }{..}}B_{{\tau }{\gamma }}^{{..}{\lambda }}
$$
of connection ${\underline{\Gamma }}_{{.}{\alpha }{\beta }}^\gamma (u)={%
\Gamma }_{{.}{\alpha }{\beta }}^\gamma (u)+B_{{\alpha }{\beta }{.}}^{{..}{%
\gamma }}(u),$ with $B_{{\alpha }{\beta }}^{{..}{\gamma }}(u)$~ taken from
(3.46), vanishes, ${\underline{R}}_{{\gamma }{.}{\alpha }{\beta }}^{{.}{%
\lambda }}=0.$ So, we can conclude that la--space $\xi $ admits a tensor
integral structure on ${\Cal {C}}({\xi })$ for d--tensors associated to
deformation tensor $B_{{\alpha }{\beta }}^{{..}{\gamma }}(u)$ if the
nc--image ${\underline{\xi }}$~ is locally parallelizable. That way we
generalize the one space tensor integral constructions in \cite{Gottlieb and
 Vacaru 1996},
were the possibility to introduce tensor integral structure on a curved
space was restricted by the condition that this space is locally
parallelizable. For $q=n+m$~ relations (3.43), written for d--tensor ${\Cal N%
}_{\underline{\alpha }}^{{.}{\underline{\beta }}{\underline{\gamma }}}$ (we
change indices ${\overline{\alpha }},{\overline{\beta }},{\ldots }$ into ${%
\underline{\alpha }},{\underline{\beta }},{\ldots })$ extend the Gauss
formula on ${\Cal {C}}({\xi })$:
$$
I_{S_{q-1}}{\Cal N}_{\underline{\alpha }}^{{.}{\underline{\beta }}{%
\underline{\gamma }}}d{\Cal S}_{\underline{\gamma }}=I_{{\underline{S}}_q}{^{%
\underline{q-1}}D}_{{\underline{\tau }}}{\Cal N}_{{\underline{\alpha }}}^{{.}%
{\underline{\beta }}{\underline{\tau }}}d{\underline{V}},\tag3.47
$$
where $d{\underline{V}}={\sqrt{{|}{\underline{G}}_{{\alpha }{\beta }}{|}}}d{%
\underline{u}}^1{\ldots }d{\underline{u}}^q$ and
$$
{^{\underline{q-1}}D}_{{\underline{\tau }}}{\Cal N}_{\underline{\alpha }}^{{.%
}{\underline{\beta }}{\underline{\tau }}}=D_{{\underline{\tau }}}{\Cal N}_{%
\underline{\alpha }}^{{.}{\underline{\beta }}{\underline{\tau }}}-T_{{.}{%
\underline{\tau }}{\underline{\epsilon }}}^{{\underline{\epsilon }}}{\Cal N}%
_{{\underline{\alpha }}}^{{\underline{\beta }}{\underline{\tau }}}-B_{{%
\underline{\tau }}{\underline{\alpha }}}^{{..}{\underline{\epsilon }}}{\Cal N%
}_{{\underline{\epsilon }}}^{{.}{\underline{\beta }}{\underline{\tau }}}+B_{{%
\underline{\tau }}{\underline{\epsilon }}}^{{..}{\underline{\beta }}}{\Cal N}%
_{{\underline{\alpha }}}^{{.}{\underline{\epsilon }}{\underline{\tau }}}.%
\tag3.48
$$

Let consider physical values $N_{{\underline{\alpha }}}^{{.}{\underline{%
\beta }}}$ on ${\underline{\xi }}$~ defined on its density ${\Cal N}_{{%
\underline{\alpha }}}^{{.}{\underline{\beta }}{\underline{\gamma }}},$ i. e.
$$
N_{{\underline{\alpha }}}^{{.}{\underline{\beta }}}=I_{{\underline{S}}_{q-1}}%
{\Cal N}_{{\underline{\alpha }}}^{{.}{\underline{\beta }}{\underline{\gamma }%
}}d{\Cal S}_{{\underline{\gamma }}}\tag3.49
$$
with this conservation law (due to (3.47)):%
$$
{^{\underline{q-1}}D}_{{\underline{\gamma }}}{\Cal N}_{{\underline{\alpha }}%
}^{{.}{\underline{\beta }}{\underline{\gamma }}}=0.\tag3.50
$$
We note that these conservation laws differ from covariant conservation laws
for well known physical values such as density of electric current or of
energy-- momentum tensor. For example, taking density ${E}_\beta ^{{.}{%
\gamma }},$ with corresponding to (3.48) and (3.50) conservation law,
$$
{^{\underline{q-1}}D}_{{\underline{\gamma }}}{E}_{{\underline{\beta }}}^{{%
\underline{\gamma }}}=D_{{\underline{\gamma }}}{E}_{{\underline{\beta }}}^{{%
\underline{\gamma }}}-T_{{.}{\underline{\epsilon }}{\underline{\tau }}}^{{%
\underline{\tau }}}{E}_{{\underline{\beta }}}^{{.}{\underline{\epsilon }}%
}-B_{{\underline{\tau }}{\underline{\beta }}}^{{..}{\underline{\epsilon }}}{E%
}_{\underline{\epsilon }}^{{\underline{\tau }}}=0,\tag3.51
$$
we can define values (see (3.47) and (3.49))
$$
{\Cal P}_\alpha =I_{{\underline{S}}_{q-1}}{E}_{{\underline{\alpha }}}^{{.}{%
\underline{\gamma }}}d{\Cal S}_{{\underline{\gamma }}}.
$$
The defined conservation laws (3.51) for ${E}_{{\underline{\beta }}}^{{.}{%
\underline{\epsilon }}}$ have nothing to do with those for energy--momentum
tensor $E_\alpha ^{{.}{\gamma }}$ from Einstein equations for the almost
Hermitian gravity \cite{Miron and Anastasieie 1994}
 or with ${\tilde E}_{\alpha \beta }$ from
(3.45) with vanishing divergence $D_\gamma {\tilde E}_\alpha ^{{.}{\gamma }%
}=0.$ So ${\tilde E}_\alpha ^{{.}{\gamma }}{\neq }{E}_\alpha ^{{.}{\gamma }}.
$ A similar conclusion was made in \cite{Gottlieb and Vacaru 1996}
 for unispacial locally
isotropic tensor integral. In the case of multispatial tensor integration we
have another possibility, namely, to identify ${E}_{{\underline{\beta }}}^{{.}%
{\underline{\gamma }}}$ from (3.51) with the na-image of ${E}_\beta ^{{.}{%
\gamma }}$ on la--space $\xi .$ We shall consider this construction in the
next section.
\vskip25pt

\head{III.5\ NA--Conservation Laws in LA--Gravity}\endhead

It is well known that the standard pseudo--tensor description of the
energy--momen\-tum values for the Einstein gravitational fields is full of
ambiguities. Some light can be shed by introducing additional geometrical
structures on curved space--time (bimetrics, biconnections, by taking into
 account background spaces, or formulating variants of general relativity
 theory on flat space).  We emphasize here that rigorous mathematical
 investigations
based on two (fundamental and background) locally anisotropic, or isotropic,
spaces should use well--defined, motivated from physical point of view,
mappings of these spaces. Our na--model largely contains both attractive
features of the mentioned approaches; na--maps establish a local 1--1
correspondence between the fundamental la--space and auxiliary la--spaces on
which biconnection (or even multiconnection) structures are induced. But
these structures are not a priory postulated as in a lot of gravitational
theories, we tend to specify them to be locally reductible to the locally
isotropic Einstein theory.

Let us consider a fixed background la--space $\underline {\xi}$ with given
metric ${\underline G}_{\alpha \beta} = ({\underline g}_{ij} , {\underline h}%
_{ab} )$ and d--connection ${\underline {\tilde {\Gamma}}}^{\alpha}_{\cdot
\beta\gamma},$ for simplicity in this subsection we consider compatible
metric and connections  being torsionless and with vanishing curvatures.
Supposing that there is an nc--transform from the fundamental la--space $\xi$
to the auxiliary $\underline {\xi} .$ we are interested in the equivalents
of the Einstein equations (3.45) on $\underline {\xi} .$

We consider that a part of gravitational degrees of freedom is "pumped out"
into the dynamics of deformation d--tensors for d--connection, ${P^{\alpha}}%
_{\beta \gamma},$ and metric, $B^{\alpha \beta} =  ( b^{ij} , b^{ab} ) .$
The remained part of degrees of freedom is coded into the metric ${%
\underline G}_{\alpha \beta}$ and d--connection ${\underline {\tilde
{\Gamma}}}^{\alpha}_{\cdot \beta \gamma} .$

We apply the first order formalism and consider $%
B^{\alpha \beta }$ and ${P^\alpha }_{\beta \gamma }$ as independent
variables on $\underline{\xi }.$ Using notations
$$
P_\alpha ={P^\beta }_{\beta \alpha },\quad {\Gamma }_\alpha ={{\Gamma }%
^\beta }_{\beta \alpha },\
{\hat B}^{\alpha \beta }=\sqrt{|G|}B^{\alpha \beta },{\hat G}^{\alpha \beta
}=\sqrt{|G|}G^{\alpha \beta },{\underline{\hat G}}^{\alpha \beta }=\sqrt{|%
\underline{G}|}{\underline{G}}^{\alpha \beta }
$$
and making identifications
$$
{\hat B}^{\alpha \beta }+{\underline{\hat G}}^{\alpha \beta }={\hat G}%
^{\alpha \beta },{\quad }{\underline{\Gamma }}_{\cdot \beta \gamma }^\alpha -%
{P^\alpha }_{\beta \gamma }={{\Gamma }^\alpha }_{\beta \gamma },
$$
we take the action of la--gravitational field on $\underline{\xi }$ in this
form:
$$
{\underline{{\Cal S}}}^{(g)}=-{(2c{\kappa }_1)}^{-1}\int {\delta }^qu{}{%
\underline{{\Cal L}}}^{(g)},\tag3.52
$$
where
$$
{\underline{{\Cal L}}}^{(g)}={\hat B}^{\alpha \beta }(D_\beta P_\alpha
-D_\tau {P^\tau }_{\alpha \beta })+({\underline{\hat G}}^{\alpha \beta }+{%
\hat B}^{\alpha \beta })(P_\tau {P^\tau }_{\alpha \beta }-{P^\alpha }%
_{\alpha \kappa }{P^\kappa }_{\beta \tau })
$$
and the interaction constant is taken ${\kappa }_1={\frac{4{\pi }}{{c^4}}}k,{%
\quad }(c$ is the light constant and $k$ is Newton constant) in order to
obtain concordance with the Einstein theory in the locally isotropic limit.

We construct on $\underline{\xi }$ a la--gravitational theory with matter
fields (denoted as ${\varphi }_A$ with $A$ being a general index)
interactions by postulating this Lagrangian density for matter fields
$$
{\underline{{\Cal L}}}^{(m)}={\underline{{\Cal L}}}^{(m)}[{\underline{\hat G}%
}^{\alpha \beta }+{\hat B}^{\alpha \beta };{\frac \delta {\delta u^\gamma }}(%
{\underline{\hat G}}^{\alpha \beta }+{\hat B}^{\alpha \beta });{\varphi }_A;{%
\frac{\delta {\varphi }_A}{\delta u^\tau }}].\tag3.53
$$

Starting from (3.52) and (3.53) the total action of la--gravity on $%
\underline{\xi }$ is written as
$$
{\underline{{\Cal S}}}={(2c{\kappa }_1)}^{-1}\int {\delta }^qu{\underline{%
{\Cal L}}}^{(g)}+c^{-1}\int {\delta }^{(m)}{\underline{{\Cal L}}}^{(m)}.%
\tag3.54
$$
Applying variational procedure on $\underline{\xi },$
 locally adapted to N--connection by using
derivations (1.4) instead of partial derivations, we derive from (3.54) the
la--gra\-vi\-ta\-ti\-on\-al field equations
$$
{\bold {\Theta }}_{\alpha \beta }={{\kappa }_1}
({\underline{{\bold t}}}_{\alpha
\beta }+{\underline{{\bold T}}}_{\alpha \beta })\tag3.55
$$
and matter field equations
$$
{\frac{{\triangle }{\underline{{\Cal L}}}^{(m)}}{\triangle {\varphi }_A}}=0,%
\tag3.56
$$
where $\frac{\triangle }{\triangle {\varphi }_A}$ denotes the variational
derivation.

In (3.55) we have introduced these values: the energy--momentum d--tensor
for la--gravi\-ta\-ti\-on\-al field
$$\multline
{\kappa }_1{\underline{{\bold t}}}_{\alpha \beta }=({\sqrt{|G|}})^{-1}{\frac{%
\triangle {\underline{{\Cal L}}}^{(g)}}{\triangle G^{\alpha \beta }}}%
=K_{\alpha \beta }+{P^\gamma }_{\alpha \beta }P_\gamma -{P^\gamma }_{\alpha
\tau }{P^\tau }_{\beta \gamma }+\\
{\frac 12}{\underline{G}}_{\alpha \beta }{\underline{G}}^{\gamma \tau }({%
P^\phi }_{\gamma \tau }P_\phi -{P^\phi }_{\gamma \epsilon }{P^\epsilon }%
_{\phi \tau }), \endmultline \tag3.57
$$
(where
$$
K_{\alpha \beta }={\underline{D}}_\gamma K_{\alpha \beta }^\gamma ,
$$
$$
2K_{\alpha \beta }^\gamma =-B^{\tau \gamma }{P^\epsilon }_{\tau (\alpha }{%
\underline{G}}_{\beta )\epsilon }-B^{\tau \epsilon }{P^\gamma }_{\epsilon
(\alpha }{\underline{G}}_{\beta )\tau }+
$$
$$
{\underline{G}}^{\gamma \epsilon }h_{\epsilon (\alpha }P_{\beta )}+{%
\underline{G}}^{\gamma \tau }{\underline{G}}^{\epsilon \phi }{P^\varphi }%
_{\phi \tau }{\underline{G}}_{\varphi (\alpha }B_{\beta )\epsilon }+{%
\underline{G}}_{\alpha \beta }B^{\tau \epsilon }{P^\gamma }_{\tau \epsilon
}-B_{\alpha \beta }P^\gamma {\quad }),
$$
$$
2{\bold \Theta }={\underline{D}}^\tau {\underline{D}}_{tau}B_{\alpha \beta }+{%
\underline{G}}_{\alpha \beta }{\underline{D}}^\tau {\underline{D}}^\epsilon
B_{\tau \epsilon }-{\underline{G}}^{\tau \epsilon }{\underline{D}}_\epsilon {%
\underline{D}}_{(\alpha }B_{\beta )\tau }
$$
and the energy--momentum d--tensor of matter
$$
{\underline{{\bold T}}}_{\alpha \beta }=2{\frac{\triangle {\Cal L}^{(m)}}{%
\triangle {\underline{\hat G}}^{\alpha \beta }}}-{\underline{G}}_{\alpha
\beta }{\underline{G}}^{\gamma \delta }{\frac{\triangle {\Cal L}^{(m)}}{%
\triangle {\underline{\hat G}}^{\gamma \delta }}}.\tag3.58
$$
As a consequence of (3.56)--(3.58) we obtain the d--covariant on $\underline{%
\xi }$ conservation laws
$$
{\underline{D}}_\alpha ({\underline{{\bold t}}}^{\alpha \beta }+{\underline{%
{\bold T}}}^{\alpha \beta })=0.\tag3.59
$$
We have postulated the Lagrangian density of matter fields (3.53) in a form
as to treat ${\underline{{\bold t}}}^{\alpha \beta }+{\underline{{\bold T}}}%
^{\alpha \beta }$ as the source in (3.55).

Now we formulate the main results of this section:

\proclaim{\bf Proposition 3.1}
The dynamics of the Einstein la--gravitational fields, modeled as solutions
of equations (3.45) and matter fields on la--space $\xi ,$ can be
equivalently locally modeled on a background la--space $\underline{\xi }$
provided with a trivial d-connection and metric structures having zero
d--tensors of torsion and curvature by field equations (3.55) and (3.56) on
condition that deformation tensor ${P^\alpha }_{\beta \gamma }$ is a
solution of the Cauchy problem posed for basic equations for a chain of
na--maps from $\xi $ to $\underline{\xi }.$
\endproclaim

\proclaim{\bf Proposition 3.2}
Local, d--tensor, conservation laws for Einstein la--gravita\-ti\-on\-al
fields can be written in form (3.59) for la--gravita\-ti\-on\-al (3.57) and
matter (3.58) energy--momentum d--tensors. These laws are d--covariant on
the background space $\underline{\xi }$ and must be completed with invariant
conditions of type (3.12)--(3.15) for every deformation parameters of a chain
of na--maps from $\xi $ to $\underline{\xi }.$
\endproclaim

The above presented considerations consist proofs of both propositions.

We emphasize that nonlocalization of both locally anisotrop\-ic and
isotropic gravitational energy--momentum values on the fundamental (locally
anisotrop\-ic or isotropic) space $\xi $ is a consequence of the absence of
global group automorphisms for generic curved spaces. Considering
gravitational theories from view of multispaces and their mutual maps {\quad
} (directed by the basic geometric structures on $\xi \quad $ such as {\quad
} N--connection, d--connection, d--torsion and d--curvature components, see
coefficients for basic na--equations (3.9)--(3.11)), we can formulate local
d--tensor conservation laws on auxiliary globally automorphic spaces being
related with space $\xi $ by means of chains of na--maps. Finally, we remark
that as a matter of principle we can use d--connection deformations in order
to modelate the la--gravitational interactions with nonvanishing torsion and
nonmetricity. In this case we must introduce a corresponding source in
(3.59) and define generalized conservation laws as in (3.44) \quad (see
similar details for locally isotropic generalizations of the Einstein
gravity in  \cite{Vacaru 1994}).
\vskip80pt
\newpage
\vskip50pt
\head REFERENCES \endhead
\vskip10pt

\ref \by M. Anastasiei  
    \paper Structures spinorielles sur les vari\'et\'es hilbertiennes
    \jour  C. R. Acad. Sc. Paris \vol A284 \yr 1977 \pages 943--945 \endref
\ref \by  P. L. Antonelli and R. Miron (eds) \book Lagrange and Finsler
      Geometry, Applications to Physics and Biology \publ Kluwer Academic 
      Publishers \publaddr  Dordrecht, Boston, London \yr 1996 \endref 
\ref \by  G. S. Asanov, \book Finsler Geometry, Relativity and Gauge
      Theories \publ Reidel \publaddr Boston \yr 1985 \endref
\ref \by G. S. Asanov and S. F. Ponomarenko \book Finsler Bun\-dle
      on Space--Time. Associated Gauge Fields and Con\-nec\-tions
       \publ \c{S}tiin\c{t}a \publaddr Chi\c{s}in\u{a}u 
      \yr 1988 \lang Russian \endref
\ref \by A. Ashtekar, J. D. Romano and S. Ranjet  
     \paper New variables for gravity: Inclusion of matter
     \jour  Phys. Rev. \vol D36 \yr 1989 \pages 2572--2587 \endref
\ref \by  M. F. Atiyah, R. Bott and A. Shapiro 
     \paper Clifford modules 
      \jour Topology \vol 3 \yr 1964 \pages 3--38 \endref
\ref \by  W. Barthel \paper Nichtlineare Zusammenhange und Deren Holonomie 
     Gruppen \jour J. Reine Angew. Math.\vol 212 \yr 1963 \pages 120--149
      \endref
\ref \by A. Bejancu \book Finsler Geometry and Applications
              \publ Ellis Horwood \publaddr Chichester, England \yr 1990 
       \endref
\ref \by  L. Berwald \paper Untersuchung der Kr\"mmung allgemeiner metrischer
      R\"aume auf Grund des in ihnen herrschenden Parallelismus
  \jour Math. Z. \vol 25 \yr 1926 \pages 40--73 \endref
\ref \by R. D. Bishop and R. J. Crittenden \book Geometry of Manifolds
     \publ Academic Press \publaddr London \yr 1964 \endref
\ref \by  E. Cartan \book Les Espaces de Finsler 
      \publ Hermann \publaddr Paris \yr 1934 \endref
\ref \by E. Cartan \book  Les Systems Differentielles Exterieurs et
                Lewrs Application Geometricques
      \publ Herman and Cie Editeur \publaddr  Paris \yr 1945 \endref
\ref \by C. Chevalley \book The Construction and Study of Certain
     Important Algebras \publ Publications of Mathematical Society
     \publaddr Tokyo \yr 1955 \endref
\ref \by P. Finsler \book $\ddot U$ber Kur\-ven und Fl$\ddot a$chen
              in All\-gemeiner R$\ddot a$men. Disser\-ta\-ti\-on
            \publaddr  G$\ddot o$ttingen \yr 1918 \bookinfo re\-print\-ed
              (Birkh$\ddot a$user, Basel, 1951)\endref
\ref \by R. Geroch  
      \jour  J. Math. Phys.\vol 9  \yr 1958 \pages 1739--1750 \endref
\ref \by C. Godbillon \book Elements de Topologie Algebraique
     \publ Herman \publaddr Paris \yr 1971 \endref
\ref \by I. Gottlieb and S. Vacaru \paper A. Moor's tensorial integration
              in generalized Lagrange spaces \inbook  Lagrange and Finsler
              Geometry, Applications to Physics and Biology \eds P. L.
              Antonelli and  Radu Miron  \publ Kluwer Academic Publishers
     \publaddr Dordrecht, Boston, London \yr 1996 \pages 209--216 \endref
\ref \by M. Karoubi \book K--Theory \publ Springer
      \publaddr Berlin \yr 1978 \endref
\ref \by  A. Kawaguchi \paper  On the theory of non--linear connections II. 
     \jour Tensor N. S. \vol 6 \yr 1956 \pages 165--169 \endref
\ref \by J. Kern \paper \paper Lagrange geometry \jour
      Arch. Math.\vol 25 \yr 1974 \pages 438--443 \endref
\ref \by S. Kobayashi and K. Nomizu \book Foundations of Differential
 Geometry \bookinfo Vols 1 and 2 \publ Interscience \publaddr New York
     \yr 1963 \endref
\ref \by S. Lang \book Differential Manifolds
              \publ Reading, Mass, Addison--Wesley \yr 1972 \endref
\ref \by C. P. Luehr and M. Rosenbaum \paper Spinor connection in general 
     relativity \jour  J. Math. Phys. \vol 21 \yr 1974 \pages 1120--1137 
     \endref
\ref \by M. Matsumoto \paper Foundations of Finsler Geometry and Special
       Finsler Spaces \publ Kaisisha \publaddr Shigaken \yr 1986 \endref
\ref \by R. Miron and M. Anastasiei\book Vector Bundles. Lagrange Spaces.
               Application in Relativity \publ Academiei \publaddr Romania
      \yr 1987 \lang Romanian \endref
\ref \by R. Miron and M. Anastasiei \book The Geometry of Lagrange
      Spaces: Theory and Applications  \publ Kluwer Academic Publishers
      \publaddr Dordrecht, Boston, London \yr 1994 \endref
\ref \by R. Miron and Gh. Atanasiu \book  Compendium sur les Espaces
    Lagrange D'ordre \newline Sup\'{e}rieur \bookinfo Seminarul de Mecanic\v a. 
   \publ Universitatea din Timi\c soara. Facultatea de Matematic\v a
   \publaddr Timi\c soara, Rom\^ania \yr 1994  \endref
\ref \by R. Miron and T. Kavaguchi \paper Relativistic geometrical optics
    \jour Int. J. Theor. Phys. \vol 30 \yr 1991  \pages 1521--1543 \endref
\ref \by P. Mocanu \paper Spa\c tii  Par\c tial  Proiective \vol 6, N 3-4
              \jour Academia R. P. Rom\^ane, Bucure\c sti \yr 1955 
       \lang Romannian \endref
\ref \by  A. Mo\'or \jour Acta Math. \vol 86 \yr 1951 \pages 71--80 \endref
\ref \by T. Ono  and Y. Takano \paper Remarks on the spinor gauge field
     theory \jour Tensor N. S.  \vol 52 \yr 1993 \pages 56--60 \endref
\ref \by R. Penrose and W. Rindler \book Spinors and Space-Time, vol. 1,
            Two-Spinor Calculus and Relativistic Fields
    \publ Cambridge University Press \publaddr Cambridge \yr 1984 \endref
\ref \by R. Penrose and W. Rindler \book Spinors and Space-Time, vol. 2,
              Spinor and Twistor Methods in Space-Time Geometry
      \publ Cambridge University Press \publaddr Cambridge \newline 
      \yr 1986 \endref
\ref \by  A. Z. Petrov \paper On selecting modellin of Sun's gravity
      \jour Doklady Academii Nauk SSSR \vol 190 \yr 1970  \pages 305--308
                \lang Russian \endref
\ref \by H. Poincare \book Science and Hypothesis \publ Walter Scott
    \publaddr London \yr 1905 \moreref \publ Dover \publaddr New York
    \yr 1952 \endref
\ref \by H. Poincare \book Oeuvres de Henry Poincare, vol. 9
      \publ Gauther-Villars \publaddr Paris \yr 1954 \endref
\ref \by V. N. Ponomarev, A. O. Barvinsky and Yu. N. Obukhov 
     \book  Geometrodynamical Methods and Gauge Approach to Gravity Theory 
     \publ Energoatomizdat \publaddr Moscow \yr 1985 \lang Russian \endref
\ref \by  D. A. Popov \paper On the theory  of Yang--Mills fields 
  \jour  Theor. Math. Phys. \vol 24 \yr 1975 \pages 347--356 
                 \lang Russian \endref
\ref \by D. A. Popov and L. I. Dikhin \paper Einstein spaces and Yang--Mills
      fields \jour Doklady Academii Nauk SSSR  \vol 225 \yr 1975
      \pages 790--794 \lang Russian \endref
\ref \by  H. Rund \book The Differential Geometry of Finsler Spaces
      \publ Springer--Verlag \publaddr Berlin \yr 1959 \endref
\ref \by S. A. Schouten and D. Struik \book Einfihung in die neucren
               Medoden den Differentiol Geometrie, Bund 1, 2 \yr 1938 \endref
\ref \by N. S. Sinyukov  \book Geodesic Maps of Riemannian Spaces
     \publ Nauka \publaddr Moscow \yr 1979 \lang Russian \endref
\ref \by A. A. Tseytlin \paper Poincare and de Sitter gauge theories of 
     gravity with propagating torsion \jour Phys. Rev. \vol D26
     \yr 1982 \pages 3327--3341 \endref
\ref \by A. Turtoi \book Applications of Algebra and Geometry in Spinors
              Theory \publ Editura Teh\-ni\-c\v a \publaddr Bucure\c sti
      \yr 1989 \lang Romanian \endref
\ref \by S. Vacaru \paper Twistor--gauge interpretation of the
Einstein--Hilbert equations \jour  Vestnik Moscovskogo Universiteta, Fizica i
              Astronomia \vol 28 \yr 1987 \pages 5--12 \lang Russian \endref
\ref \by S. Vacaru \paper Nearly geodesic mappings, twistors and conservation
             laws in gravitational theories \inbook  Contr. Int. Conf.
             " Lobachevski and Modern Geometry ", Part II,  \eds V. Bajanov
             et all \publ University Press \publaddr Kazani, Russia, Tatarstan
      \yr 1992 \pages 64 \endref
\ref \by S. Vacaru \paper Nearly autoparallel maps and conservation laws
        on curved spaces \jour Romanian Journal of Physics \vol 39
       \yr 1994 \pages 37--49 \endref
\ref \by S. Vacaru\paper Spinor structures and nonlinear connections
               in vector bundles, generalized Lagrange and Finsler spaces
        \jour J. Math. Phys.  \vol 37 \yr 1996 \pages 508--523 \endref
\ref \by S. Vacaru and Yu. Goncharenko \paper Yang-Mills fields and gauge
             gravity on generalized Lagrange and Finsler spaces \jour
        Int. J. Theor. Phys. \vol 34 \yr 1995 \pages 1955--1985 \endref
\ref  \by S. Vacaru and S. Ostaf\paper Twistors and nearly autoparallel
      maps \inbook  Coloquium on Differential Geometry, 25-30 July 1994 
     \publ Lajos Kossuth University \publaddr Debrecen, Hungary \yr 1994
      \pages 56 \endref
\ref \by S. Vacaru and S. Ostaf\paper Nearly autoparallel maps of
              Lagrange and Finsler spaces \inbook  Lagrange and Finsler
              Geometry \eds P. L. Antonelli and R. Miron
         \publ Kluwer Academic Publishers \publaddr Dordrecht, Boston, London
         \yr  1996a \pages 241--257 \endref
\ref \by S. Vacaru and S. Ostaf \paper Twistors and nearly autoparallel
   maps \jour Rep. Math. Phys. \vol 37 \yr 1996b \pages 309--328 \moreref
        E-print:\ gr-qc/9602010 \endref
\ref \by G. Vr\v anceanu \book Lec\c tii de Geometrie Differential\v a,
              vol. 2 \publ Ed. Didactic\v a \c si Pedagogic\v a 
     \publaddr Bucure\c sti \yr 1977 \lang Romannian \endref
\ref \by  K. Yano \paper Concircular geometry I--IV \jour  Proc. Imp. Acad.
     Tokyo \vol 16 \yr 1940 \pages 195--200 \moreref \pages 354--360 \moreref
    \pages 442--448 \moreref \pages 505--511 \endref
\ref \by K. Yano and S. I. Ishihara \book Tangent and Cotangent
              Bundles. Differential Geometry 
      \publ Marcel Dekker \publaddr New York \yr 1973 \endref
\bye